\documentclass[10.5pt]{article}
\usepackage{amsfonts}
\usepackage{graphicx}

\usepackage{amsmath}
\usepackage{amssymb}
\usepackage{latexsym}
\usepackage{amsmath, amsfonts,amssymb, amsthm, euscript,makeidx,color,mathrsfs}

\def\sqr#1#2{{\vcenter{\vbox{\hrule height.#2pt
              \hbox{\vrule width.#2pt height#1pt \kern#1pt \vrule width.#2pt}
              \hrule height.#2pt}}}}
\def\signed #1{{\unskip\nobreak\hfil\penalty50
              \hskip2em\hbox{}\nobreak\hfil#1
              \parfillskip=0pt \finalhyphendemerits=0 \par}}
\def\endpf{\signed {$\sqr69$}}

\def\3n{\negthinspace \negthinspace \negthinspace }
\def\2n{\negthinspace \negthinspace }
\def\1n{\negthinspace }

\def\dbB{\mathbb{B}}

\def\dbE{\mathbb{E}}
\def\dbF{\mathbb{F}}

\def\dbH{\mathbb{H}}

\def\dbP{\mathbb{P}}

\def\dbR{\mathbb{R}}
\def\dbS{\mathbb{S}}

\def\dbZ{\mathbb{Z}}

\def\sB{\mathscr{B}}

\def\={\buildrel \triangle \over =}

\def\ds{\displaystyle}

\def\ns{\noalign{\ss}}
\def\a{\alpha}
\def\b{\beta}
\def\g{\gamma}
\def\d{\delta}
\def\e{\varepsilon}
\def\z{\zeta}
\def\k{\kappa}
\def\l{\lambda}
\def\m{\mu}
\def\n{\nu}
\def\si{\sigma}
\def\t{\tau}
\def\f{\varphi}
\def\th{\theta}
\def\o{\omega}

\def\G{\Gamma}
\def\D{\Delta}
\def\Th{\Theta}

\def\F{\Phi}
\def\O{\Omega}

\def\cA{{\cal A}}

\def\cD{{\cal D}}

\def\cF{{\cal F}}

\def\cN{{\cal N}}

\def\cR{{\cal R}}

\def\cU{{\cal U}}

\def\BA{{\bf A}}

\def\BP{{\bf P}}

\def\Bp{{\bf p}}

\def\ss{\smallskip}
\def\ms{\medskip}
\def\bs{\bigskip}
\def\q{\quad}
\def\qq{\qquad}
\def\hb{\hbox}

\def\limsup{\mathop{\overline{\rm lim}}}
\def\liminf{\mathop{\underline{\rm lim}}}

\def\lan{\mathop{\langle}}
\def\ran{\mathop{\rangle}}

\def\esssup{\mathop{\rm esssup}}

\def\wt{\widetilde}

\def\cd{\cdot}
\def\cds{\cdots}

\def\as{\hbox{\rm a.s.}}

\def\tr{\hbox{\rm tr$\,$}}

\def\les{\leqslant}
\def\ges{\geqslant}

\def\({\Big (}
\def\){\Big )}
\def\[{\Big[}
\def\]{\Big]}
\def\bde{\begin{definition}\label}
\def\ede{\end{definition}}
\def\be{\begin{equation}}
\def\bel{\begin{equation}\label}
\def\ee{\end{equation}}
\def\bt{\begin{theorem}\label}
\def\et{\end{theorem}}
\def\bc{\begin{corollary}\label}
\def\ec{\end{corollary}}
\def\bl{\begin{lemma}\label}
\def\el{\end{lemma}}
\def\bp{\begin{proposition}\label}
\def\ep{\end{proposition}}
\def\bas{\begin{assumption}}
\def\eas{\end{assumption}}
\def\br{\begin{remark}\label}
\def\er{\end{remark}}
\def\ba{\begin{array}}
\def\ea{\end{array}}

\def\rf{\eqref}

\def\square#1{\vbox{\hrule\hbox{\vrule height#1%
     \kern#1\vrule}\hrule}}
\def\rectangle#1#2{\vbox{\hrule\hbox{\vrule height#1%
     \kern#2\vrule}\hrule}}

\font\tenbb=msbm10 \font\sevenbb=msbm7 \font\fivebb=msbm5

\newfam\bbfam
\scriptscriptfont\bbfam=\fivebb \textfont\bbfam=\tenbb
\scriptfont\bbfam=\sevenbb

\newtheorem{theorem}{Theorem}[section]
\newtheorem{corollary}[theorem]{Corollary}

\newtheorem{lemma}[theorem]{Lemma}
\newtheorem{proposition}[theorem]{Proposition}

\theoremstyle{definition}
\newtheorem{definition}[theorem]{Definition}
\newtheorem{remark}{Remark}

\makeatletter
   
   \@addtoreset{equation}{section}
\makeatother

\begin{document}

\title{\bf Equilibrium Strategies for Time-Inconsistent\\ Stochastic Switching Systems\footnote{This work is supported in part by NSF Grant
DMS-1406776.}}

\author{Hongwei Mei\thanks{Department of Mathematics, The University of Kansas, Lawrence, KS 66045, USA (hongwei.mei@ku.edu).}\,\, \ \ and \ \ Jiongmin Yong\thanks{Department of Mathematics,
University of Central Florida, Orlando, FL 32816, USA (jiongmin.yong@ucf.edu).}}

\date{}
\maketitle

\begin{abstract}

An optimal control problem is considered for a stochastic differential equation containing a state-dependent regime switching, with a recursive cost functional. Due to the non-exponential discounting in the cost functional, the problem is time-inconsistent in general. Therefore, instead of finding a global optimal control (which is not possible), we look for a time-consistent (approximately) locally optimal equilibrium strategy. Such a strategy can be represented through the solution to a system of partial differential equations, called an equilibrium Hamilton-Jacob-Bellman (HJB, for short) equation which is constructed via a sequence of multi-person differential games. A verification theorem is proved and, under proper conditions, the well-posedness of the equilibrium HJB equation is established as well.
\end{abstract}



\section{Introduction}

Let $(\O,\cF,\dbP)$ be a complete probability space on which a $d$-dimensional standard Brownian motion $\{W(s)\bigm|0\les s\les T\}$ is defined with the natural filtration $\dbF^W=\{\cF^W_s\}_{s\ges0}$. Let $N(ds,d\th)$ be a Poisson random measure on $\dbR$ with the intensity measure $\dbE [N(ds,d\th)]=\pi(d\th)ds$ and the natural filtration $\dbF^N=\{\cF^N_s\}_{s\ges0}$ which is assumed to be independent of the Brownian motion $W(\cd)$. Let
$$\cF_s=\cF^W_s\vee\cF^N_s\vee\cN_0,\qq\dbF=\{\cF_s\}_{s\ges0},$$
where $\cN_0$ is the set of all $\dbP$-null sets. For convenience, we let $\cF=\cF_T$ (if necessary, one may shrink $\cF$ to achieve this). It is clear that $s\mapsto\cF_s$ is right-continuous with left limit.

\ms

Consider the following controlled system of stochastic differential equations (SDEs, for short):
\bel{state}\left\{\2n\ba{ll}
\ds dX(s)=b(s,X(s),\a(s),u(s))ds+\si(s,X(s),\a(s),u(s))dW(s),\qq s\in[\t,T],\\
\ns\ds d\a(s)=\int_{\dbR}\m(X(s),\a(s-),\th)N(ds,d\th),\qq s\in[\t,T],\\
\ns\ds X(\t)=\xi,\qq\a(\t)=\iota,\ea\right.\ee
where $b:[0,T]\times\dbR^n\times M\times U\to\dbR^n$, $\si:[0,T]\times\dbR^n\times M\times U\to\dbR^{n\times d}$ and $\m:\dbR^n\times M\times\dbR\to\dbZ$ ($\dbZ$ is the set of all integers) are given (deterministic) maps with $U\subseteq\dbR^{\bar n}$ being a closed set
(which could be bounded or unbounded, and it could be equal to $\dbR^{\bar n}$), and $M=\{1,2,\cds,m\}$. In the above, the pair $(X(\cd),\a(\cd))$ is called the {\it state process}, valued in $\dbR^n\times M$, with $\a(\cd)$ being called the {\it regime switching process} (under a properly chosen map $\m(\cd\,,\cd\,,\cd)$, see below for more explanations). We call $(\t,\xi,\iota)$ an {\it initial triple} which is taken from the set $\cD$ of all {\it admissible initial triples}, and call $u(\cd)$ a {\it control process} which is taken from a set $\cU[\t,T]$ of all {\it admissible controls}. Both $\cD$ and $\cU[\t,T]$ will be precisely defined shortly.

\ms

Under proper conditions, for any $(\t,\xi,\iota)\in\cD$, and $u(\cd)\in\cU[\t,T]$, state equation \rf{state} admits a unique solution
$$(X(\cd),\a(\cd))\equiv\big(X(\cd\,;\t,\xi,\iota,u(\cd)),\a(\cd\,;\t,\xi,\iota,u(\cd))\big),$$
where $X(\cd)$ has continuous paths and $\a(\cd)$ has c\`adl\`ag (right-continuous with left-limit) paths. We call $(X(\cd),\a(\cd),u(\cd))$ a {\it state-control triple}.

\ms

Let us now briefly look at the main motivation of studying the above controlled system. Consider a financial market in which there are a number of assets trading, and there are some economic factors (such as interest rate, unemployment rate, economy growth rate, etc.) affecting the parameters (such as appreciation rates, volatility, etc.) of the asset price dynamics. It is well-accepted that for different situations of the market (for example, the ``bull market'', the ``bear market'', etc.), the asset price processes and the economic factors should follow different types of dynamics. If we assume that there are $m$ different situations, then we may use the following different SDEs to describe those different situations:
\bel{state-i}dX(s)=b(s,X(s),i,u(s))ds+\si(s,X(s),i,u(s))dW(s),\qq i=1,2,\cds,m.\ee
Components of the vector-valued process $X(\cd)$ could represent the prices of the financial instruments/assets, some economic factors, wealth process, etc. For convenience, we may refer to \rf{state-i} as the dynamics of type-$i$ market, and call $i$ the {\it market index}. As time goes by, the market situation might be changed. This can be described by {\it switching} the index $i\in M\equiv\{1,2,\cds,m\}$ from one to another. It is understood that the change should depend intrinsically on the state $X(\cd)$ in some manner, for example, at an economic crises, the stock market crashes. Right after that, the dynamics of the stock prices could be modeled by an SDE with a small drift. After some time, the economy is slowly recovered, and the dynamics of the stock prices would follow a different SDE for which the drift is larger. Maybe at some time, the economy suddenly grows very fast, and the stock prices should follow a new SDE with a larger drift and a larger volatility. Apparently, if the current market index is $i$, then the index $j\in M\setminus\{i\}$ to which the market will switch is not necessarily certain. To model that, one could assign a probability, say, $p_{ij}$ to the event of switching from $i$ to $j$. It turns out that we could use a regime switching process $\a(\cd)$ to model the above-described switching. Thus, we have the first equation in \rf{state} in which the regime switching process $\a(\cd)$ appears in the drift and diffusion. A natural way of describing the process $\a(\cd)$ is to use a special type of SDE driven by a Poisson process $N(\cd\,,\cd)$. This leads to the second equation in \rf{state}. Interestingly, under certain conditions, the solution $\a(\cd)$ of the second equation in \rf{state} satisfies the following (for $s\in[0,T)$ and $\D s>0$ small):
\bel{transition0}\ba{ll}
\ns\ds\dbP\big(\a(s+\D s)=j\bigm|X(s)=x,\a(s)=i\big)=\left\{\2n\ba{ll}
\ds q_{ij}(x)\D s+o(\D s),\qq\qq\hb{for $j\ne i$,} \\ [2mm]
\ns\ds 1+q_{ii}(x)\D s+o(\D s), \qq\hb{for $j=i$},\ea\right.\ea\ee
for some maps $q_{ij}:\dbR^n\to[0,\infty)$ with $i\ne j$, and $\ds q_{jj}(x)=-\sum_{i\ne j}q_{ij}(x)$.

\ms

We point out that in the above, functions $q_{ij}(x)$ are $x$-dependent.
This is not just for the mathematical generality. This is needed in real applications, as we indicated above, say, the ``bull market'' or the ``bear market'', etc. are determined by the relevant economic factors (some components of the state $X(\cd)$). Because of the state-dependence of the transition probability $q_{ij}(x)$, the resulting problem becomes mathematically delicate and challenging. We will see that some special techniques will be developed to handle the problem. Apparently, one could study the problem with $q_{ij}$ being independent of $x$. Mathematically, the problem becomes much easier. On the other hand, the applications will be
very limited, as one can imagine that the turning point of the market from
``bull'' to ``bear'' (or the other way around) is not determined {\it a priori}.

\ms

There are a number of other interesting problems that will lead to regime-switching systems similar to \rf{state} as well. See \cite{Yin-Zhu2009} and references cited therein for a detailed presentation. The readers are also referred to \cite{Sotomayor-Cadenillas2009,Donnelly2011,Donnelly-Heunis2012,Zhang-Elliott-Siu2012,Heunis2015} for some results relevant to financial applications with regime-switching systems.

\ms

Next, we consider the cost functional to be used to measure the performance of the control $u(\cd)$. In terms of classical stochastic optimal control theory, for our state-control triple $(X(\cd),\a(\cd),u(\cd))$ defined on $[\t,T]$, the cost functional (to be minimized) should look like
\bel{J0}J^0(\t,\xi,\iota;u(\cd))=\dbE_\t\[e^{-\l(T-\t)}h^0(X(T),\a(T))
+\int_\t^Te^{-\l(s-\t)}g^0(s,X(s),\a(s),u(s))ds\],\ee
with $\dbE_\t=\dbE[\,\cd\bigm|\cF_\t]$ being the conditional expectation operator, $\l\ges0$ being called the {\it discount factor}, for some (deterministic) maps $h^0:\dbR^n\times M\to\dbR$ and $g^0:[0,T]\times\dbR^n\times M\times U\to\dbR$. If we introduce the following {\it backward stochastic differential equation} (BSDE, for short) (\cite{Pardoux-Peng1990, Ma-Yong1999, Yong-Zhou1999})
\bel{BSDE0}\left\{\2n\ba{ll}
\ds dY^0(s)=-\l Y^0(s)-g^0\big(s,X(s),\alpha(s),u(s)\big)ds+Z^0(s)dW(s),\qq s\in[\t,T],\\
\ns\ds Y^0(T)=h^0(X(T),\a(T)),\ea\right.\ee
then it admits a unique {\it adapted solution} $(Y^0(\cd),Z^0(\cd))$, and the following holds
\bel{Y0}\ba{ll}
\ns\ds Y^0(t)=e^{-\l(T-t)}h^0(X(T),\a(T))+\int_t^Te^{-\l(s-t)}g^0(s,X(s),\a(s),
u(s))ds\\
\ns\ds\qq\qq\qq\qq\qq\qq-\int_t^Te^{-\l(s-t)}Z^0(s)dW(s),\qq t\in[\t,T].\ea\ee
Taking conditional expectation and setting $t=\t$, we obtain
\bel{Y0=J0}\ba{ll}
\ns\ds Y^0(\t)=\dbE_\t\[e^{-\l(T-\t)}h^0(X(T),\a(T))+\int_\t^Te^{-\l(s-t)}g^0(s,X(s),\a(s),u(s))ds
\]=J^0(\t,\xi,\iota;u(\cd)).\ea\ee
Because of the above, we may call $Y^0(\cd)$ a {\it disutility process} or a {\it cost process}. Inspired by the {\it stochastic differential utility} introduced by Duffie--Epstein \cite{Duffie-Epstein1992a,Duffie-Epstein1992b} (see also \cite{Duffie-Lions1992}), we introduce the following BSDE:
\bel{BSDE1}\left\{\2n\ba{ll}
\ds dY^1(s)=-\l Y^1(s)-g^1\big(s,X(s),\alpha(s),Y^1(s),Z^1(s),\int_\dbR\G^1(s,\th)\pi(d\th),u(s)\big)ds\\
\ns\ds\qq\qq\qq+Z^1(s)dW(s)+\int_\dbR\G^1(s-,\th)\wt N(ds,d\th),\qq s\in[\t,T],\\
\ns\ds Y^1(T)=h^1(X(T),\a(T)).\ea\right.\ee
where
$$\wt N(ds,d\th)=N(ds,d\th)-\pi(d\th)ds,$$
with $N(ds,d\th)$ being the Poisson random measure appears in \rf{state}. Under some mild conditions, BSDE \rf{BSDE1} admits a unique {\it adapted solution} $(Y^1(\cd),Z^1(\cd),\G^1(\cd\,,\cd))$ (\cite{Situ2006, Kruse-Popier2016, Kruse-Popier2017}). Similar to \rf{Y0}, we could define
\bel{Y1}\ba{ll}
\ns\ds J^1(\t,\xi,\iota;u(\cd))=Y^1(\t)=\dbE_\t\[e^{-\l(T-\t)}h^1(X(T),\a(T))\\
\ns\ds\qq\qq\qq\qq\qq\qq+\1n\int_\t^T\2n e^{-\l(s-\t)}g^1\big(s,X(s),\a(s),Y^1(s),Z^1(s),\int_\dbR\G^1(s,\th)
\pi(d\th),u(s)\big)ds\].\ea\ee
Comparing \rf{Y0} with \rf{Y1}, we see that in \rf{Y0}, the right-hand side is independent of $Y^0(\cd)$, whereas, in \rf{Y1}, the right-hand side containing $Y^1(\cd)$ (and $(Z^1(\cd),\G^1(\cd\,,\cd))$). Thus, future utility/disutility affects the current one. We refer to $Y^1(\cd)$, together with $Z^1(\cd)$ and $\G^1(\cd\,,\cd)$, as a {\it recursive cost process}.

\ms

In the above, $e^{-\l(\cd-\t)}$ is called an {\it exponential discount}, which describes people's {\it time-preference}. Studies show that people's time-preference might not necessarily be represented by exponential discounting. Therefore, the recursive cost process $(Y(\cd),Z(\cd),\G(\cd\,,\cd))$ might be the adapted solution to the following family of BSDEs (parameterized by $\t\in[0,T)$):
\bel{BSDE2}\left\{\1n\ba{ll}
\ds dY(\t;s)=-g\big(\t,s,X(s),\alpha(s),Y(\t;s),Z(\t;s),\int_\dbR\G(\t;s,\th)\pi(d\th),u(s)\big)ds\\
\ns\ds\qq\qq\qq+Z(\t;s)dW(s)+\int_\dbR\G(\t;s-,\th)\wt N(ds,d\th),\qq s\in[\t,T],\\
\ns\ds Y(\t;T)=h(\t;X(T),\a(T)).\ea\right.\ee
Clearly, \rf{BSDE1} is a special case of \rf{BSDE2}. Under some mild conditions, for any fixed $\t\in[0,T)$, BSDE \rf{BSDE2} admits a unique {\it adapted solution} $(Y(\t;\cd),Z(\t;\cd),\G(\t;\cd\,,\th))$. We define the following cost functional:
\bel{cost}J(\t,\xi,\iota;u(\cd))=Y(\t;\t).\ee
It is easy to see that the following holds:
\bel{cost2}\ba{ll}
\ns\ds J(\t,\xi;\iota;u(\cd))=Y(\t;\t)=\dbE_\t\[h(\t;X(T),\a(T))\\
\ns\ds\qq\qq\qq\qq\qq\qq+\int_\t^Tg\big(\t,s,X(s),\a(s),Y(\t;s),Z(\t;s),\int_\dbR\G(\t;s,\th)\pi(d\th),u(s)
\big)ds\].\ea\ee
We call the above defined $J(\t,\xi,\iota;u(\cd))$ a {\it recursive cost functional} (with general discounting). Our optimal control problem can be stated as follows.

\ms

\bf Problem (N). \rm For any given $(\t,\xi,\iota)\in\cD$, find a $\bar u(\cd)\in\cU[\t,T]$ such that
\bel{J=inf J}J(\t,\xi,\iota;\bar u(\cd))=\inf_{u(\cd)\in\cU[\t,T]}J(\t,\xi,\iota;u(\cd)).\ee
Any $\bar u(\cd)\in\cU[\t,T]$ satisfying the above is called an {\it open-loop optimal control}, the corresponding state process $(\bar X(\cd),\bar\a(\cd))$ is called an {\it open-loop optimal state process}, and the triple $(\bar X(\cd),\bar\a(\cd),\bar u(\cd))$ is called an {\it open-loop optimal triple}.

\ms

When the switching process $\a(\cd)$ is absent in the state equation and in the cost functional, the problem becomes the one studied by Wei--Yong--Yu (\cite{Wei-Yong-Yu2016}), which is an optimal control of an SDE with recursive cost functional having general discounting. It was shown that such an optimal control problem is time-inconsistent, i.e., an optimal control $\bar u(\cd)$ with optimal state process $\bar X(\cd)$ obtained for an initial pair $(\t,\xi)$ might not remain optimal for a later initial pair $(\t',\bar X(\t'))$ on the optimal state process path (with $\t'>\t$). Therefore, we expect that the above Problem (N) is also time-inconsistent as well. Namely, suppose $\bar u(\cd)\equiv\bar u(\cd\,;\t,\xi,\iota)$ is an open-loop optimal control for Problem (N) with the initial triple $(\t,\xi,\iota)$, and with the open-loop optimal state process
$$(\bar X(\cd),\bar\a(\cd))\equiv\big(X(\cd\,;\t,\xi,\iota,\bar u(\cd)),\a(\cd\,;\t,\xi,\iota,\bar u(\cd))\big).$$
Then it is possible that for some $\t'\in(\t,T)$, the following ay fail:
$$\bar u(\cd\,;\t,\xi,\iota)\Big|_{[\t',T]}=\bar u(\cd\,;\t',\bar X(\t'),\bar\a(\t')),\qq\as$$
This means that the restriction of an open-loop optimal control for the initial triple $(\t,\xi,\iota)$ might not remain optimal for a later initial triple $(\t',\bar X(\t'),\bar\a(\t'))$ (with $\t'\in(\t,T)$) along the optimal path for the triple $(\t,\xi,\iota)$. Due to the possible time-inconsistency of Problem (N), inspired by \cite{Yong2012b,Wei-Yong-Yu2016}, we are going to seek time-consistent equilibrium strategies. The main novelty of the current paper is to explore how the equilibrium strategies, together with the so-called equilibrium Hamilton-Jacobi-Bellman (HJB, for short) equation system will look like when the state-dependent regime switching process $\a(\cd)$ presents. We will adopt the main idea of \cite{Yong2012b, Wei-Yong-Yu2016}. However, due to the appearance of the state-dependent regime switching process $\a(\cd)$ in the state equation and the recursive cost functional, substantially modifications have to be introduced; some of which have their own interests for the switching diffusion processes. We would like to mention that when the current paper has completed, the reference \cite{Wei2017} appeared, which studied the case that the transition probabilities $q_{ij}$ are independent of the state, with non-recursive cost functional. As we have mentioned earlier that such a setting has a serious limitation of applications. Also, approximately local optimality of the equilibrium strategy was not even mentioned there.

\ms

Qualitative study of time-inconsistent problems can be traced back to the works of Hume \cite{Hume1739} and Smith \cite{Smith1759} in the 18th century (see \cite{Palacios-Huerta2003} for a survey). Quantitative investigations started from Strotz \cite{Strotz1955}, followed by the works of Pollak \cite{Pollak1968}, Laibson \cite{Laibson1997}, and so on. In the recent years, time-inconsistent optimal control problems have attracted many authors' attention. Following the basic idea of Pollak \cite{Pollak1968}, Yong introduced the equilibrium HJB equation in \cite{Yong2012b}. See also Yong \cite{Yong2014, Yong2015}, and Ekeland--Lazrak (\cite{Ekeland-Lazrak2010}), Hu--Jin--Zhou (\cite{Hu-Jin-Zhou2012}), Bj\"ork--Murgoci (\cite{Bjork-Murgoci2014}), Bj\"ork--Murguci--Zhou (\cite{Bjork-Murgoci-Zhou2014}), Bj\"ork--Khapko--Murgoci (\cite{Bjork-Khapko-Murgoci2017}), and references cited therein, for relevant works.

\ms

The rest of the paper is organized as follows. In Section 2, we make some preliminaries, mainly on the regime-switching process, the well-posedness of the state equation \rf{state}. Section 3 is devoted to a construction of a family of approximate equilibrium strategies, by means of multi-person differential games. In Section 4, we formally obtain the equilibrium HJB equation and derive the approximate local optimality of the equilibrium strategy. Section 5 is devoted to the well-posedness of the equilibrium HJB equation, for a special case. An example is given in Section 6 to illustrate our problem setting and main result. Finally, some concluding remarks will be collected in Section 7.

\ms

\section{Preliminary}

Let us first recall some standard spaces. As usual, $\dbR^n$ is the $n$-dimensional standard (real) Euclidean space, and $\dbR^{n\times d}$ is the space of all $(n\times d)$ (real) matrices. The set of all $(n\times n)$ symmetric matrices is denoted by $\dbS^n$. The transpose of a matrix $A$ is denoted by $A^\top$. For any Euclidean space $H$ (which could be $\dbR^n$, $\dbR^{n\times d}$, etc.), and $M=\{1,2,\cds,m\}$, let $p\ges1$,
$$\ba{ll}
\ns\ds L^p_\dbF(\t,T;H)=\Big\{\f:[\t,T]\times\O\to H\bigm|\f(\cd)\hb{ is $\dbF$-progressively measurale, }\dbE\int_\t^T|\f(s)|^pds<\infty\Big\},\\
\ns\ds L^p_{\cF_\t}(\O;H)=\Big\{\xi:\O\to H\bigm|\xi\hb{ is $\cF_\t$-measurable, }\dbE|\xi|^p<\infty\Big\},\\
\ns\ds L_\dbF(\t,T;M)=\Big\{\a:[\t,T]\times\O\to M\bigm|\a(\cd)\hb{ is $\dbF$-prograssively measurable}\Big\},\\
\ns\ds L_{\cF_\t}(\O;M)=\Big\{\iota:\O\to M\bigm|\iota\hb{ is $\cF_\t$-measurable}\Big\},\ea$$
Note that since $M$ is a finite set, any random variables/processes taking values in $M$ is always bounded. For any $0\les\t_1<\t_2\les T$, we denote
$$\cD^p[\t_1,\t_2]=\Big\{(\t,\xi,\iota)\bigm|\t\in[\t_1,\t_2],~\xi\in L^p_{\cF_\t}(\O;\dbR^n),~\iota\in L_{\cF_\t}(\O;M)\Big\},\qq p\ges1,$$
and simply denote $\cD^p[0,T]=\cD^p$. Any element in $\cD^p$ is called an
{\it admissible initial triple}.

\ms

Next, we define the set of all {\it admissible control processes} on $[\t_1,\t_2]$ by the following:
$$\ba{ll}
\ns\ds\cU^p[\t_1,\t_2]=\Big\{u:[\t_1,\t_2]\times\O\to U\bigm|u(\cd)\in
L^p_\dbF(t_1,t_2;\dbR^{\bar n})\Big\}.\ea$$
Also, we recall that $\dbE_\t=[\,\cd\,|\,\cF_\t]$.

\ms

\subsection{The regime switching process}

In this subsection, we look at the differential equation in \rf{state} for the switching process $\a(\cd)$.

\ms

Let  $\b_{ij}:\dbR^n\to\dbR$ ($1\les i,j\les m$) be deterministic functions satisfying the following:
\bel{beta}\ba{ll}
\ns\ds-\b_0\les\b_{10}(x)\equiv\b_{11}(x)\les\b_{12}(x)\les\cds\les\b_{1m}(x)\\
\ns\ds\qq=\b_{20}(x)\les\b_{21}(x)\equiv\b_{22}(x)\les\b_{23}(x)\les \cds\les\b_{2m}(x)\\
\ns\ds\qq=\b_{30}(x)\les\cds\\
\ns\ds\qq=\b_{i0}(x)\les \cds\les\b_{i(i-1)}(x)\equiv\b_{ii}(x)\les\b_{i(i+1)}(x)\les\cds\les\b_{im}(x)\\
\ns\ds\qq=\b_{(i+1)0}\les\cds\\
\ns\ds\qq=\b_{m0}(x)\les\b_{m1}(x)\les\cds\les\b_{m(m-1)}(x)\equiv\b_{mm}(x)\les\b_0.
\ea\ee
where $\b_0>0$ is a large fixed constant. Write ${\cR_0}=[-\beta_0,\beta_0]$.  Define
\bel{Delta}
\D_{ij}(x)=\big[\b_{i(j-1)}(x),\b_{ij}(x)\big),\qq1\les i,j\les m,\ee
and for any $\d>0$, denote
\bel{D}\D_{ij}^\d(x)=\bigcup_{|y-x|\les\d}\D_{ij}(y),\qq\D_{ij}^{-\d}(x)
=\bigcap_{|y-x|\les\d}\D_{ij}(y),\qq x\in\dbR^n.\ee
%
Note that for a given $x\in\dbR^n$, for some $i,j\in M$, the set
$\D_{ij}(x)$ could be empty, in particular, it is always true that
$$\D_{ii}=\varnothing,\qq\forall i\in M.$$
Also, it could be true, say, $\D_{14}(x)=[\b_{13}(x),\b_{14}(x))
=\varnothing$ if $\b_{13}(x)=\b_{14}(x)$, which will also lead to $\D_{14}^{-\d}(x)=\varnothing$, for any $\d>0$.

\ms

Now, let us introduce the following hypothesis.

\ms

{\bf(H0)} The L\'evy measure $\pi(\cd)$ is  non-atomic on the Borel $\si$-field $\sB(\dbR)$ of $\dbR$. 
Moreover, the functions $\b_{ij}:\dbR^n\to\dbR$ ($1\les i\les m$, $1\les j\les m$) are uniformly continuous satisfying \rf{beta}. Further, there exists a $K>0$ such that
$$\pi\big(\D_{ij}(x)\big)\les K,\qq\forall x\in\dbR^n,~1\les i,j\les m.$$
and
\bel{D<Kd}\pi\big(\D_{ij}^\d(x)\setminus\D_{ij}(x)\big)+\pi\big(\D_{ij}(x)
\setminus\D_{ij}^{-\d}(x)\big)\les K\d,\qq\forall\d>0,~1\les i,j\les m.\ee

Note that if all the maps $\b_{ij}(\cd)$ ($1\les i\les m,~1\les j\les m$) are uniformly Lipschitz continuous, say, $|\b_{ij}(x)-\b_{ij}(y)|\les L|x-y|$, for some $L>0$, and $\pi(\cd)$ is absolutely continuous with respect to the Lebesgue measure, say, $\pi(ds)=\pi_0(s)ds$, with $0\les\pi_0(s)\les K$ for some measurable function $\pi_0(\cd)$ and constant $K>0$, then (H0) holds. In fact, for such a case,
$$\D^\d_{ij}(x)=\bigcup_{|y-x|\les\d}\D_{ij}(y)
=\[\inf_{|y-x|\les\d}\b_{i(j-1)}(y),\sup_{|y-x|\les\d}\b_{ij}(y)\).$$
$$\D^{-\d}_{ij}(x)=\bigcap_{|y-x|\les\d}\D_{ij}(y)=\[\sup_{|y-x|\les\d}
\b_{i(j-1)}(y),\inf_{|y-x|\les\d}\b_{ij}(y)\).$$
Hence, one has
$$\pi\big(\D^\d_{ij}(x)\setminus\D_{ij}(x)\big)\les KL\(|\b_{i,{j-1}}(x)-\inf_{|y-x|\les\d}\b_{i,{j-1}}(y)|+|\sup_{|y-x|\les\d}\b_{ij}(y)-
\b_{ij}(x)|\)\les2KL\d.$$
Likewise,
$$\pi\big(\D_{ij}(x)\setminus\D^{-\d}_{ij}(x)\big)\les2KL\d.$$
The above observation shows that \rf{D<Kd} is a very mild and reasonable condition. We keep in mind that under (H0), due to the uniform continuity of the functions $\b_{ij}(\cd)$, the length
$$|\D_{ij}(x)|\equiv|\b_{ij}(x)-\b_{i(j-1)}(x)|$$
is also uniformly continuous. Therefore, if $\D_{ij}(x)=\varnothing$,
then $|\D^\d_{ij}(x)|$ will also be small.

\ms

Next, we define
\bel{q(ij)}\left\{\1n\ba{ll}
\ds q_{ij}(x)=\pi\big(\D_{ij}(x)\big),\qq\qq\qq\qq1\les i,j\les m,~i\ne j,\\
\ns\ds q_{jj}(x)=-\sum_{i\ne j}q_{ij}(x),\qq\qq\qq\qq1\les i\les m,\ea\right.\qq x\in\dbR^n,\ee
%
Now, we take the map $\m(\cd\,,\cd\,,\cd)$ in the equation for $\a(\cd)$ as follows:
\bel{m=sum}\m(x,i,\th)=\sum_{j=1}^m(j-i)I_{\D_{ij}(x)}(\th),\qq(x,i,\th)\in\dbR^n\times M\times\dbR.\ee
Then the second equation in \rf{state} reads
\bel{da}d\a(s)=\sum_{j=1}^m\int_\dbR\big[j-\a(s-)\big]I_{\D_{\a(s-)j}(X(s))}(\th)N(ds,d\th)
=\sum_{j=1}^m\big[j-\a(s-)\big]N\big(ds,\D_{\a(s-)j}(X(s))\big).\ee
Thus, the equation for $\a(\cd)$ over $[\t,T]$ with initial condition $\a(\t)=i$ can also be written as
\bel{a(t)}\a(t)=i+\sum_{j=1}^m\int_{(\t,t]}\big[j-\a(s-)\big]N(ds,\D_{\a(s-)j}(X(s))
\big),\qq t\in[\t,T].\ee
The well-posedness of the state equation \rf{state} will be carried out in the next subsection. Let us assume that for a given $\dbF$-adapted continuous process $X(\cd)$, and initial condition $\a(\t)=i$, the above admits a unique solution $\a(\cd)$ on $[\t,T]$. We have the following result.

\bp{a} \sl Let {\rm(H0)} hold. Suppose $X(\cd)$ is an $\dbF$-adapted
process with $X(s)=x$ satisfying
\bel{|X-X|}\dbE\[\sup_{s\les\t\les s+\D s}|X(\t)-X(s)|^p\]\les K|\D s|^{p\over2},\ee
for some constants $p>2$ and $K>0$. Then the solution $\a(\cd)$ to the equation \rf{da} satisfies
\bel{P(a)}\dbP\big(\a(s+\D s)=j\bigm|\a(s)=i,X(s)=x\big)=q_{ij}(x)\D s+o(\D s),~~\hb{for~} j\ne i.\ee

\ep

\it Proof. \rm Let $\k>0$ be small so that the second condition in \rf{beta} holds, and let $\d\in(0,\k]$. Denote
$$\ba{ll}
\ns\ds A_0=\(N\big((s,s+\D s],\cR_0\big)=0\),\q A_1=\(N\big((s,s+\D s],\cR_0)=1\),\q A_2=\(N\big((s,s+\D s],\cR_0\big)\ges2\).\ea$$
%
Then one has (noting $\pi(\dbR)=1$)
$$\dbP(A_0)=e^{-\D s},\q\dbP(A_1)=e^{-\D s}\D s,\q\dbP(A_2)=\sum_{k=2}^\infty{(\D s)^k\over k!}e^{-\D s}\les(\D s)^2.$$
Now, let $j\in M\setminus\{i\}$ and $\a(\cd)$ be the solution to \rf{da} with $\a(s)=i$. Then
$$\ba{ll}
\ns\ds\dbP\(\big(\a(s+\D s)=j\bigm|X(s)=x,\a(s)=i\big)\cap A_0\)=\dbP(\varnothing)=0,\\
\ns\ds\dbP\(\big(\a(s+\D s)=j\bigm|X(s)=x,\a(s)=i\big)\cap A_2\)\les\dbP(A_2)\les(\D s)^2.\ea$$
Hence,
$$\ba{ll}
\ns\ds\dbP\(\a(s+\D s)=j\bigm|X(s)=x,\a(s)=i\)=\dbP\(\big(\a(s+\D s)=j,X(s)=x,\a(s)=i\big)\cap A_1\)+o(\D s).\ea$$
On $A_1\cap(\a(s+\D s)=j\bigm|X(s)=x,\a(s)=i)$, there exists a stopping time $\bar s\in(s,s+\D s]$ so that the underline L\'evy process has no jumps on $[s,\bar s)\cup(\bar s,s+\D s]$, and only has a jump at $\bar s$. Then (still on this set)
$$\ba{ll}
\ns\ds j-i=\int_{(s,s+\D s]}\sum_{k=1}^m(k-i)N\big(ds,\D_{ik}(X(t))\big)
=\sum_{k=1}^m(k-i)N\big((s,s+\D s],\D_{ik}(X(\bar s))\big)\\
\ns\ds\qq=\sum_{k=1}^mkN\big((s,s+\D s],\D_{ik}(X(\bar s))\big)-i.\ea$$
This implies that
$$\ba{ll}
\ns\ds A_1\cap\big(\a(s+\D s)=j\bigm|\a(s)=i,X(s)=x\big)=A_1\cap\(N\big((s,s+\D s],\D_{ij}
(X(\bar s))\big)=1\)\\
\ns\ds=\[A_1\cap\(N\big((s,s+\D s],\D_{ij}(X(\bar s))\big)=1,|X(\bar s)-x|<\d\)\]\\
\ns\ds\qq\qq\qq\bigcup\[A_1\cap
\(N\big((s,s+\D s],\D_{ij}(X(\bar s))\big)=1,|X(\bar s)-x|\ges\d\)\].\ea$$
Note that
$$\ba{ll}
\ns\ds\dbP\(A_1\cap
\big(N\big((s,s+\D s],\D_{ij}(X(\bar s))\big)=1,|X(\bar s)-x|\ges\d\big)\)\les\dbP\(|X(\bar s)-x|\ges\d\)\les{K\over\d^p}(\D s)^{p\over2}.\ea$$
Also, since for $|X(\bar s)-x|<\d$, $\D_{ij}(X(\bar s))\subseteq\D_{ij}^\d(x)$, we have
$$N\big((s,s+\D s],\D_{ij}^\d(x)\big)=N\big((s,s+\D s],\D_{ij}(X(\bar s))\big)+N\big((s,s+\D s],\D_{ij}^\d(x)\setminus\D_{ij}(X(\bar s))\big).$$
Hence, one has the following:
$$\ba{ll}
\ns\ds A_1\cap\(N\big((s,s+\D s],\D_{ij}\big(X(\bar s)\big)=1,|X(\bar s)-x|<\d\)\\
\ns\ds\subseteq\[A_1\cap\(N\big((s,s+\D s],\D_{ij}^\d(x)\big)=1,|X(\bar s)-x|<\d\)\]\bigcup\[A_1\cap\(N\big((s,s+\D s],\D_{ij}^\d(x)\big)\ges2\)\]\\
\ns\ds=\[A_1\cap\(N((s,s+\D s],\D_{ij}(x))=1\)\]\bigcup
\[A_1\cap\(N\big((s,s+\D s],\D_{ij}^\d(x)\setminus\D_{ij}(x)\big)=1\)\],\ea$$
noting that $A_1\cap\(N\big((s,s+\D s],\D_{ij}^\d(x)\big)\ges2\)
=\varnothing$. On the other hand, by the choice of $\k>0$ (see \rf{beta}), for $|x-y|<\d<\k$, $\D_{ij}(x)\cap\D_{ij}(y)\ne\varnothing$. Then,
making use of (H0), we have
$$\dbP\(N\big((s,s+\D s],\D_{ij}^\d(x)\setminus\D_{ij}(x)\big)=1\)=\pi\big(\D_{ij}^\d(x)\setminus\D_{ij}(x)\big)\D s\les K\d\D s.$$
Now, we choose $0<\e<{p-2\over2}$ and $(\D s)^{p-2-2\e\over p}\les\d<\k$. Then
$$\ba{ll}
\ns\ds\dbP\(A_1\cap\big(N\big([s,s+\D s],\D_{ij}(X(\bar s))\big)=1,
|X(\bar s)-x|\ges\d\big)\)\les\dbP\big(|X(\bar s)-x|\ges\d\big)\les K(\D s)^{1+\e}=o(\D s),\\
\ns\ds\dbP\(N([s,s+\D s],\D_{ij}^\d(x)\setminus\D_{ij}(x))=1\)
\les\pi\big(\D_{ij}^\d(x)\setminus\D_{ij}(x)\big)\D s\les K(\D s)^{1+{p-2-2\e\over p}}=o(\D s).\ea$$
Hence,
$$\ba{ll}
\ns\ds\dbP\(A_1\cap\big(\a(s+\D s)=j\bigm|\a(s)=i,X(s)=x\big)\)
=\dbP\(N\big([s,s+\D s],\D_{ij}(x)\big)=1\)+o(\D s)\\
\ns\ds=\pi(\D_{ij}(x))(\D s)e^{-\pi(\D_{ij}(x))\D s}+o(\D s)=q_{ij}(x)\D s+o(\D s).\ea$$
This completes the proof. \endpf

\ms

According to the above, we see that for any $\dbF$-adapted continuous process $X(\cd)$ satisfying \rf{|X-X|}, the solution $\a(\cd)$ of the second equation in \rf{state} satisfies the following:
\bel{transition2}\ba{ll}
\ds\dbP\big(\a(s+\D s)=j\bigm|X(s),\a(s)\hb{ given}\big)=\left\{\2n\ba{ll}
\ds q_{\a(s)j}\big(X(s)\big)\D s+o(\D s),\qq\qq\q\hb{for $j\ne\a(s)$;} \\ [2mm]
\ns\ds 1+q_{\a(s)\a(s)}\big(X(s)\big)\D s+o(\D s), \qq\hb{for $j=\a(s)$}.\ea\right.\ea\ee
In what follows, we call $Q(x)=[q_{ij}(x)]_{m\times m}$ the {\it probability transition matrix} of process $\a(\cd)$. If $Q(x)$ is independent of $x$, $\a(\cd)$ is a Makov process by itself. See \cite{Yin-Zhu2009} for some detailed discussion on $\a(\cd)$.

\ms

\subsection{Well-posedness of the state equation and the recursive cost functional}

In this subsection, we present the well-posedness of our state equation \rf{state} and some relevant results. Due to the presence of the regime switching process, some results have their own interest. Let us first introduce the following hypothesis.

\ms

{\bf (H1)} The maps $b:[0,T]\times\dbR^n\times M\times U\to\dbR^n$ and $\si:[0,T]\times \dbR^n\times M\times U\to\dbR^{n\times d}$ are continuous and there exists a constant $L>0$ and a fixed $u_0\in U$ such that for any $(s,i)\in [0,T]\times M\times U$, $(x_1,u_1),(x_2,u_2)\in\dbR^n\times U$,
\bel{|b-b|+|si-si|}\left\{\2n\ba{ll}
\ns\ds|b(s,x_1,i,u_1)-b(s,x_2,i,u_2)|+|\si(s,x_1,i,u_1)-\si(s,x_2,i,u_2)|\les L\big(|x_1-x_2|+|u_1-u_2|\big),\\
\ns\ds|b(s,0,i,u_0)|+|\si(s,0,i,u_0)|\les L.\ea\right.\ee

We will see that for the well-posedness of the state equation, the Lipschitz continuity of the maps $x\mapsto(b(s,x,i,u),\si(s,x,i,u))$ will be enough. The Lipschitz continuity of the map $u\mapsto(b(s,x,i,u),\si(s,x,i,u))$ will be used in proving the approximate local optimality of the equilibrium strategy (in Section 4).

\ms

The following is the main result of this subsection.

\bt{well-posedness} \sl Let {\rm(H0)--(H1)} hold. Then for any initial
triple $(\t,\xi,\iota)\in\cD^r$ and $u(\cd)\in\cU^p[\t,T]$ with $p>2$,
state equation \rf{state} admits a unique solution $(X(\cd),\a(\cd))$. Moreover, the following estimates hold:
\bel{|X|}\dbE_\t\[\sup_{\t\les s\les T}|X(s)|^p\]\les K\(1+|\xi|^p
+\dbE_\t\int_\t^T|u(s)|^pdr\),\ee
\bel{|X-X|}\dbE_\t\[\sup_{\t\les r\les s}|X(r)-\xi|^p\]\les K(s-\t)^{{p\over2}-1}\[(s-\t)(1+|\xi|^p)+\int_\t^s|u(r)|^pdr\].\ee
Further, if $(X_1(\cd),\a_1(\cd))$ and $(X_2(\cd),\a_2(\cd))$ are solutions of the state equation \rf{state} corresponding to the initial triples
$(\t,\xi_1,\iota),(\t,\xi_2,\iota)\in\cD^p$ (only $\xi_1$ and $\xi_2$
could possibly be different), then
\bel{P(Ac)+|X-X|}\(\dbE_\t[I_{A_T^c}]\)^2+\dbE_\t\[\sup_{\t\les s\les T}|X_1(s)-X_2(s)|^2I_{A_s}\]\les K|\xi_1-\xi_2|^2,\ee
where
$$A_s=\{\o\in\O\bigm|\a_1(t,\o)=\a_2(t,\o),\q\t\les t\les s\},\qq s\in[\t,T].$$

\et

We would like to point out that estimates \rf{|X|} and \rf{|X-X|} are standard for SDEs. Whereas, estimate \rf{P(Ac)+|X-X|} seems to be new for SDEs with regime switchings.

\ms

\it Proof. \rm The proof is lengthy and is split into several steps.

\ms

\it Step 1. A priori estimate for the solutions. \rm

\ms

Suppose $(X(\cd),\a(\cd))$ is a solution to the state equation
\rf{state} corresponding the initial triple $(\t,\xi,\iota)\in\cD^p$
and the control $u(\cd)\in\cU^p[\t,T]$. Then we see that
$$\ba{ll}
\ns\ds\dbE_\t\[\sup_{\t\les t\les s}|X(t)|^p\]\les K\[|\xi|^p+\dbE_\t\(\int_\t^s|b(t,X(t),\a(t),u(t))|dt\)|^p
+\dbE_\t\Big|\int_\t^s \si(t,X(t),\a(t),u(t))dW(t)\Big|^p\]\\
\ns\ds\qq\qq\qq\qq\les K\[|\xi|^p+\dbE_\t\(\int_\t^s
L\big(1+|X(t)|+|u(t)|\big)dt\)^p+\dbE_\t\(\int_\t^sL^2\big(1+|X(t)|
+|u(t)|\big)^2dt\)^{p\over2}\]\\
\ns\ds\qq\qq\qq\qq\les
K\[1+|\xi|^p+\int_\t^s\dbE_\t\(\sup_{
\t\les r\les t}|X(r)|^p\)dt+\dbE_\t\int_\t^s|u(t)|^pdt\].\ea$$
By Gronwall's inequality, we have
$$\dbE_\t\[\sup_{\t\les t\les s}|X(t)|^p\]\les K\[1+|\xi|^p+
\dbE\int_\t^s|u(t)|^pdt\],\qq s\in[\t,T].$$
Thus, \rf{|X|} follows. Also,
\bel{E|X-X|}\ba{ll}
\ns\ds\dbE_\t\[\sup_{\t\les t\les s}|X(t)-\xi|^p\]\\
\ns\ds\les K\[\dbE_\t\(\int_\t^s |b(t,X(t),\a(t),u(t))|dt\)^p+\dbE_\t\(\int_\t^s |\si(t,X(t),\a(t),u(t))|^2dt\)^{p\over2}\]\\
\ns\ds\les K\[\dbE_\t\(\int_\t^sL\big(1+|X(t)|+|u(t)|\big)dt\)^p
+\dbE_\t\(\int_\t^sL^2\big(1+|X(t)|+|u(t)|\big)^2dt\)^{p\over2}\]\\
\ns\ds\les K\big[(s-\t)^{p-1}+(s-\t)^{p-2\over2}\big]
\dbE_\t\int_\t^s\big(1+|X(t)|^p+|u(t)|^p\big)dt\\ [1mm]
\ns\ds\les K(s-\t)^{p-2\over2}\[(s-\t)(1+|\xi|^p)
+\dbE_\t\int_\t^s|u(t)|^pdt\],\qq s\in[\t,T].\ea\ee
Thus, \rf{|X-X|} follows.

\ms

\it Step 2. Uniqueness and continuous dependence with respect to the initial state. \rm

\ms

Let $(X_1(\cd),\a_1(\cd))$ and $(X_2(\cd),\a_2(\cd))$ be the solutions to the state equation with the initial triple $(\t,\xi_1,\iota)$, $ (\t,\xi_2,\iota)\in\cD^p$, and control $u(\cd)\in\cU^p[\t,T]$. Note that
$$\a_1(\t)=\a_2(\t)=\iota.$$
Denote
$$\n(s_1,s_2)=\dbE\int_{s_1}^{s_2}|u(t)|^pdt,\qq\t\les s_1<s_2\les T.$$

\ms

We discretize $[\t,T]$ into $\ell$ equal length intervals, each of
which is of length $\eta={T-\t\over\ell}$. Let $t_k=\t+k\eta$ for $k=0,\cds,\ell$. Let $i,j\in M$ with $j\ne i$. Define
$$\ba{ll}
\ns\ds B^i_{k+1}=\big\{\a_1(t_k)=\a_2(t_k)=i,~N((t_k,t_{k+1}],\cR_0)=1\big\},\\
\ns\ds C^i_{k+1}=\big\{\a_1(t_k)=\a_2(t_k)=i,~N((t_k,t_{k+1}],\cR_0)\ges2\big\},\qq k\ges1.\ea$$
We want to estimate the (conditional) probability of the set
$$B^i_{k+1}\cap\big\{\a_1(t_{k+1})=j,\a_2(t_{k+1})\ne j\big\}.$$
To this end, we denote
$$\ba{ll}
\ns\ds A^{ij}_1(X_1(\cd))=B^i_{k+1}\bigcap\Big\{\int_{(t_k,t_{k+1}]}N\big(ds,\D_{ij}
(X_1(s))\big)=1\Big\},\\
\ns\ds A^{ij}_1(X_1(t_k))=B^i_{k+1}\bigcap\Big\{N\big((t_k,t_{k+1}],\D_{ij}
(X_1(t_k))\big)=1\Big\},\\
\ns\ds A^{ij}_0(X_1(t_k))=B^i_{k+1}\bigcap\Big\{N\big((t_k,t_{k+1}],\D_{ij}
(X_1(t_k))\big)=0\Big\},\\
\ns\ds A^{ij}_1(X_2(t_k))=B^i_{k+1}\bigcap\Big\{N\big((t_k,t_{k+1}],\D_{ij}
(X_2(t_k))\big)=1\Big\},\\
\ns\ds A^{ij}_0(X_2(t_k))=B^i_{k+1}\bigcap\Big\{N\big((t_k,t_{k+1}],\D_{ij}
(X_2(t_k))\big)=0\Big\},\\
\ns\ds A^{ij}_0(X_2(\cd))=B^i_{k+1}\bigcap\Big\{\int_{(t_k,t_{k+1}]}N\big(ds,\D_{ij}
(X_2(s))\big)=0\Big\}.\ea$$
Note that the above defined sets are all subsets of $B^i_{k+1}$, and
$$B^i_{k+1}=A^{ij}_0(X_1(t_k))\bigcup A^{ij}_1(X_1(t_k))=A^{ij}_0(X_2(t_k))\bigcup A^{ij}_1(X_2(t_k)).$$
Also, we denote
$$G_{i'}(\d)=\Big\{\sup_{s\in[t_k,t_{k+1}]}|X_{i'}(s)-X_{i'}(t_k)|\ges\d\Big\},\q i'=1,2,\q\d>0.$$
Then
$$\ba{ll}
\ns\ds B^i_{k+1}\cap\big\{\a_1(t_{k+1})=j,\a_2(t_{k+1})\ne j\big\}= A^{ij}_1(X_1(\cd))\cap A^{ij}_0(X_2(\cd))\\
\ns\ds=\[A_1^{ij}(X_1(\cd))\cap A^{ij}_0(X_2(\cd))\]\bigcap\[A_0^{ij}(X_1(t_k))\cup A_1^{ij}(X_1(t_k))\]\\
\ns\ds=\[A_1^{ij}(X_1(\cd))\cap A^{ij}_0(X_2(\cd))\cap A_0^{ij}(X_1(t_k))\]\bigcup\[A_1^{ij}(X_1(\cd))\cap A^{ij}_0(X_2(\cd))\cap A_1^{ij}(X_1(t_k))\]\\
\ns\ds\subseteq\[A_1^{ij}(X_1(\cd))\cap A_0^{ij}(X_1(t_k))\]\bigcup \[A_1^{ij}(X_1(t_k))\cap A^{ij}_0(X_2(\cd))\]\\
\ns\ds=\[A_1^{ij}(X_1(\cd))\cap A_0^{ij}(X_1(t_k))\]\bigcup \Big\{\[A_1^{ij}(X_1(t_k))\cap A^{ij}_0(X_2(\cd))\]\bigcap\[A_0^{ij}(X_2(t_k))\cap A^{ij}_1(X_2(t_k))\]\Big\}\\
\ns\ds\subseteq\[A_1^{ij}(X_1(\cd))\cap A_0^{ij}(X_1(t_k))\]\bigcup \[A_1^{ij}(X_1(t_k))\cap A^{ij}_0(X_2(t_k))\]\bigcup\[A_1^{ij}(X_2(t_k))\cap A^{ij}_0(X_2(\cd))\]\\
\ns\ds\subseteq\(B^i_{k+1}\cap\Big\{N\big((t_k,t_{k+1}],\D^\d_{ij}
(X_1(t_k))\setminus\D_{ij}(X_1(t_k)\big)=1\Big\}\cap G_1(\d)^c\)\bigcup
G_1(\d)\\
\ns\ds\q\bigcup\(B^i_{k+1}\cap\Big\{N\big((t_k,t_{k+1}],\D_{ij}
\big(X_1(t_k))\setminus\D_{ij}(X_2(t_k))\big)=1\Big\}\)\\
\ns\ds\q\bigcup\(B^i_{k+1}\cap\Big\{N\big((t_k,t_{k+1}],\D_{ij}
\big(X_2(t_k))
\setminus\D_{ij}^{-\d}(X_2(t_k)\big)=1\Big\}\cap G_2(\d)^c\)\bigcup G_2(\d).\ea$$
We now estimate the conditional probabilities of the sets on the right-hand side of the above one-by-one. By Chebyshev's inequality, one has
$$\ba{ll}
\ns\ds\dbE_{t_k}\big[I_{G_1(\d)\cup G_2(\d)}\big]\les\dbP_{t_k}\(\sup_{s\in[t_k,t_{k+1}]}|X_1(s)-X_1(t_k)|\ges\d\)
+\dbP_{t_k}\(\sup_{s\in[t_k,t_{k+1}]}|X_2(s)-X_2(t_k)|\ges\d\)\\
%
%
\ns\ds\les{1\over\d^p}\[\dbE_{t_k}\(\sup_{s\in[t_k,t_{k+1}]}|X_1(s)-X_1(t_k)|^p\)
+\dbE_{t_k}\(\sup_{s\in[t_k,t_{k+1}]}|X_2(s)-X_2(t_k)|^p\)\]\\
\ns\ds\les{K\over\d^p}(t_{k+1}-t_k)^{p-2\over2}\[(t_{k+1}-t_k)(1+|X_1(t_k)|^p
+|X_2(t_k)|^p\big)+\dbE_{t_k}\int_{t_k}^{t_{k+1}}|u(t)|^pdt\]\\
\ns\ds=K\eta^{p-2\over2}\d^{-p}\[\eta\big(1+|X_1(t_k)|^p
+|X_2(t_k)|^p\big)+\dbE_{t_k}\big[\n(t_k,t_{k+1})\big]\].\ea$$
Hereafter, $\dbP_t[A]\equiv\dbE_t[I_A]$. Next,
$$\ba{ll}
\ns\ds\dbP_{t_k}\(B^i_{k+1}\cap\Big\{N\big((t_k,t_{k+1}],\D_{ij}
(X_1(t_k))\setminus\D_{ij}(X_2(t_k)\big)=1\Big\}\)\\
\ns\ds\les K\pi\big(\D_{ij}(X_1(t_k))\setminus\D_{ij}(X_2(t_k))\big)(t_{k+1}-t_k)
\les K\eta|X_1(t_k)-X_2(t_k)|,\ea$$
and by (H0),
$$\ba{ll}
\ns\ds\dbP_{t_k}\(B^i_{k+1}\cap\Big\{N\big((t_k,t_{k+1}],\D^\d_{ij}
(X_1(t_k))\setminus\D_{ij}(X_1(t_k)\big)=1\Big\}\cap G_2(\d)^c\)\\
\ns\ds\les K\pi\big(\D_{ij}^\d(X_1(t_k))\setminus\D_{ij}(X_1(t_k))\big)
(t_{k+1}-t_k)\les K\d\eta,\ea$$
$$\ba{ll}
\ns\ds\dbP_{t_k}\(B^i_{k+1}\cap\Big\{N\big((t_k,t_{k+1}],\D_{ij}
(X_2(t_k))\setminus\D^{-\d}_{ij}(X_2(t_k)\big)=1\Big\}\cap G_2(\d)^c\)\\
\ns\ds\les K\pi\big(\D_{ij}(X_2(t_k))\setminus\D_{ij}^{-\d}(X_2(t_k))\big)(t_{k+1}-t_k)
\les K\d\eta.\ea$$
Hence,
\bel{2.18}\ba{ll}
\ns\ds\dbP_{t_k}\(B^i_{k+1}\cap\big\{\a_1(t_{k+1})=j,\a_2(t_{k+1})\ne j\big\}\)\\
\ns\ds\les K\eta\[\(1\1n+\1n|X_1(t_k)|^p\2n+\1n|X_2(t_k)|^p\)\eta^{p-2\over2}\d^{-p}\1n
+\1n\d\1n+\1n|X_1(t_k)\1n-\1n X_2(t_k)|\1n+\1n\eta^{r-4\over2}\d^{-p}\dbE_{t_k}\2n\int_{t_k}^{t_{k+1}}
\2n|u(t)|^pdt\],\ea\ee
We note that
$$A_k^c\subseteq A_{k+1}^c\subseteq A_k^c\bigcup C_{k+1}\bigcup\(\bigcup_{i,j=1}^m B_{k+1}^i\cap\Big\{\a_1(t_{k+1})=j,\a_2(t_{k+1})
\ne j\Big\}\cap A_k\).$$
Hence,
$$\ba{ll}
\ns\ds\dbP_{t_k}(A_k^c)\les\dbP_{t_k}(A_{k+1}^c)\les\dbP_{t_k}(A_k^c)+\dbP_{t_k}(C_{k+1})
+\sum_{i,j=1}^m\dbP_{t_k}\(B_{k+1}^i\cap\Big\{\a_1(t_{k+1})=j,\a_2(t_{k+1})
\ne j\Big\}\cap A_k\)\\
\ns\ds\les I_{A_k^c}+K\eta^2+K\eta\[\(1+|X_1(t_k)|^p
+|X_2(t_k)|^p\)\eta^{p-2\over2}\d^{-p}+\d\\
\ns\ds\qq\qq\qq\qq\qq+|X_1(t_k)-X_2(t_k)|+\eta^{p-4\over2}\d^{-p}
\dbE_{t_k}\int_{t_k}^{t_{k+1}}|u(t)|^pdt\]I_{A_k}\ea$$
Applying $\dbE_\t$ on both sides of the above, making use of the a priori estimate \rf{|X|}, we have
\bel{probaA-1}\ba{ll}
\ns\ds0\les\dbP_\t(A^c_{k+1})-\dbP_\t(A^c_k)\les  K\eta^2+K\eta\Big\{\[\eta^{p-2\over2}\d^{-p}\dbE_\t\(1+|X_1(t_k)|^p
+|X_2(t_k)|^p)|\)+\d\\
\ns\ds\qq\qq\qq\qq\qq\qq+\dbE_\t\(|X_1(t_k)-X_2(t_k)|I_{A_k}\)
+\eta^{p-4\over2}\d^{-p}\dbE_\t\int_{t_k}^{t_{k+1}}|u(t)|^pdt\]\Big\}\\
\ns\ds\q\les K\eta^2+K\eta\Big\{\eta^{p-2\over2}\d^{-p}
+\d+\dbE\(|X_1(t_k)-X_2(t_k)|I_{A_k}\)+\eta^{p-4\over2}\d^{-p}
\dbE\big[\n(t_k,t_{k+1})\big]\Big\}.\ea\ee
Note that
$$\n(t_0,t_1)+\n(t_1,t_2)+\cds+\n(t_k,t_{k+1})=\n(t_0,t_{k+1}).$$
Thus
\bel{P(Ac)}\dbP(A^c_{k+1})\les K\[k\eta^2+k\eta\d+\eta\sum_{j=1}^k
\dbE\(|X_1(t_j)-X_2(t_j)|I_{A_{t_j}}\)+\eta^{p-2\over2}\d^{-p}
\dbE\big[\n(t_0,t_{k+1})\big]\]\ee
On the other hand, let
$$\bar\t=T\land\inf\{s\in[\t,T]\bigm|\a_1(s)\ne\a_2(s)\}.$$
Then
\bel{stability-3}\ba{ll}
\ns\ds\dbE_\t\[\sup_{\t\les t\les s}|X_1(t\land\bar\t)
-X_2(t\land\bar\t)|\]^2\\
\ns\ds\les K\Big\{|\xi_1-\xi_2|^2+\dbE_\t\(\int_\t^{s\land\bar\t}
|b(t',X_1(t'),\a_1(t'),u(t'))-b(t',X_2(t'),\a_2(t'),u(t'))|dt'\)^2\\
\ns\ds\qq+\dbE_\t\sup_{\t\les t\les s}\Big|\int_\t^{t\land\bar\t}
[\si(t',X_1(t'),\a_1(t'),u(t'))-\si(t',X_2(t'),\a_2(t'),u(t'))]dW(t')
\Big|^2\Big\}\\
\ns\ds\les K\Big\{|\xi_1-\xi_2|^2+\int_\t^s\dbE_\t\[\sup_{\t\les t'\les t}|X_1(t'\land\bar\t)
-X_2(t'\land\bar\t))|^2\]dt'\Big\}.\ea\ee
By Grownwall's inequality, one has
\bel{stability-2}\dbE_\t\[\sup_{\t\les t\les T}|X_1(t\land\bar\t)- X_2(t\land\bar\t)|^2\]\les K|\xi_1-\xi_2|^2.\ee
Since $A_t=\{\bar\t>t\}$, we see that
\bel{2.21}\ba{ll}
\ns\ds K|\xi_1-\xi_2|^2\ges\dbE_\t\[\sup_{\t\les t\les T}|X_1(t\land\bar\t)
-X_2(t\land\bar\t)|^2\]\ges\dbE_\t\[\sup_{\t\les t\les T}|X_1(t\land
\bar\t)-X_2(t\land\bar\t)|^2I_{A_t}\]\\
\ns\ds\qq=\dbE_\t\[\sup_{\t\les t\les T}|X_1(t)-X_2(t)|^2I_{A_t}\]
\ges\[\dbE_\t\(|X_1(t)-X_2(t)|I_{A_t}\)\]^2,\qq t\in[\t,T].\ea\ee
Plugging into \eqref{P(Ac)}, we have (noting $\ell\eta=T-\t$)
\bel{P(Ac)*}\dbP_\t(A^c_T)\les K\(\ell\eta^2+\ell\eta\d+\ell\eta|\xi_1
-\xi_2|+\eta^{p-2\over2}\d^{-p}\dbE\big[\n(\t,T)\big]\)\les K\(\eta+\d+|\xi_1-\xi_2|+\eta^{p-2\over2}\d^{-p}\).\ee
In the above, $\eta$ and $\d$ are free to choose, with $K$ being an
absolute constant. By a truncation argument, we may assume that $\xi_1$ and $\xi_2$ are bounded. On the set $(|\xi_1-\xi_2|=0)$, one has
$P_\t(A^c_T)=0$ since we may first choose $\d>0$ arbitrarily small,
then choose $\eta>0$ even smaller. On the set $(|\xi_1-\xi_2|>0)$,
we may first choose $\d=|\xi_1-\xi_2|$. Then choose $\eta>0$ smaller
so that
$$\eta+\eta^{p-2\over2}\d^{-p}\les|\xi_1-\xi_2|.$$
Since we assume that $p>2$, the above can be achieved. Thus,
\bel{2.23}\dbP_\t(A^c_T)\les K|\xi_1-\xi_2|.\ee
Combining \rf{2.21} and \rf{2.23}, we obtain \rf{P(Ac)+|X-X|}.

\ms

\it Step 3. Existence of solutions. \rm

\ms

Without loss of generality, assume $T-\t=1$. Let $(\t,\xi,\iota)\in
\cD^p$ and $u(\cd)\in\cU^p[\t,T]$ with $p>2$. Define a sequence
$\{(\wt X_\ell(\cd),\wt\a_\ell(\cd))\}_{\ell\ges1}$ on $[\t,T]$ by
$$(\wt X_\ell(s),\wt\a_\ell(s))=(X_\ell(k),\a_\ell(s)),\qq{k\over2^\ell}\les s<{k+1\over2^\ell},$$
where the sequence $\{(X_\ell(k),\a_\ell(k)):k=0,\cds,2^\ell\}$ is
iteratively defined by the following:
$$X_\ell(0)=x,\qq\a_\ell(0)=\iota,$$
and
\bel{X_ell}\left\{\ba{ll}
\ds X_\ell(k+1)=X_\ell(k)+\int_{k\over2^\ell}^{k+1\over2^\ell}b(s,X_\ell(k),
\a_\ell(s),u(s))ds+\int_{k\over2^\ell}^{k+1\over2^\ell}\si(s,X_\ell(k),
\a_\ell(s),u(s))dW(s),\\
\ns\ds\a_\ell(s)=\a_\ell\({k\over2^\ell}\)+\int_{({k\over2^\ell},s]}
\int_\dbR\m(X_\ell(k),\a_\ell(t-),\th)N(dt,d\th),\qq s\in\big[{k\over2^\ell},{k+1\over2^\ell}\big).\ea\right.\ee
We are going to show that the sequence $(\wt X_\ell(\cd),
\wt\a_\ell(\cd))$ has a limit which is a solution to the state equation
\rf{state}.

\ms

First of all, by induction, we are able to show that for some constant $K>0$, the following holds:
\bel{|X|<K}\dbE\[\sup_{0\les k\les2^\ell}|X_\ell(k)|^p\]\les K,\qq
\forall\ell\ges1.\ee
Now we estimate the difference between $(\wt X_\ell(\cd),\wt\a_\ell
(\cd))$ and $(\wt X_{\ell+1}(\cd),\wt\a_{\ell+1}(\cd))$ on time
interval $[{k\over2^\ell}, {k+1\over2^\ell})$. Note that on this interval,
$(X_\ell(\cd),\a_\ell(\cd))$ is given by \rf{X_ell}, whereas $(X_{\ell+1}(\cd),\a_{\ell+1}(\cd))$ has two different forms on the first half $[{k\over2^\ell},{2k+1\over2^{\ell+1}})$ and the second half $[{2k+1\over2^{\ell+1}},{k+1\over2^\ell})$. More precisely, we have the following:
$$\left\{\ba{ll}
\ds X_{\ell+1}(2k+1)=X_{\ell+1}(2k)+\int_{k\over{2^\ell}}^{2k+1\over2^{\ell+1}}
b(s,X_{\ell+1}(2k),\a_{\ell+1}(s),u(s))ds\\ [4mm]
\ns\ds\qq\qq\qq\qq\qq\qq+\int_{k\over2^\ell}^{2k+1\over2^{\ell+1}}
\si(s,X_{\ell+1}(2k),\a_{\ell+1}(s),u(s))dW(s),\\ [4mm]
\ns\ds\a_{\ell+1}(s)=\a_{\ell+1}\big({k\over2^\ell}\big)
+\int_{({k\over2^\ell},s]}\int_\dbR\m(X_{\ell+1}(2k),\a_{\ell+1}(s-),\th)
N(ds,d\th),\qq s\in[{k\over2^\ell},{2k+1\over2^{\ell+1}}),\\
\ea\right.$$
and
$$\left\{\ba{ll}
\ds X_{\ell+1}(2k+2)=X_{\ell+1}(2k+1)+\int_{2k+1\over2^{\ell+1}}^{k+1
\over2^\ell}b(s,X_{\ell+1}(2k+1),\a_{\ell+1}(s),u(s))ds\\
\ns\ds\qq\qq\qq\qq\qq\qq+\int_{2k+1\over2^{\ell+1}}^{k+1\over2^\ell}
\si(s,X_{\ell+1}(2k+1),\a_{\ell+1}(s),u(s))dW(s),\\
\ns\ds\a_{\ell+1}(s)=\a_{\ell+1}({2k+1\over2^{\ell+1}})+\2n\int_{( {2k+1\over{2^{\ell+1}}},s]}\int_\dbR\mu(X_{\ell+1}(2k+1),\a_{\ell+1}(s-),
\th)N(ds,d\th),\q s\in\big({2k+1\over2^{\ell+1}},{k+1\over2^\ell}\big].\\
\ea\right.$$
By a standard calculation (see \rf{E|X-X|}), we have
\bel{ddd-0}\dbE_{k\over2^\ell}|X_{\ell+1}(2k+1)-X_{\ell+1}(2k)|^p\les K\({1\over2^\ell}\)^{p-2\over2}\[{1\over2^\ell}(1+|X_{\ell+1}(2k)|^p)
+\dbE_{k\over2^\ell}\int_{k\over2^\ell}^{2k+1\over2^{\ell+1}}|u(s)|^pds\].\ee
Define
$$\ba{ll}
\ns\ds A_\ell(s)=\{\wt\a_\ell(t)=\wt\a_{\ell+1}(t)\bigm|\t\les t\les s\},\\ %
\ns\ds B_j=\Big\{N\({j+1\over2^{\ell+1}},\cR_0\)-N\({j\over2^{\ell+1}},\cR_0\)=1\Big\}.\ea$$
Note that
$$A_\ell\({k+1\over2^\ell}\)^c=A_\ell\({k\over2^\ell}\)^c\cup \[A_\ell\({k\over2^\ell}\)\cap A_\ell\({2k+1\over2^{\ell+1}}\)^c\]\cup\[A_\ell\({2k+1\over2^{\ell+1}}
\)\cap A_\ell\({k+1\over2^\ell}\)^c\]$$
and in the equations for $\a_{\ell+1}(\cd)$, either $X_{\ell+1}(2k)$
or $X_{\ell+1}(2k+1)$ appears, instead of $X_{\ell+1}(s)$.
Because of this, similar to \rf{2.18}, we have, with the similar
terms involving $\D_{ij}(X_1(t_k))\setminus\D_{ij}(X_1(t)k)$ and $\D_{ij}(X_2(t_k))\setminus\D_{ij}^{-\d}(X_2(t_k))$ absent
(noting \rf{|X|<K}),
$$\ba{ll}
\ns\ds\dbP_{k\over2^\ell}\[A_\ell\({k\over2^\ell}\)\cap A_\ell\({2k+1\over2^{\ell+1}}\)^c\]\les I_{A_\ell({k\over2^\ell})}\dbP_{k\over2^\ell}
\[\Big\{\a_\ell\({2k+1\over2^{\ell+1}}\)\ne\a_{\ell+1}
\({2k+1\over2^{\ell+1}}\)\Big\}\cap B_{2k}\]+{K\over2^{2\ell}}\\
\ns\ds\les{K\over2^\ell}\[\({1\over2^\ell}\)^{p-2\over2}\d^{-p}
+|X_{\ell+1}(2k)-X_\ell(k)|I_{A_\ell({k\over2^\ell})}
+\({1\over2^\ell}\)^{p-4\over2}\d^{-p}\dbE_{{k\over2^\ell}}
\int_{k\over2^\ell}^{2k+1\over2^{\ell+1}}|u(t)|^pdt\]
+{K\over2^{2\ell}}.\ea$$
Likewise,
$$\ba{ll}
\ns\ds\dbP_{2k+1\over2^{\ell+1}}\[A_\ell\({2k+1\over2^{\ell+1}}\)
\cap A_\ell\({k+1\over2^\ell}\)^c\]\les I_{A_\ell({2k+1\over2^{\ell+1}})}\dbP_{2k+1\over2^{\ell+1}}
\[\Big\{\a_\ell\({k+1\over2^\ell}\)\ne\a_{\ell+1}\(
{k+1\over2^\ell}\)\Big\}\cap B_{2k+1}\]+{K\over2^{2\ell}}\\
\ns\ds\les{K\over2^\ell}\[\({1\over2^\ell}\)^{p-2\over2}\d^{-p}
+|X_{\ell+1}(2k+1)-X_\ell(k)|I_{A_\ell({2k+1\over2^{\ell+1}})}
+\({1\over2^\ell}\)^{p-4\over2}\d^{-p}\dbE_{2k+1\over2^{\ell+1}}
\int_{2k+1\over2^{\ell+1}}^{k+1\over2^\ell}|u(t)|^pdt\]
+{K\over2^{2\ell}}.\ea$$
Thus we have
$$\ba{ll}
\ns\ds\dbP_{k\over2^\ell}\[A_\ell\({k+1\over2^\ell}\)^c\]
=I_{A_\ell({k\over2^\ell})^c}+\dbP_{k\over2^\ell}
\[A_\ell\({k\over2^\ell}\)\cap A_\ell\({2k+1\over2^{\ell+1}}\)\]
+\dbP_{k\over2^\ell}\[A_\ell\({2k+1\over2^{\ell+1}}\)\cap A_\ell\({k+1\over2^\ell}\)^c\]\\
\ns\ds\les I_{A_\ell({k\over2^\ell})}+{K\over2^\ell}
\Big\{\({1\over2^\ell}\)^{p-2\over2}\d^{-p}
+\({1\over2^\ell}\)^{p-4\over2}\d^{-p}\dbE_{{k\over2^\ell}}
\int_{k\over2^\ell}^{k+1\over2^\ell}|u(t)|^pdt
+{1\over2^\ell}\\
\ns\ds\qq\qq\qq\qq+|X_{\ell+1}(2k)-X_\ell(k)|I_{A_\ell({k\over2^\ell})}
+\dbE_{k\over2^\ell}\[|X_{\ell+1}(2k+1)-X_\ell(k)|I_{A_\ell({2k+1\over2^{\ell+1}})}\]\Big\}.\ea$$
%
%
%
%
%
%
Taking expectation on both sides of the above, we have
\bel{PAK-0}\ba{ll}
\ns\ds\dbP\[A_\ell\({k+1\over2^\ell}\)^c\]\les \dbP\[A_\ell\({k\over2^\ell}\)^c\]+{K\over2^\ell}\Big\{{1\over2^\ell}
+\({1\over2^\ell}\)^{p-2\over2}\d^{-p}+\({1\over2^\ell}\)^{p-4\over2}
\d^{-p}\dbE\int_{k\over2^\ell}^{k+1\over2^\ell}|u(t)|^pdt\\
\ns\ds\qq\qq\qq\qq\q+\dbE\[|X_{\ell+1}(2k')- X_\ell(k')|I_{A_\ell({k\over2^\ell})}\]+\dbE\[|X_{\ell+1}(2k'+1)- X_\ell(k')|I_{A_\ell({2k+1\over
2^{\ell+1}})}\]\Big\}.\ea\ee
Define stopping time
$$\t_\ell=T\land\inf\{s:\wt\a_\ell(s)\ne\wt\a_{\ell+1}(s)\}.$$
Note that (making use of \rf{ddd-0} and \rf{|X|<K})
$$\ba{ll}
\ns\ds\dbE\[\sup_{0\les k'\les k+1}\Big|\wt X_\ell\big({k'\over2^\ell}\land \t_\ell\big)-\wt X_{\ell+1}\big({k'\over2^\ell}\land\t_\ell\big)\Big|^2\]\les
\dbE\[\sup_{0\les k'\les k}\Big|\wt X_\ell\big({k'\over2^\ell}\land\t_\ell\big)-\wt X_{\ell+1}\big({k'\over2^\ell}\land\t_\ell\big)\Big|^2\]\\
\ns\ds\q+K\dbE\int_{k\over2^\ell}^{2k+1\over2^{\ell+1}}\Big|\wt X_\ell\big({k\over2^\ell}\land\t_\ell\big)-\wt X_{\ell+1}\big({k\over2^\ell}\land \t_\ell\big)\Big|^2ds+K\dbE\int_{2k+1\over2^{\ell+1}}^{k+1\over2^\ell}\Big|\wt X_\ell\big({k\over2^\ell}\land\t_\ell\big)-\wt X_{\ell+1}\big({2k+1\over2^{\ell+1}}\land\t_\ell\big)\Big|^2ds\\
\ns\ds\q\les\(1+{K\over2^\ell}\)\dbE\[\sup_{0\les k'\les k}\Big|\wt X_\ell\big({k'\over2^\ell}\land\t_\ell\big)-\wt X_{\ell+1}\big({k'\over2^\ell}\land\t_\ell\big)\Big|^2\]+{K\over2^\ell}
\dbE\Big|\wt X_{\ell+1}\big({2k+1\over2^{\ell+1}}\land\t_\ell\big)-\wt X_{\ell+1}\big({k\over2^\ell}\land\t_\ell\big)\Big|^2\\
\ns\ds\q\les\(1+{K\over2^\ell}\)\dbE\[\sup_{0\les k'\les k}\Big|\wt X_\ell\big({k'\over2^\ell}\land\t_\ell\big)-\wt X_{\ell+1}\big({k'\over2^\ell}\land\t_\ell\big)\Big|^2\]
+{K\over2^\ell}\[{1\over2^\ell}\(1+\dbE\big|\wt X_{\ell+1}
\big({k\over2^\ell}\land\t_\ell\big)\big|^p\)\\
\ns\ds\qq\qq\qq\qq\qq\qq\qq\qq\qq\qq\qq\qq\qq\qq\qq\qq+\dbE\int_{k\over2^\ell}^{2k+1\over2^{\ell+1}}|u(s)|^pds\]\\
\ns\ds\q\les\(1+{K\over2^\ell}\)\dbE\[\sup_{0\les k'\les k}\Big|\wt X_\ell\big({k'\over2^\ell}\land\t_\ell\big)-\wt X_{\ell+1}\big({k'\over2^\ell}\land\t_\ell\big)\Big|^2\]
+{K\over2^\ell}\[{1\over2^\ell}
+\dbE\big[\n\big({k\over2^\ell},{2k+1\over2^{\ell+1}}\big)\big]\].\ea$$
If we denote the left hand side of the above by $\f_{k+1}$ and denote $\z_k=\dbE\big[\n\big(\t,{k\over2\ell}\big)\big]$, then the above reads
$$\f_{k+1}\les\(1+{K\over2^\ell}\)\f_k+{K\over2^\ell}\[{1\over2^\ell}
+\z_{k+1}-\z_k\],\qq 0\les k\les2^\ell.$$
It is clear that
$$\ba{ll}
\ns\ds\f_{2^\ell}\les\(1+{K\over2^\ell}\)^{2^\ell}\f_0
+{K\over2^\ell}\sum_{k=0}^{2^\ell-1}\(1+{K\over2^\ell}\)^k
\[{1\over2^\ell}+\z_{k+1}-\z_k\]\\
\ns\ds\q={K\over2^\ell}{1\over2^\ell}
\Big\{{2^\ell\over K}\[\(1+{K\over2^\ell}\)^{2^\ell}-1\]\Big\}
+{K\over2^\ell}\sum_{k=0}^{2^\ell-1}\(1+{K\over2^\ell}\)^k
[\z_{k+1}-\z_k]\\
\ns\ds\q\les{e^K\over2^\ell}+{K\over2^\ell}\[\sum_{k=0}^{2^\ell-1}
\(1+{K\over2^\ell}\)^k\z_{k+1}-\sum_{k=1}^{2^\ell-1}
\(1+{K\over2^\ell}\)^k\z_k\]\\
\ns\ds\q={e^K\over2^\ell}+{K\over2^\ell}\[
\(1+{K\over2^\ell}\)^{2^\ell-1}\z_{2^\ell}+\sum_{k=1}^{2^\ell-1}
\(1+{K\over2^\ell}\)^{k-1}\z_k-\(1+{K\over2^\ell}\)
\sum_{k=1}^{2^\ell-1}\(1+{K\over2^\ell}\)^{k-1}\z_k\]\\
\ns\ds\q\les{e^K\over2^\ell}+{K\over2^\ell}e^K\dbE\int_t^T|u(s)|^pds
\les{K\over2^\ell}.\ea$$
Since $A_{k'\over2^\ell}=\{\t_\ell\ges{k'\over2^\ell}\}$, we obtain
$$\ba{ll}
\ns\ds{K\over2^\ell}\ges\dbE\[
\sup_{0\les k'\les2^\ell}\Big|\wt X_\ell\big({k'\over2^\ell}\land\t_\ell
\big)-\wt X_{\ell+1}\big({k'\over2^\ell}\land\t_\ell\big)\Big|^2\\
\ns\ds\ges\dbE\[\sup_{0\les k'\les2^\ell}\Big|\wt X_\ell\big({k'\over
2^\ell}\land\t_\ell\big)-\wt X_{\ell+1}\big({k'\over2^\ell}\land \t_\ell\big)\Big|^2I_{A_\ell({k'\over2^\ell})}\]
=\dbE\[\sup_{0\les k'\les2^\ell}\Big|X_\ell(k')-X_{\ell+1}(2k')\Big|^2I_{A_\ell(
{k'\over2^\ell})}\].\ea$$
Likewise,
$${K\over2^\ell}\ges\dbE\[\sup_{0\les k'\les2^\ell}\Big|X_\ell(k')-X_{\ell+1}(2k'+1)\Big|^2I_{A_\ell(
{2k'+1\over2^\ell+1})}\].$$
Plugging the above into \eqref{PAK-0}, we have
\bel{PAK-1}\dbP\[A_\ell\({k+1\over2^\ell}\)^c\]\les \dbP\[A_\ell\({k\over2^\ell}\)^c\]+{K\over2^\ell}\Big\{{1\over2^\ell}
+\({1\over2^\ell}\)^{p-2\over2}\d^{-p}+\({1\over2^\ell}\)^{p-4\over2}
\d^{-p}(\z_{k+1}-\z_k)+\({1\over2^\ell}\)^{1\over2}\Big\}.\ee
Hence,
\bel{PA(T)c}\ba{ll}
\ns\ds\dbP\[A_\ell(T)^c\]\les K\[{1\over2^\ell}
+\({1\over2^\ell}\)^{p-2\over2}\d^{-p}+\({1\over2^\ell}\)^{1\over2}
+\({1\over2^\ell}\)^{p-2\over2}
\d^{-p}\z_{2^\ell}\]\les K\[\({1\over2^\ell}\)^{1\over2}+\({1\over2^\ell}\)^{p-2\over2}\d^{-p}\].\ea\ee
Since $p>2$, we see that
\bel{sumP}\sum_{\ell=1}^\infty\dbP\big[A_\ell(T)^c\big]<\infty,\ee
and
\bel{sum|X-X|}\ba{ll}
\ns\ds\sum_{\ell=1}^\infty\dbE\[\sup_{\t\les s\les T}\big|\wt X_\ell(s)
-\wt X_{\ell+1}(s)\big|^2\]\les\sum_{\ell=1}^\infty\dbE\[\sup_{\t\les s\les T}\big|\wt X_\ell(s)-\wt X_{\ell+1}(s)|^2I_{A_\ell(T)}\]\\
\ns\ds\qq\qq\qq\qq\qq\qq+\sum_{\ell=1}^\infty\dbE\[\sup_{\t\les s\les T}
\big(1+|\wt X_\ell(s)|^2+|\wt X_{\ell+1}(s)|^2
+\int_\t^T|u(s)|^2ds\big)I_{A_\ell(T)^c}\]\\
\ns\ds\les K\sum_{\ell=1}^\infty{1\over2^\ell}+K\sum_{\ell=1}^\infty
\(1+\dbE\[\sup_{\t\les s\les T}(|\wt X_\ell(s)|^p+|\wt X_{\ell+1}(s)|^p)+\int_\t^T|u(s)|^pds\]\)^{2\over p}\big\{\dbP\big[A_\ell(T)^c\big]\big\}^{1-{2\over p}}\\
\ns\ds\les K\sum_{\ell=1}^\infty\Big\{{1\over2^\ell}+\[\({1\over2^\ell}\)^{1\over2}
+\({1\over2^\ell}\)^{p-2\over2}\d^{-p}\]^{p-2\over p}\Big\}<\infty.
\ea\ee
Clearly, \rf{sum|X-X|} implies that there exist an $\wt X(\cd)$ such that
$$\lim_{\ell\to\infty}\dbE\[\sup_{\t\les s\les T}|\wt X_\ell(s)-\wt X(s)|^2\]=0.$$
By Borel-Cantelli's lemma, condition \rf{sumP} implies that for some $\wt\a(\cd)$,
\bel{}\lim_{\ell\to\infty}\wt\a_\ell(s)=\wt\a(s),\qq\forall s\in[\t,T],~\as\ee
Actually, if we let
$$\BA(T)^c=\limsup_{\ell\to\infty}A_\ell(T)^c\equiv\bigcap_{k=1}^\infty
\bigcup_{\ell=k}^\infty A_\ell(T)^c,$$
Then
$$\dbP\big(\BA(T)^c\big)\les\sum_{\ell=k}^\infty\dbP\big[A_\ell(T)^c\big]\to0,
\qq\forall k\ges1.$$
Thus $\BA(T)^c$ is a $\dbP$-null set. Now, any $\o\in\O$ with
$$\o\in\BA(T)=\[\bigcap_{k=1}^\infty\bigcup_{\ell=k}^\infty A_\ell(T)^c\]^c
=\bigcup_{k=1}^\infty\bigcap_{k=\ell}^\infty A_\ell(T),$$
if and only if there exists a $k=k(\o)\ges1$ such that
$$\o\in\bigcap_{\ell=k(\o)}^\infty A_\ell(T)\q\iff\q\wt\a_\ell(s,\o)=\wt\a_{\ell+1}(s,\o),\q\forall\ell\ges k(\o).$$
Therefore, the limit $\wt\a(\cd)$ can be defined by
$$\wt\a(s,\o)=\left\{\1n\ba{ll}
\ds\wt\a_\ell(s,\o),\qq\ell\ges k(\o),\qq\o\in\BA(T),\\
\ns\ds1,\qq\qq\qq\qq\qq\q\o\in\BA(T)^c.\ea\right.$$
Now we want to show that $(\wt X(\cd),\wt\a(\cd))$ is a solution of \rf{state}. Knowing that
$$\big\{\wt\a_\ell(s)\ne\wt\a(s),~\t\les s\les T\big\}\subseteq\bigcup_{k\ges\ell}A_\ell(T)^c\equiv\BA_\ell(T)^c,\qq
\lim_{\ell\to\infty}\dbP\big(\BA(T)^c\big)=0,$$
we have
$$\ba{ll}
\ns\ds\dbE\[\sup_{\t\les s\les T}\Big|\wt X(s)-\xi-\int_\t^sb(t,\wt X(t),\wt\a(t),u(t))dt-\int_\t^s\si(t,\wt X(t),\wt\a(t),u(t))dW(t)\Big|^2\]\\
\ns\ds\les K\dbE\[\sup_{\t\les s\les T}\big|\wt X(s)-\wt X_\ell(s)|^2+\sup_{\t\les s\les T}\dbE\Big|\wt X_\ell(s)-\xi-\int_\t^sb(t,\wt X_\ell(t),\wt\a_\ell(t),u(t)dt\\
\ns\ds\qq\qq\qq\qq\qq\qq\qq\qq\qq\qq\qq-\int_\t^s\si(t,\wt X_\ell(t),\wt \a_\ell(t),u(t)dW(t)\Big|^2\\
\ns\ds\qq+K\int_\t^T|\wt X(t)-\wt X_\ell(t)|^2dt+K\dbE\int_\t^T(1+|\wt X(t)|^2+|\wt X_\ell(t)|^2+|u(t)|^2)I_{\BA_\ell(T)^c}dt\]\\
\ns\ds\les\dbE\[\sup_{\t\les s\les T}|\wt X(s)-\wt X_\ell(s)|^2+\dbE\[\sup_{\t\les s\les T}\Big|\wt X_\ell(s)-\xi-\int_\t^sb(t,\wt X_\ell(t),\wt \a_\ell(t),u(t)dt\\
\ns\ds\qq\qq\qq\qq\qq\qq\qq\qq\qq\qq\qq-\int_\t^s\si(t,\wt X_\ell(t),\wt\a_\ell(t),u(t)dW(t)\Big|^2\]\\
\ns\ds\qq+K\int_\t^T|\wt X(t)-\wt X_\ell(t)|^2dt+K\(\dbE\[\sup_{\t\les t\les T}(1+|\wt X(t)|^p+|\wt X_\ell(t)|^p)\]\)^{p\over2}\(\dbP[\BA_\ell(T)^c]\)^{1-{2\over p}}\to0,\ea$$
as $\ell\to\infty$. Thus $\wt X(\cd)$ is the solution of the first equation in \rf{state}. Moreover, by Step 1, we see that $\wt X(\cd)$ satsifies \rf{|X|} and \rf{|X-X|}. Further, noting the second equation in \rf{X_ell}, the convergence of $\wt X_\ell(\cd)$ to $\wt X(\cd)$ and $\dbP[\BA_\ell(T)^c]\to0$, together with \rf{D<Kd}, we obtain the following:
$$\ba{ll}
\ns\ds\dbE\[\sup_{\t\les s\les T}|\wt\a(s)-\iota-\int_{(\t,s]}\m(\wt X(t),\wt\a(t-),\th)N(dt,d\th)|\]\\
\ns\ds\les\dbE\[\sup_{\t\les s\les T}|\wt\a(s)-\wt\a_\ell(s)|\]+\dbE\[\sup_{t\les s\les T}\Big|\wt \a_\ell(s)-\iota-\int_{(\t,s]}\m(\wt X(t),\wt\a(t-),\th)N(dt,d\th)\Big|\,\]\\
\ns\ds\les\dbE\[\sup_{\t\les s\les T}|\wt\a(s)-\wt\a_\ell(s)|\]+\dbE\[\sup_{\t\les s\les T}\Big|\int_{(\t,s]}\(\m(\wt X(t),\wt\a(t-),\th)-\m(\wt X_\ell(t),\wt\a_\ell(t-),\th)\)N(dt,d\th)\Big|\,\]\to0,\ea$$
as $\ell\to\infty$. Thus, $\wt\a(\cd)$ satisfies the second equation in \rf{state}. \endpf

\ms

The following result is the It\^o's formula for the solution $(X(\cd),\a(\cd))$ of our state equation \rf{state}, which is called a
controlled {\it stochastic diffusion with regime switching}. The proof
is omitted here.

\bt{well-posedness-state equation} \sl Let {\rm(H0)--(H1)} hold. Let $(\t,\xi,\iota)\in\cD^p$, and $u(\cd)\in\cU^p[\t,T]$, with $p\ges2$.
Let $(X(\cd),\a(\cd))$ be the unique solution to \rf{state}.
Let $\Th:[t,T]\times\dbR^n\times M\to\dbR$ with $(t,x)\mapsto\Th(t,x,i)$ being $C^{1,2}$. Then the following holds:
\bel{Ito}\ba{ll}
\ns\ds d\Th(s,X(s),\a(s))=\[\Th_s(s,X(s)\a(s))+\cA^{u(s)}\Th(s,X(s),\a(s))\]ds\\ [3mm]
\ns\ds\qq\qq\qq\qq+\Th_x(s,\1n X(s),\a(s))\si(s,\1n X(s),\a(s),u(s))dW(s)\\ [2mm]
\ns\ds\qq\qq\qq\qq+\int_\dbR\[\Th\big(s,X(s),\a(s-)\1n+\1n\m(X(s),\a(s-),\th)\big)\1n
-\1n\Th(s,X(s),\a(s-))\]\wt N(ds,d\th),\ea\ee
where $\cA^u$ is a differential operator defined by the following:
\bel{cA}\ba{ll}
\ns\ds\cA^u\Th(s,x,i)={1\over2}\tr\[\si(s,x,i,u)^\top\Th_{xx}(s,x,i)
\si(s,x,i,u)\]\1n+\1n\Th_x(s,x,i)b(s,x,i,u)+\big[Q(x)\Th(s,x,\cd)\big]_i,\ea\ee
where $Q(x)=[q_{ij}(x)]_{m\times m}$, and $[Q(x)\Th(s,x,\cd)]_i$ stands for the $i$-th component of the vector $Q(x)\Th(s,x,\cd)$:
$$Q(x)\Th(s,x,\cd)=Q(x)\begin{pmatrix}\Th(s,x,1)\\ \Th(s,x,2)\\ \vdots\\ \Th(s,x,m)\end{pmatrix}.$$

\et

As a convention, for a differentiable function $\Th:\dbR^n\to\dbR$, the gradient $\Th_x(x)$, taking values in
$\dbR^{1\times n}$, is a row vector and the Hessian $\Th_{xx}(x)$, taking values in $\dbS^n$, is a symmetric matrix.

\ms

In what follows, we will fix a $p>2$, and simply denote
$$\cD[\t_1,\t_2]=\cD^p[\t_1,\t_2],\q\cD=\cD^p[0,T],\q
\q\cU[\t_1,\t_2]=\cU^p[\t_1,\t_2],\q\cU=\cU^p[0,T].$$

Next, for the recursive cost functional (\ref{BSDE1})--(\ref{cost}), we introduce the following hypothesis.

\ms

{\bf(H2)} Let $h:[0,T]\times\dbR^n\times M\to\dbR$ and $g:[0,T]\times[0,T]\times\dbR^n\times M\times\dbR\times\dbR^d\times\dbR\times U\to\dbR$,
\bel{Lip-h-g}\left\{\2n\ba{ll}
\ns\ds|h(t_1,x_1,i)-h(t_2,x_2,i)|\les L\big(|t_1-t_2|+|x_1-x_2|\big)\\
\ns\ds|g(t_1,s,x_1,i,y_1,z_1,\g_1,u)-g(t_2,s,x_2,i,y_2,z_2,\g_2,u)|\\
\ns\ds\qq\les L\big(|t_1-t_2|+|x_1-x_2|+|y_1-y_2|+|z_1-z_2|+|\g_1-\g_2|\big).\ea\right.\ee

\ms

Under (H0)--(H2), for any initial triple $(\t,\xi,\iota)\in\cD$, and $u(\cd)\in\cU[\t,T]$, state equation (\ref{state}) admits a unique
solution $(X(\cd),\a(\cd))$, and BSDE (\ref{BSDE1}) admits a unique
adapted solution $(Y(\t;\cd),Z(\t;\cd),\G(\t;\cd,\cd))$ (see \cite{Situ2006} Theorem 234, p.216, and also \cite{Kruse-Popier2016, Kruse-Popier2017}). Thus, $s\mapsto Y(\t;s)$ is $\dbF$-adapted (on $[\t,T]$) and depends on $(\t,\xi,\iota,u(\cd))$ through $(X(\cd),\a(\cd))$; and for $t\in[\t,T]$, one has
\bel{cost3}\ba{ll}
\ns\ds Y(\t;t)=\dbE_t\[h(\t;X(T),\a(T))\1n+\2n\int_t^T\2n g\big(\t;s,X(s),\a(s),Y(\t;s),Z(\t;s),\int_\dbR\G(\t,s,\th)\pi(d\th),u(s)\big)ds\].\ea\ee
We regard $Y(\t;t)$ as a utility/disutility (associated with the triple $(X(\cd),\a(\cd),u(\cd))$) at time $t$, which depends on its ``future'' value $Y(\t;s)$ (for $s\in(t,T]$). Because of this, we call $(Y(\t;\cd),Z(\t;\cd)$, $\G(\t;\cd\,,\cd))$ a {\it recursive utility/disutility process}. This is a generalization of the so-called {\it stochastic differential utility} introduced by Duffie--Epstein \cite{Duffie-Epstein1992a, Duffie-Epstein1992b} (see also \cite{Duffie-Lions1992, Lazrak-Quenez2003, Lazrak2004, Wei-Yong-Yu2016}). By taking $t=\t$ in (\ref{cost3}), we obtain our cost functional (\ref{cost2}).

\ms

From the above, we see that under (H1)--(H2), $J(\t,\xi,\iota;u(\cd))$ is well-defined for all $(\t,\xi,\iota)\in\cD$ and $u(\cd)\in\cU[\t,T]$. Hence, Problem (N) is well-formulated.

\ms

\subsection{Time-consistent case}

We let $\t=t\in[0,T)$ be fixed and consider Problem (N) on $[t,T]$. In such a case, the problem will be time-consistent. Actually, for state equation \eqref{state}, when the cost is non-recursive, the dynamic programming approach, involving Hamilton-Jacobi-Bellman (HJB for short) equation, can be found in \cite{Oksendal-Sulem2005} and the related maximum principle can be found in \cite{Framstad-Oksenda-Sulem2004} and \cite{Donnelly2011}. In this subsection, we will present a verification theorem and a representation of the optimal recursive process for the case that the cost functional is recursive. These results will be useful below.

\ms

For fixed $t\in[0,T)$. Consider the state equation \eqref{state}, with the initial triple $(\t,\xi,\iota)\in\cD[t,T]$, the control $u(\cd)\in\cU[\t,T]$, and with the recursive cost functional
\bel{cost4}\ba{ll}
\ns\ds\wt J(t;\t,\xi,\iota;u(\cd))=Y(t;\t)=\dbE_\t\[h(t;X(T),\a(T))\\
\ns\ds\qq\qq\qq\qq\qq\qq+\int_\t^Tg\big(t,s,X(s),\a(s),Y(t;s),Z(t;s),\int_\dbR\G(t;s,\th)\pi(d\th),u(s)\big)ds
\],\ea\ee
Here, $(Y(t;\cd),Z(t;\cd),\G(t;\cd\,,\cd))$ is the adapted solution of the following BSDE:
\bel{BSDE2*}\left\{\2n\ba{ll}
\ns\ds dY(t;s)=-g\big(t,s,X(s),\a(s),Y(t;s),Z(t;s),\int_\dbR\G(t;s,\th)\pi(d\th),u(s)\big)ds\\
\ns\ds\qq\qq\qq+Z(t;s)dW(s)+\int_\dbR\G(t;s-,\th)\wt N(ds,d\th),\qq s\in[\t,T],\\
\ns\ds Y(t;T)=h(t;X(T),\a(T)).\ea\right.\ee
Note that the above equation is the same as \eqref{BSDE2}, with $\t$ replaced by $t$. However, we only solve it on $[\t,T]$ since $(X(s),\a(s))$ are only defined on this interval (and $\t\ges t$). Note that
$$\wt J(t;t,\xi,\iota;u(\cd))=J(t;\xi,\iota;u(\cd)).$$
We introduce the following problem.

\ms

\bf Problem (C$_t$). \rm For any given $(\t,\xi,\iota)\in\cD[t,T]$, find a $\bar u(\cd)\in\cU[\t,T]$ such that
$$\wt J(t;\t,\xi,\iota;\bar u(\cd))=\inf_{u(\cd)\in\cU[\t,T]}\wt J(t;\t,\xi,\iota;u(\cd))\equiv V(t;\t,\xi,\iota).$$

\ms

For the above problem, we introduce the Hamiltonian $\dbH:[0,T]\times[0,T]\times
\dbR^n\times M\times\dbR^m\times\dbR^{m\times n}\times(\dbS^n)^m\times U\to\dbR$ as follows:
\bel{H}\ba{ll}
\ns\ds\dbH(t;s,x,i,v,\Bp,\BP,u)=\Bp^ib(s,x,i,u)+{1\over2}\tr\big[\si(s,x,i,u)^\top \BP^i\si(s,x,i,u)\big]\\
\ns\ds\qq\qq\qq\qq\qq\qq+\sum_{j=1}^mq_{ij}(x)v^j+g(t,s,x,i,v^i,\Bp^i\si(s,x,i),
\sum_{j=1}^mq_{ij}(x)v^j,u)\\
\ns\ds\qq\qq\qq\qq\q\equiv\Bp^ib(s,x,i,u)+{1\over2}\tr\big[\si(s,x,i,u)^\top \BP^i\si(s,x,i,u)\big]
+\big[Q(x)v\big]_i\\
\ns\ds\qq\qq\qq\qq\qq\qq+g\big(t,s,x,i,v^i,\Bp^i\si(s,x,i),
\big[Q(x)v\big]_i,u\big),\ea\ee
where $v=(v^1,\cds,v^m)^\top\in\dbR^n$, $\Bp=\big((\Bp^1)^\top,(\Bp^2)^\top,\cds,(\Bp^m)^\top
\big)^\top\in\dbR^{m\times n}$ with $\Bp^i\in\dbR^{1\times n}$, and
$\BP=(\BP^1,\BP^2,\cds$, $\BP^n)\in(\dbS^n)^m$ with $\BP^i\in\dbS^n$.
Comparing the above with the definition of $\cA^u$ (see \rf{cA}), we see that
\bel{cA=H}\ba{ll}
\ns\ds\dbH\big(t;s,x,i,V(s,x,\cd),V_x(s,x,i),V_{xx}(s,x,i),u\big)\\
\ns\ds=V_x(s,x,i)b(s,x,i,u)+{1\over2}\tr\[\si(s,x,i,u)^\top V_{xx}(s,x,i)\si(s,x,i,u)\]\\
\ns\ds\qq+\big[Q(x)V(s,x,\cd)\big]_i+g(t,s,x,i,V(s,x,i),V_x(s,x,i)\si(s,x,i,u),
[Q(x)V(s,x,\cd)]_i,u\big)\\
\ns\ds\equiv\cA^uV(s,x,i)+g(t,s,x,i,V(s,x,i),V_x(s,x,i)\si(s,x,i,u),
[Q(x)V(s,x,\cd)]_i,u\big).\ea\ee
For convenience, we may abbreviate the above as
\bel{abb}\dbH=\cA V+g.\ee
We now introduce the following hypothesis.

\ms

{\bf(H3)} The set $U$ is a closure of a domain in $\dbR^{\bar n}$. There exists a map $\psi:[0,T]\times[0,T]\times\dbR^n\times M\times\dbR^m\times\dbR^{1\times n}\times\dbS^n\to U$ such that
$$\ba{ll}
\ns\ds\psi(t;s,x,i,v,\Bp,\BP)\in\hb{argmin}~\dbH(t;s,x,i,v,\Bp,\BP,u)\\ [2mm]
\ns\ds:=\Big\{\bar u\in U\bigm|\dbH(t;s,x,i,v,\Bp,\BP,\bar u)=\min_{u\in U}
\dbH(t;s,x,i,v,\Bp,\BP,u)\Big\}\\ [3mm]
\ns\ds~~~~~~~(t,s,x,i,v,\Bp,\BP)\in[0,T]\times[0,T]\times\dbR^n\times M\times\dbR^m\times\dbR^{1\times n}\times\dbS^n,\ea$$
and $\psi$ is $C^1$ with bounded derivatives.

\ms

For convenience, in (H3), the map $\psi(\cd)$ is assumed to be $C^1$ with bounded derivatives. This will make some of the derivations simpler later, and such a condition might be more than enough. But, for now, we prefer not to get into the most generality in this aspect. More general cases are still open (see Section 6 for more comments). The above can also be written as
\bel{psi}\dbH\big(t;s,x,i,v,\Bp,\BP,\psi(t;s,x,i,v,\Bp,\BP)\big)=\min_{u\in U}\dbH(t;s,x,i,v,\Bp,\BP,u).\ee
Note that in the case that $\si$ is independent of $u$, we have
\bel{si=si}\psi(t;s,x,i,v,\Bp,\BP)=\psi(t;s,x,i,v,\Bp),\ee
which is independent of $\BP$. Such a case will be assumed in Section 5. The following is called a {\it verification theorem}.

\bt{verification theorem} \sl Let {\rm(H1)--(H2)} hold. Suppose that $V(t;\cd\,,\cd\,,\cd):[0,T]\times \dbR^n\times M\to\dbR$ is the unique classical solution to the following HJB equation:
\bel{HJB1}\left\{\2n\ba{ll}
\ns\ds V_s(t;s,x,i)+\inf_{u\in U}\dbH(t;s,x,i,V(t;s,x,\cd),V_x(t;s,x,i),V_{xx}(t;s,x,i),u)=0, \\ [2mm]
\ns\ds V(t;T,x,i)=h(t,x,i).\ea\right.\ee
Then
\bel{V<J}V(t;\t,\xi,\iota)\les\wt J(t;\t,\xi,\iota;u(\cd)),\q\forall(\t,\xi,\iota)\in\cD[t,T]\times\dbR^n\times M,\q u(\cd)\in\cU[\t,T].\ee
If, in addition, {\rm(H3)} also holds, $(\t,\xi,\iota)\in\cD[t,T]$ is given and $(\bar X(t;\cd),\bar\a(t;\cd),\bar u(t;\cd))$ is a state-control triple of \eqref{state} with
$(\bar X(t;\t),\bar\a(t;\t))=(\xi,\iota)$ such that
\bel{bar u=Psi}\bar u(t;s)=\bar\Psi(t;s,\bar X(t;s),\bar\a(t;s)),\qq s\in[\t,T],\ee
where
\bel{Psi}\bar\Psi(t;s,x,i)=\psi(t;s,x,i,V(t;s,x,\cd),V_x(t;s,x,i),V_{xx}(t;s,x,i)),\q s\in[\t,T].\ee
Then
\bel{V=J}V(t;\t,\xi,\iota)=\wt J(t;\t,\xi,\iota,\bar u(t;\cd)),\ee
and $(\bar X(t;\cd),\bar\a(t;\cd);\bar u(t;\cd))$ is an open-loop optimal triple of Problem (C$_t$) for the initial triple $(\t,\xi,\iota)$.

\et

\it Proof. \rm For given $(\t,\xi,\iota)\in\cD[t,T]$ and any $u(\cd)\in\cU[\t,T]$, let $(X(\cd),\a(\cd))$ be the corresponding state process. By generalized It\^o's formula, 
$$\ba{ll}\ns\ds h(t;X(T),\a(T))-V(t;\t,\xi,\iota)=V(t;T,X(T),\a(T))-V(t;\t,\xi,\iota)\\
\ns\ds=\int_\t^TV_s(t;s,X(s),\a(s))+\cA^{u(s)}V(t;s,X(s),\a(s))ds\\
\ns\ds\qq+\int_\t^TV_x(t;s,X(s),\a(s))\si(s,X(s),\a(s),u(s))dW(s)\\
\ns\ds\qq+\int_{(\t,T]}[V(t;s,X(s),\a(s-)+\m(X(s),\a(s-),\th))-V(t;s,X(s),\a(s-))]\wt N(ds,d\th).\ea$$
Suppressing $s$ in $X(s)$, $\a(s)$, $u(s)$, etc., and suppressing all the arguments but $u$ in $\dbH$ and $g$, we have (recalling \eqref{abb})
$$\ba{ll}
\ns\ds V_s(t;s,X,\a)+\cA^uV(t;s,X,\a)=V_s(t;s,X,\a)+\dbH(u)-g(u)\\
\ns\ds\qq\qq\qq\qq\qq\qq\qq\ges\big[V_s(t;s,X,\a)+\inf_{u\in U}\dbH(u)\big]-g(u)=-g(u).\ea$$
Thus,
$$\ba{ll}
\ns\ds V(t;\t,\xi,\iota)=\dbE_\t\[h(t;X(T),\a(T))-\int_\t^TV_s(t;s,X(s),\a(s))+\cA^{u(s)}V(t;s,X(s),\a(s))ds\]\\
\ns\ds\qq\qq\q\les\dbE_\t\[h(t,X(T),\a(T)+\int_t^Tg(\,\cds,u(s))ds\]=\wt J(t;\t,\xi,\iota,u(\cd)).\ea$$
Next, if (H3) holds and $\bar u(\cd)$ is defined by \eqref{bar u=Psi}, then with $u(\cd)=\bar u(\cd)$, one has equalities in the above derivation. Hence, \eqref{V=J} holds. \endpf

\ms

Note that although state equation \eqref{state} is independent of $t$, due to the fact that $\bar\Psi(\cd)$ is depending on $t$ which leads to $\bar u(t;\cd)$ defined by \eqref{bar u=Psi} depending on $t$. Consequently, $(\bar X(t;\cd),\bar\a(t;\cd))$ will also depend on $t$. Now, for the optimal 6-tuple $\big(\bar X(t;\cd),\bar\a(t;\cd),\bar Y(t;\cd),\bar Z(t;\cd),\bar\G(t;\cd\,,\cd),\bar\Psi(t;\cd\,,\cd\,,\cd)\big)$ obtained above, we have the following coupled forward-backward stochastic differential equation (FBSDE, for short), called the {\it optimality system} on $[t,T]$:
\bel{FBSDE1}\left\{\2n\ba{ll}
\ns\ds d\bar X(t;s)=b\big(s,\bar X(t;s),\bar\a(t;s),\bar\Psi(t;s,\bar X(t;s),\bar\a(t;s))
\big)ds\\ [2mm]
\ns\ds\qq\qq\qq+\si\big(s,\bar X(t;s),\bar\a(t;s),\bar\Psi(t;s,\bar X(t;s),\bar\a(t;s))\big)dW(s),  \\
\ns\ds d\bar\a(t;s)=\int_\dbR\m(\bar X(t;s),\bar\a(t;s-),\th)N(ds,d\th),\\
\ns\ds d\bar Y(t;s)=-g\big(t,s,\bar X(t;s),\bar\a(t;s),\bar Y(t;s),\bar Z(t;s),\int_\dbR\bar\G (t;s,\th)\pi(d\th),\bar\Psi(t;s,\bar X(t;s),\bar\a(t;s))\big)ds\\ [1mm]
\ns\ds\ ~~~~~~~~~~~~+\bar Z(t;s)dW(s)+\int_\dbR\bar\G(t;s-,\th)\wt N(ds,d\th),\\ [1mm]
\ns\ds\bar X(t;\t)=\xi,\qq\bar\a(t;\t)=\iota,\qq\bar Y(t;T)=h(t,\bar X(t;T),\bar\a(t;T)).\ea\right.\ee
Inspired by \cite{Ma-Protter-Yong1994, Ma-Yong1999}, we have the following representation theorem.

\bt{representation} \sl Let {\rm(H1)--(H3)} hold. Let $(\bar X(t;\cd),\bar\a(t;\cd),\bar Y(t;\cd),\bar Z(t;\cd),
\bar\G(t;\cd\,,\cd))$ be the adapted solution to the FBSDE \eqref{FBSDE1}. Let $\Th(t;\cd\,,\cd\,,\cd)$
be the classical solution to the following partial differential equation (PDE, for short):
\bel{TH-1}\left\{\2n\ba{ll}
\ns\ds\Th_s(t;s,x,i)+\bar\cA\Th(t;s,x,i)+g\big(t,s,x,i,\Th(t;s,x,i),\Th_x(t;s,x,i)\si(s,x,i,\bar\Psi(t;s,x,i)),\\
\ns\ds\qq\qq\qq\qq[Q(x)\Th(t;s,x,\cd)]_i,\bar\Psi(t;s,x,i)\big)=0,\qq(s,x,i)\in[t,T]\times\dbR^n\times M,\\
\ns\ds\Th(t;T,x,i)=h(t,x,i),\qq(x,i)\in\dbR^n\times M,\ea\right.\ee
where
$$\ba{ll}
\ns\ds\bar\cA\Th(t;s,x,i)={1\over2}\tr\[\si\big(s,x,i,\bar\Psi(t;s,x,i)\big)^\top
\Th_{xx}(t;s,x,i)\si\big(s,x,i,\bar\Psi(t;,s,x,i)\big)\]\\
\ns\ds\qq\qq\qq\qq+\Th_x(t;s,x,i)b(s,x,i,\bar\Psi(t;s,x,i))+\big[Q(x)\Th(t;s,x,\cd))\big]_i\,.\ea$$
%
%
%
%
Then the following representation hold:
\bel{Y=Th}\left\{\2n\ba{ll}
\ns\ds\bar Y(t;s)=\Th(t;s,\bar X(t;s),\bar\a(t;s)),\\
\ns\ds\bar Z(t;s)\1n=\1n\Th_x(t;s,\bar X(t;s),\bar\a(t;s))\si\big(s,\bar X(t;s),\bar\a(t;s),\bar\Psi(t;s,\bar X(t;s),\bar\a(t;s))\big),\\
\ns\ds\bar\G(t;s,\th)=\Th\big(t;s,\bar X(t;s),\bar\a(t;s)+\m(\bar X(t;s),\bar\a(t;s),\th)\big)-\Th(t;s,\bar X(t;s),\bar\a(t;s)),\ea\right.\q s\in[t,T].\ee

\et

\it Proof. \rm Let $(\bar X(t;\cd),\bar\a(t;\cd))$ be the solution to the first two equations in \rf{FBSDE1} and $\Th(t;\cd\,,\cd\,,\cd)$ be the solution to the PDE \rf{TH-1}. Then we define $\bar Y(t;\cd)$, $\bar Z(t;\cd)$, $\bar\G(t;\cd\,,\cd)$ by \rf{Y=Th}. Note that
$$i+\m(x,i,\th)=i+\sum_{j=1}^m(j-i)I_{\D_{ij}(x)}(\th)=\sum_{j=1}^mjI_{\D_{ij}(x)}(\th).$$
Thus,
$$\ba{ll}
\ns\ds\int_\dbR\[\Th(t;s,x,i+\m(x,i,\th))-\Th(t;s,x,i)\]\pi(d\th)=\2n\int_\dbR\[\Th\big(t;s,x,\sum_{j=1}^mjI_{\D_{ij}(x)}(\th)\big)-\Th(t;s,x,i)\]\pi(d\th)\\
\ns\ds=\sum_{j=1}^m\Th(t;s,x,j)\pi(\D_{ij}(x))-\Th(t;s,x,i)=\sum_{j\ne i}\Th(t;s,x,j)q_{ij}(x)+\(\pi\big(\D_{ii}(x)\big)-1\)\Th(t;s,x,i)\\
\ns\ds=\sum_{j=1}^m\Th(t;s,x,j)q_{ij}(s)=\big[Q(x)\Th(t;s,x,\cd)\big]_i.\ea$$
This leads to (note that $\bar\G(\t;\cd\,,\cd)$ is defined by \rf{Y=Th})
$$\int_\dbR\bar\G(t;s,\th)\pi(d\th)=\big[Q(\bar X(t;s))\Th(t;s,\bar X(t;s),\cd)\big]_{\bar\a(t;s)}.$$
Let us denote
$$\bar\si(t;s)=\si\big(s,\bar X(t;s),\bar\a(t;s),\bar\Psi(t;s,\bar X(t;s),\bar\a(t;s))\big).$$
Then by generalized It\^o's formula and PDE \rf{TH-1}, one has
$$\ba{ll}
\ns\ds d\bar Y(t;s)=d\Th(t;s,\bar X(t;s),\bar\a(t;s))\\ [1mm]
\ns\ds=[\Th_s(t;s,\bar X(t;s),\bar\a(t;s))+\bar\cA\Th(t;s,\bar X(t;s),\bar\a(t;s))]ds+\Th_x(t;s,\bar X(t;s),\bar\a(t;s))\bar\si(t;s)dW(s)\\ [2mm]
\ns\ds\q+\int_\dbR\big[\Th\big(t;s,\bar X(t;s),\bar\a(t;s-)+\m(\bar X(t;s),\bar\a(t;s-),\th)\big)-\Th(t;s,\bar X(t;s),\bar\a(t;s-))\big]\wt N(ds,d\th)\ea$$
$$\ba{ll}
\ns\ds=-g\big(t,s,\bar X(t;s),\bar\a(t;s),\Th(t;s,\bar X(t;s),\bar\a(t;s)),\Th_x(t;s,\bar X(t;s),\bar\a(t;s))\bar\si(t;s),\\ [1mm]
\ns\ds\qq\qq\qq\qq\qq\qq\qq\qq\qq\big[Q(x)\Th(t;s,\bar X(t;s),\cd)\big]_{\bar\a(t;s)},\bar\Psi(t;\bar X(t;s),\bar\a(t;s))\big)ds\\ [1mm]
\ns\ds\q+\Th_x(t;s,\bar X(t;s),\bar\a(t;s))\bar\si(t;s)dW(s)\\ [2mm]
\ns\ds\q+\int_\dbR\big[\Th\big(t;s,\bar X(t;s-),\bar\a(t;s-)+\m(\bar X(t;s-),\bar\a(t;s-),\th)\big)-\Th(t;s,\bar X(t;s-),\bar\a(t;s-))\big]\wt N(ds,d\th)\\ [2mm]
\ns\ds=-g\big(t,s,\bar X(t;s),\bar\a(t;s),\bar Y(t;s),\bar Z(t;s),\int_\dbR\bar\G(t;s,\th)\pi(d\th),\bar\Psi(t;\bar X(t;s),\bar\a(t;s)\big)ds\\ [1mm]
\ns\ds\q+\bar Z(t;s)dW(s)+\int_\dbR\bar\G(t;s,\th)\wt N(ds,d\th).\ea$$
By the uniqueness of adapted solution to the FBSDE, we obtain the representation \eqref{representation}.
\endpf

\ms

We call \rf{TH-1} the {\it representation PDE} of FBSDE \rf{FBSDE1}, because the backward components $(\bar Y(\cd),\bar Z(t;\cd),$ $\bar\G(t;\cd\,,\cd))$ of \rf{FBSDE1} are represented in terms of the forward components $(\bar X(t;\cd)$, $\bar\a(t;\cd))$ through the solution $\Th(t;\cd\,,\cd\,,\cd)$ of \rf{TH-1}. We may also call the solution $\Th(t;\cd\,,\cd\,,\cd)$ of PDE \rf{TH-1} a {\it decoupling field} of our FBSDE (See \cite{Ma-Wu-Zhang-Zhang2015, Zhang2017}).

\section{Approximate Equilibrium Strategy}

Inspired by \cite{Yong2012b, Wei-Yong-Yu2016}, we investigate Problem (N) by looking at a family of approximate problems determined by partitions of the time interval $[0,T]$. More precisely, let $\Pi$ be a partition of $[0,T]$:
$$\Pi:\q0=t_0<t_1<\cds<t_{N-1}<t_N=T,$$
whose {\it mesh size} $\|\Pi\|$ is defined by
$$\|\Pi\|=\max_k|t_k-t_{k-1}|.$$
We imagine that there are $N$ players involved in a multi-person differential game. Player $k$ takes over the system at time $t_{k-1}$, controls the system on $[t_{k-1},t_k)$ and hand it over to Player $(k+1)$ at $t_k$. There are two rules for the players:

\ms

(i) Every player is playing optimally (in a proper sense, see below).

\ms

(ii) Every player (stubbornly) uses his/her own way of discounting for the future cost, although they will not control the system then.

\ms

In this subsection, we construct an {\it equilibrium strategy} for these $N$ players (or associated with the partition $\Pi$), making use of the dynamic programming idea.

\ms

We begin with Player $N$ on $[t_{N-1},T]$. Let $(\t,\xi,\iota)\in\cD[t_{N-1},T]$ and $u^N(\cd)\in\cU[\t,T]$. Consider the state equation \rf{state}, whose solution is denoted by $(X^N(\cd),\a^N(\cd))$ with the recursive cost functional
\bel{J(N)}J^N(\t,\xi,\iota;u^N(\cd))=Y^N(\t),\ee
where $(Y^N(\cd),Z^N(\cd),\G^N(\cd\,,\cd))$ is the adapted solution to the following BSDE on $[\t,T]$:
\bel{BSDE(N)}\left\{\2n\ba{ll}
\ns\ds dY^N(s)=-g\big(t_{N-1},s,X^N(s),\a^N(s),Y^N(s),Z^N(s),
\int_\dbR\G^N(s,\th)\pi(d\th),u^N(s)\big)ds\\
\ns\ds\qq\qq\qq\qq+Z^N(s)dW(s)+\int_\dbR\G^N(s-,\th)\wt N(ds,d\th),\qq s\in[\t,T),\\
\ns\ds Y^N(T)=h(t_{N-1},X^N(T),\a^N(T)).\ea\right.\ee
Note that
\bel{J(N)=J}J^N(t_{N-1},\xi,\iota;u^N(\cd))=J(t_{N-1},\xi,\iota;u^N(\cd)),
\qq\forall(t_{N-1},\xi,\iota)\in\cD,~u^N(\cd)\in\cU[t_{N-1},T].\ee
Player $N$ wants to solve the following problem.

\ms

\bf Problem (C$_N$). \rm For $(\t,\xi,\iota)\in\cD[t_{N-1},T]$, find a
$\bar u^N(\cd)\in\cU[\t,T]$ such that
$$J^N(\t,\xi,\iota;\bar u^N(\cd))=\inf_{u^N(\cd)\in\cU[\t,T]}J^N(\t,\xi,\iota;u^N(\cd))\equiv V^N(\t,\xi,\iota).$$

It is known that under (H1)--(H2), the value function $V^N(\cd\,,\cd\,,\cd)$ is the classical solution to the following HJB equation system:
\bel{HJB(N)}\left\{\2n\ba{ll}
\ns\ds V^N_s(s,x,i)+\inf_{u\in U}\dbH^N(s,x,i,V^N(s,x,\cd),V^N_x(s,x,i),V^N_{xx}(s,x,i),u)=0,\\
\ns\ds\qq\qq\qq\qq\qq\qq\qq\qq\qq(s,x,i)\in[t_{N-1},T]\times\dbR^n\times M,\\
\ns\ds V^N(T,x,i)=h(t_{N-1},x,i),\qq(x,i)\in\dbR^n\times M,\ea\right.\ee
where
$$\ba{ll}
\ns\ds\dbH^N(s,x,i,v,\Bp,\BP,u)={1\over2}\tr\big[\si(s,x,i,u)^\top \BP\si(s,x,i,u)\big]+\Bp b(s,x,i,u)+\big[Q(x)v\big]_i\\
\ns\ds\qq\qq\qq\qq\qq\qq+g\big(t_{N-1},s,x,i,v^i,\Bp\si(s,x,i,u),\big[Q(x)v\big]_i,u\big),\ea$$
with $v=(v^1,\cds,v^m)^\top\in\dbR^n$. Now, under (H3), we may define
\bel{bar Psi(N)}\ba{ll}
\ns\ds\bar\Psi^N(s,x,i)=\psi (t_{N-1};s,x,i,V^N(s,x,\cd),V_x^N(s,x,i),V_{xx}^N(s,x,i)),\\
\ns\ds\qq\qq\qq\qq\qq\qq\qq\qq(s,x,i)\in[t_{N-1},T]\times\dbR^n\times M.\ea\ee
Then \rf{HJB(N)} is equivalent to the following:
\bel{HJB(N)*}\left\{\2n\ba{ll}
\ns\ds V^N_s(s,x,i)+\bar\cA^NV^N(s,x,i)+g(t_{N-1},s,x,i,V^N(s,x,i),V^N_x(s,x,i)
\si(s,x,i,\bar\Psi^N(s,x,i)),\\
\ns\ds\qq\qq\qq\qq[Q(s)V^N(s,x,\cd)]_i,\bar\Psi^N(s,x,i))=0,
\qq(s,x,i)\in[t_{N-1},T]\times\dbR^n\times M,\\
\ns\ds V^N(T,x,i)=h(t_{N-1},x,i),\qq(x,i)\in\dbR^n\times M,\ea\right.\ee
where
\bel{A(N)}\ba{ll}
\ns\ds\bar\cA^NV^N(s,x,i)={1\over2}\tr\[\si\big(s,x,i,\bar\Psi^N(s,x,i)\big)^\top
V^N_{xx}(s,x,i)\si\big(s,x,i,\bar\Psi^N(s,x,i)\big)\]\\
\ns\ds\qq\qq\qq\qq\qq+V^N_x(s,x,i)b(s,x,i,\bar\Psi^N(s,x,i))
+\big[Q(x)V^N(s,x,\cd))\big]_i\,.\ea\ee
Now, for any $(t_{N-1},\xi,\iota)\in\cD[t_{N-1},T]$, let $(\bar X^N(\cd),\bar\a^N(\cd))$ be the solution to the following {\it closed-loop system}:
\bel{closed-(N)}\left\{\2n\ba{ll}
\ns\ds d\bar X^N(s)=b\big(s,\bar X^N(s),\bar\a^N(s),\bar\Psi^N(s,\bar X^N(s),\bar\a^N(s))
\big)ds\\ [2mm]
\ns\ds\qq\qq\qq+\si\big(s,\bar X^N(s),\bar\a^N(s),\bar\Psi^N(s,\bar X^N(s),\bar\a^N(s))\big)dW(s),\qq s\in[t_{N-1},T],\\
\ns\ds d\bar\a^N(s)=\int_\dbR\m(\bar X^N(s),\bar\a^N(s-),\th)N(ds,d\th),\qq s\in[t_{N-1},T],\\
\ns\ds\bar X^N(t_{N-1})=\xi,\q\bar\a^N(t_{N-1})=\iota,\ea\right.\ee
and define
\bel{bar u(N)}\bar u^N(s)=\bar\Psi^N(s,\bar X^N(s),\bar\a^N(s)),\qq s\in[t_{N-1},T].\ee
Then by Theorem \ref{verification theorem},
\bel{V(N)=J(N)}V^N(t_{N-1},\xi,\iota)=J^N(t_{N-1},\xi,\iota,\bar u^N(\cd)).\ee
and $(\bar X^N(\cd),\bar\a^N(\cd),\bar u^N(\cd))$ is the optimal triple of Problem (C$_N$) for the initial triple $(t_{N-1},\xi,\iota)$.

\ms

Next, for Player $(N-1)$ on $[t_{N-2},t_{N-1}]$, let $(\t,\xi,\iota)\in\cD[t_{N-2},t_{N-1}]$ and $u^{N-1}(\cd)\in\cU[\t,t_{N-1}]$. Consider the state equation \rf{state} on $[\t,t_{N-1}]$, whose solution is denoted by $(X^{N-1}(\cd),\a^{N-1}(\cd))$, with the recursive cost functional
\bel{J(N-1)}J^{N-1}(\t,\xi,\iota;u^{N-1}(\cd))=Y^{N-1}(\t),\ee
where $(Y^{N-1}(\cd),Z^{N-1}(\cd),\G^{N-1}(\cd\,,\cd))$ is the adapted solution to the following BSDE on $[\t,T]$:
\bel{BSDE(N-1)}\2n\left\{\2n\ba{ll}
\ns\ds dY^{N-1}(s)=-g\big(t_{N-2},s,\bar X^N(s),\bar\a^N(s),Y^{N-1}(s),Z^{N-1}(s),\\
\ns\ds\qq\qq\qq\qq\qq\qq\qq\qq\int_\dbR\G^{N-1}(s;\th)\pi(d\th),\bar\Psi^N(s,\bar X^N(s),\bar\a^N(s))\big)ds\\
\ns\ds\qq\qq\qq\qq+Z^{N-1}(s)dW(s)+\int_\dbR\G^{N-1}(s-,\th)\wt N(ds,d\th),\qq s\in[t_{N-1},T),\\
\ns\ds dY^{N\1n-\1n1}(s)\1n=\1n-g\big(t_{N\1n-\1n2},s,X^{N\1n-\1n1}(s),\a^{N\1n-\1n1}(s),Y^{N\1n-\1n1}(s),
Z^{N\1n-\1n1}(s),\1n\int_\dbR\2n\G^{N\1n-\1n1}(s;\th)\pi(d\th),u^{N\1n-\1n1}(\cd)\big)ds\\
\ns\ds\qq\qq\qq\qq+Z^{N-1}(s)dW(s)+\int_\dbR\G^{N-1}(s-,\th)\wt N(ds,d\th),\qq s\in[t_{N-2},t_{N-1}),\\
\ns\ds Y^{N-1}(T)=h(t_{N-2},\bar X^N(T),\bar\a^N(T)),\q Y^{N-1}(t_{N-1}-)=Y^{N-1}(t_{N-1}),\ea\right.\ee
where $(\bar X^N(\cd),\bar\a^N(\cd),\bar u^N(\cd))$ is the solution to the closed-loop system \rf{closed-(N)} (on $[t_{N-1},T]$) with the initial condition at $t_{N-1}$ given by
$$\bar X^N(t_{N-1})=X^{N-1}(t_{N-1}),\q\bar\a^N(t_{N-1})=\a^{N-1}(t_{N-1}).$$
The cost functional \rf{J(N-1)} is designed based on the two rules mentioned above: Player $(N-1)$ assumes
that Player $N$ will play optimally, based on the initial triple $(t_{N-1},X^{N-1}(t_{N-1}),
\a^{N-1}(t_{N-1}))$, which leads to the dependence of $(Y^{N-1}(\cd),Z^{N-1}(\cd),\G^{N-1}(\cd\,,\cd))$ on
$(\bar X^N(\cd),\bar\a^N(\cd))$ on the interval $[t_{N-1},T]$; Player $(N-1)$ insists to discount in his/her own way, which leads to that $t_{N-2}$ appears in the generator $g(\cd)$ of the above BSDE (instead of $t_{N-1}$). See \cite{Yong2012b} for more details. The problem that Player $(N-1)$ wants to solve is the following:

\ms

\bf Problem (C$_{N-1}$). \rm For $(\t,\xi,\iota)\in\cD^p[t_{N-2},t_{N-1})$, find a
$\bar u^{N-1}(\cd)\in\cU[\t,t_{N-1}]$ such that
$$J^{N-1}(\t,\xi,\iota;\bar u^{N-1}(\cd))=\inf_{u^{N-1}(\cd)\in\cU[\t,t_{N-1}]}J^{N-1}(\t,\xi,\iota;u^{N-1}(\cd))\equiv V^{N-1}(\t,\xi,\iota).$$

\ms

We point out that similar to \rf{J(N)=J}, the following holds:
\bel{J(N-1)=J}\ba{ll}
\ns\ds J^{N-1}\big(t_{N-2},\xi,\iota;u^{N-1}(\cd)\big)=J\big(t_{N-2},\xi,\iota;u^{N-1}(\cd)\oplus \bar u^N(\cd)\big),\\
\ns\ds\qq\qq\qq\qq\forall(t_{N-2},\xi,\iota)\in\cD,~u^{N-1}(\cd)\in\cU[t_{N-2},t_{N-1}],\ea\ee
where
$$u^{N-1}(\cd)\oplus\bar u^N(\cd)=u^{N-1}(\cd)I_{[t_{N-2},t_{N-1})}(\cd)+\bar u^N(\cd)I_{[t_{N-1},T]}(\cd).$$

Note that Problem (C$_{N-1}$) is not classical in the following sense. Player $(N-1)$ takes control $u^{N-1}(\cd)\in\cU[t_{N-2},t_{N-1}]$ controlling the system on $[t_{N-2},t_{N-1}]$. But, the cost functional $J^{N-1}(\t,\xi\iota;u^{N-1}(\cd))$ involves a BSDE on $[t_{N-1},T]$, depending on $(\bar X^N(\cd),\bar\a^N(\cd))$. We now make a proper manipulation so that the cost functional becomes a classical one. To this end, let us look at the following decoupled FBSDE on $[t_{N-1},T]$:
\bel{FBSDE(N)}\left\{\2n\ba{ll}
\ns\ds d\bar X^N(s)=b\big(s,\bar X^N(s),\bar\a^N(s),\bar\Psi^N(s,\bar X^N(s),\bar\a^N(s))
\big)ds\\ [2mm]
\ns\ds\qq\qq\qq+\si\big(s,\bar X^N(s),\bar\a^N(s),\bar\Psi^N(s,\bar X^N(s),\bar\a^N(s))\big)dW(s),\\
\ns\ds d\bar\a^N(s)=\int_\dbR\m(\bar X^N(s),\bar\a^N(s-),\th)N(ds,d\th),\\
\ns\ds dY^{N-1}(s)=-g\big(t_{N-2},s,\bar X^N(s),\bar\a^N(s),Y^{N-1}(s),Z^{N-1}(s),\int_\dbR\G^{N-1}(s;\th)\pi(d\th),\\
\ns\ds\qq\qq\qq\qq\qq\qq\qq\qq\qq\qq\qq\qq\qq\bar\Psi^N(s,\bar X^N(s),\bar\a^N(s))\big)ds\\
\ns\ds\qq\qq\qq+Z^{N-1}(s)dW(s)+\int_\dbR\G^{N-1}(s-,\th)\wt N(ds,d\th),\qq s\in[t_{N-1},T),\\
\ns\ds\bar X^N(t_{N-1})=\xi,\q\bar\a^N(t_{N-1})=\iota,\q Y^{N-1}(T)=h(t_{N-2},\bar X^N(t_N),\bar\a^N(t_N)).\ea\right.\ee
Let $\Th^{N-1}(\cd\,,\cd\,,\cd)$ be the solution to the following representation PDE:
\bel{TH(N-1)}\left\{\2n\ba{ll}
\ns\ds\Th^{N-1}_s(s,x,i)+\bar\cA^N\Th^{N-1}(s,x,i)\\
\ns\ds\q+g\big(t_{N-2},s,x,i,\Th^{N-1}(s,x,i),\Th^{N-1}_x(s,x,i)\si(s,x,i,\bar\Psi^N(s,x,i)),\\
\ns\ds\qq\qq\qq\qq\qq\q[Q(x)\Th^{N-1}(s,x,\cd)]_i,\bar\Psi^N(s,x,i)\big)=0;\qq(s,x,i)\1n\in\1n[t_{N-1},T]\1n
\times\1n\dbR^n\1n\times\1n M,\\
\ns\ds\Th^{N-1}(T,x,i)=h(t_{N-2},x,i),\qq(x,i)\in\dbR^n\times M,\ea\right.\ee
where $\bar\cA^N$ is defined by \rf{A(N)}. Then, by Theorem \ref{representation}, the following representation hold:
\bel{Y=Th(N-1)}\left\{\2n\ba{ll}
\ns\ds Y^{N-1}(s)=\Th^{N-1}(s,\bar X^N(s),\bar\a^N(s)),\\
\ns\ds Z^{N-1}(s)=\Th^{N-1}_x(s,\bar X^N(s),\bar\a^N(s))\si\big(s,\bar X^N(s),\bar\a^N(s),\bar\Psi^N(s,\bar X^N(s),\bar\a^N(s))\big),\\
\ns\ds\bar\G^{N-1}(s,\th)=\Th^{N-1}(s,\bar X^N(s),\bar\a^N(s)+\m(\bar X^N(s),\bar\a^N(s),\th))-\Th^{N-1}(s,\bar X^N(s),\bar\a^N(s)),\\
\ns\ds\qq\qq\qq\qq\qq\qq\qq\qq\qq\qq\qq\qq\qq\qq\qq\q s\in[t_{N-1},T].\ea\right.\ee
Thus, on $[t_{N-2},t_{N-1}]$, we have
\bel{BSDE(N-1)*}\left\{\2n\ba{ll}
\ns\ds dY^{N\1n-\1n1}(s)\1n=\1n-g\big(t_{N\1n-\1n2},s,X^{N\1n-\1n1}(s),\a^{N\1n-\1n1}(s),Y^{N\1n-\1n1}(s),Z^{N\1n-\1n1}(s),
\1n\int_\dbR\2n\G^{N\1n-\1n1}(s;\th)\pi(d\th),u(s)\big)ds\\
\ns\ds\qq\qq\qq\qq+Z^{N-1}(s)dW(s)+\int_\dbR\G^{N-1}(s-,\th)\wt N(ds,d\th),\qq s\in[t_{N-2},t_{N-1}),\\
\ns\ds Y^{N-1}(t_{N-1})=\Th^{N-1}\big(t_{N-1},X^{N-1}(t_{N-1}),\a^{N-1}(t_{N-1})\big).\ea\right.\ee
With the above manipulation, we see that the cost functional $J^{N-1}(\t,\xi,i;u^{N-1}(\cd))$ defined by \rf{J(N-1)} is completely determined by the above BSDE (on $[t_{N-2},t_{N-1}]$). Then Problem (C$_{N-1}$) becomes a classical optimal control problem for a regime switching system with recursive cost functional. We call
$J^{N-1}(\t,\xi,\iota;u^{N-1}(\cd))$ defined by \rf{J(N-1)}, \rf{TH(N-1)} and \rf{BSDE(N-1)*} a {\it sophisticated recursive cost functional}.

\ms

Similar to Problem (C$_N$), under (H1)--(H2), the value function $V^{N-1}(\cd\,,\cd\,,\cd)$ of Problem (C$_{N-1}$) is the classical solution of the following HJB equation:
\bel{HJB(N-1)}\left\{\2n\ba{ll}
\ns\ds V^{N-1}_s(s,x,i)+\inf_{u\in U}\dbH^{N-1}(s,x,i,V^{N-1}(s,x,\cd),V^{N-1}_x(s,x,i),V^{N-1}_{xx}(s,x,i),u)=0,\\
\ns\ds\qq\qq\qq\qq\qq\qq\qq\qq\qq(s,x,i)\in[t_{N-2},t_{N-1}]\times\dbR^n\times M,\\
\ns\ds V^{N-1}(t_{N-1},x,i)=\Th^{N-1}(t_{N-1},x,i),\qq(x,i)\in\dbR^n\times M,\ea\right.\ee
where
$$\ba{ll}
\ns\ds\dbH^{N-1}(s,x,i,v,\Bp,\BP,u)=\Bp b(s,x,i,u)+{1\over2}\tr\big[\si(s,x,i,u)^\top \BP\si(s,x,i,u)\big]+\big[Q(x)v\big]_i\\
\ns\ds\qq\qq\qq\qq\qq\qq+g\big(t_{N-2},s,x,i,v^i,\Bp\si(s,x,i),\big[Q(x)v\big]_i,u\big),\ea$$
with $v=(v^1,\cds,v^m)^\top\in\dbR^n$. Further, under (H3), we may define
\bel{bar Psi(N-1)}\ba{ll}
\ns\ds\bar\Psi^{N-1}(s,x,i)=\psi (t_{N-2};s,x,i,V^{N-1}(s,x,\cd),V_x^{N-1}(s,x,i),V_{xx}^{N-1}(s,x,i)),\\
\ns\ds\qq\qq\qq\qq\qq\qq\qq\qq(s,x,i)\in[t_{N-2},t_{N-1}]\times\dbR^n\times M,\ea\ee
Then \rf{HJB(N-1)} is equivalent to the following:
\bel{HJB(N-1)*}\left\{\2n\ba{ll}
\ns\ds V^{N-1}_s(s,x,i)+\bar\cA^{N-1}V^{N-1}(s,x,i)+g\big(t_{N-2},s,x,i,V^{N-1}(s,x,i),\\
\ns\ds\qq\qq V^{N-1}_x(s,x,i)\si(s,x,i,\bar\Psi^{N-1}(s,x,i)),[Q(s)V^{N-1}(s,x,\cd)]_i,\bar\Psi^{N-1}(s,x,i)\big)=0,\\
\ns\ds\qq\qq\qq\qq\qq\qq\qq\qq\qq\qq\qq(s,x,i)\in[t_{N-2},t_{N-1}]\times\dbR^n\times M,\\
\ns\ds V^{N-1}(t_{N-1},x,i)=\Th^{N-1}(t_{N-1},x,i),\qq(x,i)\in\dbR^n\times M,\ea\right.\ee
where
\bel{A(N-1)}\ba{ll}
\ns\ds\bar\cA^{N-1}V^{N-1}(s,x,i)={1\over2}\tr\[\si\big(s,x,i,\bar\Psi^{N-1}(s,x,i)\big)^\top
V^{N-1}_{xx}(s,x,i)\si\big(s,x,i,\bar\Psi^{N-1}(s,x,i)\big)\]\\
\ns\ds\qq\qq\qq\qq\qq+V^{N-1}_x(s,x,i)b(s,x,i,\bar\Psi^{N-1}(s,x,i))
+\big[Q(x)V^{N-1}(s,x,\cd))\big]_i\,.\ea\ee
Now, for any $(t_{N-2},\xi,\iota)\in\cD[t_{N-2},t_{N-1})$, let $(\bar X^{N-1}(\cd),\bar\a^{N-1}(\cd))$ be the solution to the following closed-loop system on $[t_{N-2},t_{N-1}]$:
\bel{closed-(N-1)}\left\{\2n\ba{ll}
\ns\ds d\bar X^{N-1}(s)=b\big(s,\bar X^{N-1}(s),\bar\a^{N-1}(s),\bar\Psi^{N-1}(s,\bar X^{N-1}(s),\bar\a^{N-1}(s))\big)ds\\ [2mm]
\ns\ds\qq\qq\qq+\si\big(s,\bar X^{N-1}(s),\bar\a^{N-1}(s),\bar\Psi^{N-1}(s,\bar X^{N-1}(s),\bar\a^{N-1}(s))\big)dW(s),\\
\ns\ds d\bar\a^{N-1}(s)=\int_\dbR\m(\bar X^{N-1}(s),\bar\a^{N-1}(s-),\th)N(ds,d\th),\\
\ns\ds\bar X^{N-1}(t_{N-2})=\xi,\q\bar\a^{N-1}(t_{N-2})=\iota,\ea\right.\ee
and define
\bel{bar u(N-1)}\bar u^{N-1}(s)=\bar\Psi^{N-1}(s,\bar X^{N-1}(s),\bar\a^{N-1}(s)),\qq s\in[t_{N-2},t_{N-1}].\ee
Then by Theorem \ref{verification theorem},
\bel{V(N-1)=J(N-1)}V^{N-1}(t_{N-2},\xi,\iota)=J^{N-1}(t_{N-2},\xi,\iota,\bar u^{N-1}(\cd)).\ee
and $(\bar X^{N-1}(\cd),\bar\a^{N-1}(\cd),\bar u^{N-1}(\cd))$ is an optimal triple of Problem (C$_{N-1}$) (with the sophisticated recursive cost functional) for initial triple $(t_{N-2},\xi,\iota)$.

\ms

In general, for Player $k$ on $[t_{k-1},t_k)$, let $(\t,\xi,\iota)\in\cD[t_{k-1},t_k)$ and $u^k(\cd)\in\cU[\t,t_k]$. Consider the state equation \rf{state} on $[\t,t_k]$, whose solution is denoted by $(X^k(\cd),\a^k(\cd))$, with the recursive cost functional
\bel{J(k)}J^k(\t,\xi,\iota;u^k(\cd))=Y^k(\t),\ee
where $(Y^k(\cd),Z^k(\cd),\G^k(\cd\,,\cd))$ is the adapted solution to the following BSDE on $[t_{k-1},T]$:
\bel{BSDE(k)}\left\{\2n\ba{ll}
\ns\ds dY^k(s)=-g\big(t_{k-1},s,\bar X^\Pi(s),\bar\a^\Pi(s),Y^k(s),Z^k(s),\int_\dbR\G^k(s;\th)\pi(d\th),\bar\Psi^\Pi(s,\bar X^\Pi(s),\bar\a^\Pi(s))\big)ds\\
\ns\ds\qq\qq\qq\qq+Z^k(s)dW(s)+\int_\dbR\G^k(s-,\th)\wt N(ds,d\th),\qq s\in[t_k,T),\\
\ns\ds dY^k(s)\1n=\1n-g\big(t_{k-1},s,X^k(s),\a^k(s),Y^k(s),Z^k(s),\1n\int_\dbR\2n\G^k(s;\th)\pi(d\th),
u^k(\cd)\big)ds\\
\ns\ds\qq\qq\qq\qq+Z^k(s)dW(s)+\int_\dbR\G^k(s-,\th)\wt N(ds,d\th),\qq s\in[t_{k-1},t_k),\\
\ns\ds Y^k(T)=h(t_{k-1},\bar X^\Pi(T),\bar\a^\Pi(T)),\q Y^k(t_\ell-)=Y^k(t_\ell),\q k\les\ell\les N-1.\ea\right.\ee
In the above, $(\bar X^\Pi(\cd),\bar\a^\Pi(\cd))$ is the solution to the following closed-loop system on
$[t_k,T]$:
\bel{closed-loop(Pi)}\left\{\2n\ba{ll}
\ns\ds d\bar X^\Pi(s)=b\big(s,\bar X^\Pi(s),\bar\a^\Pi(s),\bar\Psi^\Pi(s,\bar X^\Pi(s),\bar\a^\Pi(s))\big)ds\\
\ns\ds\qq\qq\qq+\si\big(s,\bar X^\Pi(s),\bar\a^\Pi(s),\bar\Psi^\Pi(s,\bar X^\Pi(s),\bar\a^\Pi(s))\big)dW(s),\q
s\in[t_k,T],\\
\ns\ds d\bar\a^\Pi(s)=\int_\dbR\m(s,\bar X^\Pi(s),\bar\a^\Pi(s-),\th)N(ds,d\th),\qq s\in[t_k,T],\\
\ns\ds\bar X^\Pi(t_k)=X^k(t_k),\q\bar\a^\Pi(t_k)=\a^k(t_k),\ea\right.\ee
where
\bel{Psi(Pi)*}\ba{ll}
\ns\ds\bar\Psi^\Pi(s,x,i)=\psi\big(t^\Pi(s);s,x,i,V^\Pi(s,x,\cd),V^\Pi_x(s,x,i),V^\Pi_{xx}(s,x,i)\big),\\
\ns\ds\qq\qq\qq\qq\qq\qq\qq\qq\qq(s,x,i)\in[t_k,T]\times\dbR^n\times M,\ea\ee
with
\bel{V(Pi)}V^\Pi(s,x,i)=\sum_{j=k+1}^{N-1}V^j(s,x,i)I_{[t_{j-1},t_j)}(s)+V^N(s,x,i)I_{[t_{N-1},T]}(s),\q s\in[t_k,T],\ee
and
\bel{t(Pi)}t^\Pi(s)=\sum_{j=1}^{N-1}t_{j-1}I_{[t_{j-1},t_j)}(s)+t_{N-1}I_{[t_{N-1},T]}(s),\qq s\in[0,T].\ee
Note that $V^j(\cd\,,\cd\,,\cd)$ is the classical solution to the corresponding HJB equation system on $[t_{j-1},t_j]$ (see \rf{HJB(N)} and \rf{HJB(N-1)} for $V^N(\cd\,,\cd\,,\cd)$ and $V^{N-1}(\cd\,,\cd\,,\cd)$, respectively; and see \rf{HJB(k)} below for $V^k(\cd\,,\cd\,,\cd)$).
Player $k$ tries to solve the following problem.

\ms

\bf Problem (C$_k$). \rm For $(\t,\xi,\iota)\in\cD[t_{k-1},t_k)$, find a
$\bar u^k(\cd)\in\cU[\t,t_k]$ such that
$$J^k(\t,\xi,\iota;\bar u^k(\cd))=\inf_{u^k(\cd)\in\cU[\t,t_k]}J^k(\t,\xi,\iota;u^k(\cd))\equiv V^k(\t,\xi,\iota).$$

\ms

Keep in mind that similar to \rf{J(N)=J} and \rf{J(N-1)=J}, we have
\bel{J(k)=J}J^k\big(t_{k-1},\xi,\iota;u^k(\cd)\big)=J(t_{k-1},\xi,\iota;u^k(\cd)\oplus
\bar\Psi^\Pi(\cd)\big),\qq(t_{k-1},\xi,\iota)\in\cD,~u^k(\cd)\in\cU[t_{k-1},t_k],\ee
where
$$u^k(\cd)\oplus\bar\Psi^\Pi(\cd)=u^k(\cd)I_{[t_{k-1},t_k)}(\cd)+\bar \Psi^\Pi(\cd\,,\bar X^\Pi(\cd),\bar\a^\Pi(\cd))I_{[t_k,T]}(\cd).$$
Similar to Problem (C$_{N-1}$), we need to get a sophisticated recursive cost functional for Player $k$. To this end, we look at the following decoupled FBSDE on $[t_k,T]$:
\bel{FBSDE(k)}\left\{\2n\ba{ll}
\ns\ds d\bar X^\Pi(s)=b\big(s,\bar X^\Pi(s),\bar\a^\Pi(s),\bar\Psi^\Pi(s,\bar X^\Pi(s),\bar\a^\Pi(s))\big)ds\\ [2mm]
\ns\ds\qq\qq\qq+\si\big(s,\bar X^\Pi(s),\bar\a^\Pi(s),\bar\Psi^\Pi(s,\bar X^\Pi(s),\bar\a^\Pi(s))\big)dW(s),\q s\in[t_k,T],\\
\ns\ds d\bar\a^\Pi(s)=\int_\dbR\m(\bar X^\Pi(s),\bar\a^\Pi(s-),\th)N(ds,d\th),\q s\in[t_k,T],\\
\ns\ds dY^k(s)=-g\big(t_{k-1},s,\bar X^\Pi(s),\bar\a^\Pi(s),Y^k(s),Z^k(s),\int_\dbR\G^k(s;\th)\pi(d\th),\bar\Psi^\Pi(s,\bar X^\Pi(s),\bar\a^\Pi(s))\big)ds\\
\ns\ds\qq\qq\qq\qq+Z^k(s)dW(s)+\int_\dbR\G^k(s-,\th)\wt N(ds,d\th),\qq s\in[t_k,T),\\
\ns\ds\bar X^\Pi(t_k)=\xi,\q\bar\a^\Pi(t_k)=\iota,\q Y^k(T)=h(t_{k-1},\bar X^\Pi(T),\bar\a^\Pi(T)).\ea\right.\ee
Let $\Th^k(\cd\,,\cd\,,\cd)$ be the solution to the following representation PDE:
\bel{TH(k)}\left\{\2n\ba{ll}
\ns\ds\Th^k_s(s,x,i)+\bar\cA^\Pi\Th^k(s,x,i)+g\big(t_{k-1},s,x,i,\Th^k(s,x,i),\Th^k_x(s,x,i)\si(s,x,i,\bar\Psi^\Pi(s,x,i)),\\
\ns\ds\qq\qq\qq\qq\qq\qq\qq\qq\qq[Q(x)\Th^k(s,x,\cd)]_i,\bar\Psi^\Pi(s,x,i)\big)=0;\\
\ns\ds\qq\qq\qq\qq\qq\qq\qq\qq\qq(s,x,i)\in[t_k,T]\times\dbR^n\times M,\\
\ns\ds\Th^k(T,x,i)=h(t_{k-1},x,i),\qq(x,i)\in\dbR^n\times M,\ea\right.\ee
where $\bar\cA^\Pi$ is defined by the following:
\bel{cA(Pi)}\ba{ll}
\ns\ds\bar\cA^\Pi\Th^k(s,x,i)={1\over2}\tr\[\si\big(s,x,i,\bar\Psi^\Pi(s,x,i)\big)^\top
\Th^k_{xx}(s,x,i)\si\big(s,x,i,\bar\Psi^\Pi(s,x,i)\big)\]\\
\ns\ds\qq\qq\qq\qq\qq+\Th^k_x(s,x,i)b(s,x,i,\bar\Psi^\Pi(s,x,i))+\big[Q(x)\Th^k(s,x,\cd))\big]_i\,.\ea\ee
Then by Theorem \ref{representation}, the following representation hold:
\bel{Y=Th(k)}\left\{\2n\ba{ll}
\ns\ds Y^k(s)=\Th^k(s,\bar X^\Pi(s),\bar\a^\Pi(s)),\\
\ns\ds Z^k(s)=\Th^k_x(s,\bar X^\Pi(s),\bar\a^\Pi(s))\si(s,\bar X^\Pi(s),\bar\a^\Pi(s),\bar\Psi^\Pi(s,\bar X^\Pi(s),\bar\a^\Pi(s)),\\
\ns\ds\bar\G^k(s,\th)\1n=\1n\Th^k(s,\bar X^\Pi(s),\bar\a^\Pi(s)\1n+\1n\m(\bar X^\Pi(s),\bar\a^\Pi(s),\th))\1n-\1n\Th^k(s,\bar X^\Pi(s),\bar\a^\Pi(s)),\ea\right.\q s\in[t_k,T].\ee
Thus, on $[t_{k-1},t_k]$, we have
\bel{BSDE(k)*}\left\{\2n\ba{ll}
\ns\ds dY^k(s)\1n=\1n-g(t_{k-1},s,X^k(s),\a^k(s),Y^k(s),Z^k(s),\1n\int_\dbR\2n\G^k(s;\th)\pi(d\th),u(s))ds\\
\ns\ds\qq\qq\qq\qq+Z^k(s)dW(s)+\int_\dbR\G^k(s-,\th)\wt N(ds,d\th),\qq s\in[t_{k-1},t_k),\\
\ns\ds Y^k(t_k)=\Th^k\big(t_k,X^k(t_k),\a^k(t_k)\big).\ea\right.\ee
Consequently, Player $k$ has his/her sophisticated recursive cost functional \rf{J(k)} with \rf{BSDE(k)*}, and Problem (C$_k$) becomes a classical optimal control problem on $[t_{k-1},t_k]$.

\ms

Now, to solve Problem (C$_k$), we let $V^k(\cd\,,\cd\,,\cd)$ be the classical solution to the following HJB equation:
\bel{HJB(k)}\left\{\2n\ba{ll}
\ns\ds V^k_s(s,x,i)+\inf_{u\in U}\dbH^k(s,x,i,V^k(s,x,\cd),V^k_x(s,x,i),V^k_{xx}(s,x,i),u)=0,\\
\ns\ds\qq\qq\qq\qq\qq\qq\qq\qq\qq(s,x,i)\in[t_{k-1},t_k]\times\dbR^n\times M,\\
\ns\ds V^k(t_k,x,i)=\Th^k(t_k,x,i),\qq(x,i)\in\dbR^n\times M,\ea\right.\ee
where
$$\ba{ll}
\ns\ds\dbH^k(s,x,i,v,\Bp,\BP,u)=\Bp b(s,x,i,u)+{1\over2}\tr\big[\si(s,x,i,u)^\top \BP\si(s,x,i,u)\big]+\big[Q(x)v\big]_i\\
\ns\ds\qq\qq\qq\qq\qq+g\big(t_{k-1},s,x,i,v^i,\Bp\si(s,x,i,u),\big[Q(x)v\big]_i,
u\big),\ea$$
with $v=(v^1,\cds,v^m)^\top\in\dbR^n$. Then we define
\bel{bar Psi(N-2)}\ba{ll}
\ns\ds\bar\Psi^k(s,x,i)=\psi(t_{k-1};s,x,i,V^k(s,x,\cd),V_x^k(s,x,i),V_{xx}^k(s,x,i)),
\qq(s,x,i)\in[t_{k-1},t_k]\times\dbR^n\times M.\ea\ee
With the above, \rf{HJB(k)} can be rewritten as follows:
\bel{HJB(k)*}\left\{\2n\ba{ll}
\ns\ds V^k_s(s,x,i)+\bar\cA^kV^k(s,x,i)+g(t_{k-1},s,x,i,V^k(s,x,i),V^k_x(s,x,i)
\si(s,x,i,\bar\Psi^k(s,x,i)),\\
\ns\ds\qq\qq\qq\qq[Q(s)V^k(s,x,\cd)]_i,\bar\Psi^k(s,x,i))=0,
\qq(s,x,i)\in[t_{k-1},t_k]\times\dbR^n\times M,\\
\ns\ds V^k(t_k,x,i)=\Th^k(t_k,x,i),\qq(x,i)\in\dbR^n\times M,\ea\right.\ee
where $\bar\cA^\Pi$ is defined by \rf{cA(Pi)}. Now, for any $(t_{k-1},\xi,\iota)\in\cD[t_{k-1},t_k)$, let $(\bar X^k(\cd),\bar\a^k(\cd))$ be the solution to the following closed-loop system on $[t_{k-1},t_k]$:
\bel{closed-(k)}\left\{\2n\ba{ll}
\ns\ds d\bar X^k(s)=b\big(s,\bar X^k(s),\bar\a^k(s),\bar\Psi^k(s,\bar X^k(s),\bar\a^k(s))\big)ds\\ [2mm]
\ns\ds\qq\qq\qq+\si\big(s,\bar X^k(s),\bar\a^k(s),\bar\Psi^k(s,\bar X^k(s),\bar\a^k(s))\big)dW(s),\\
\ns\ds d\bar\a^k(s)=\int_\dbR\m(\bar X^k(s),\bar\a^k(s-),\th)N(ds,d\th),\\
\ns\ds\bar X^k(t_{k-1})=\xi,\q\bar\a^k(t_{k-1})=\iota.\ea\right.\ee
Let
\bel{bar u(k)}\bar u^k(s)=\bar\Psi^k(s,\bar X^k(s),\bar\a^k(s)),\qq s\in[t_{k-1},t_k].\ee
Then, by Theorem \ref{verification theorem},
\bel{V(k)=J(k)}V^k(t_{k-1},\xi,\iota)=J^k(t_{k-1},\xi,\iota;\bar u^k(\cd))=\inf_{u^k(\cd)\in\cU[t_{k-1},t_k]} J^k(t_{k-1},\xi,\iota;u^k(\cd)).\ee
and $(\bar X^k(\cd),\bar\a^k(\cd),\bar u^k(\cd))$ is an optimal triple of Problem (C$_k$) for the initial triple $(t_{k-1},\xi,\iota)$. We refer to \rf{V(k)=J(k)} as the {\it local optimality} of $\bar u^k(\cd)$ or $\bar\Psi^k(\cd\,,\cd\,,\cd)$.

\ms

By induction, we can construct the sequences
$$\left\{\2n\ba{ll}
\ns\ds V^k(s,x,i),\qq(s,x,i)\in[t_{k-1},t_k)\times\dbR^n\times M,\q k=N,N-1,\cds,2,1,\\
\ns\ds\bar\Psi^k(s,x,i),\qq(s,x,i)\in[t_{k-1},t_k)\times\dbR^n\times M,\q k=N,N-1,\cds,2,1,\\
\ns\ds\Th^k(s,x,i),\qq(s,x,i)\in[t_k,T]\times\dbR^n\times M,\q k=N-1,N-2,\cds,2,1,0.\ea\right.$$
We denote
$$\left\{\2n\ba{ll}
\ns\ds V^\Pi(s,x,i)=\sum_{k=1}^{N-1}V^k(s,x,i)I_{[t_{k-1},t_k)}(s)+V^N(s,x,i)I_{[t_{N-1},T]}(s),\qq(s,x,i)\in[0,T)
\times\dbR^n\times M,\\
\ns\ds\bar\Psi^\Pi(s,x,i)=\sum_{k=1}^{N-1}\bar\Psi^k(s,x,i)I_{[t_{k-1},t_k)}(s)+\bar\Psi^N(s,x,i)I_{[t_{N-1},T]}(s)\\
\ns\ds\qq\qq~=\psi\big(\t^\Pi(s),x,i,V^\Pi(s,x,\cd),V^\Pi_x(s,x,i),V^\Pi_{xx}(s,x,i)\big),\q(s,x,i)\in[0,T]
\times\dbR^n\times M,\\
\ns\ds\bar\cA^\Pi\Th(s,x,i)={1\over2}\tr\[\si(s,x,i,\bar\Psi^\Pi(s,x,i))^\top
\Th_{xx}(s,x,i)\si(s,x,i,\bar\Psi^\Pi(s,x,i))\]\\ [2mm]
\ns\ds\qq\qq\qq\qq\qq+\Th_x(s,x,i)b(s,x,i,\bar\Psi^\Pi(s,x,i))
+\big[Q(x)\Th(s,x,\cd))\big]_i\,,\\
\ns\ds\qq\qq\qq\qq\qq\qq\qq\qq\qq\qq\qq(s,x,i)\in[t_{k-1},T]\times\dbR^n\times M.\ea\right.$$
The recursive construction for the above sequences can be summarized as follows:

\ms

\it Cycle $N$. \rm

\ms

$\bullet$ Define $\Th^N(t_N,\cd\,,\cd)=h(t_{N-1},\cd\,,\cd)$.

\ms

$\bullet$ Solve HJB equation system on $[t_{N-1},T]$ with terminal condition $\Th^N(t_N,\cd\,,\cd)$ to get $V^N(\cd\,,\cd\,,\cd)$.

\ms

$\bullet$ Define $\bar\Psi^\Pi(\cd\,,\cd\,,\cd)$ on $[t_{N-1},T]$, using $V^N(\cd\,,\cd\,,\cd)$,

\ms

\it Cycle $(N-1)$. \rm

\ms

$\bullet$ Solve representation PDE on $[t_{N-1},T]$ to get $\Th^{N-1}(\cd\,,\cd\,,\cd)$, using $\bar\Psi^\Pi(\cd\,,\cd\,,\cd)$.

\ms

$\bullet$ Solve HJB equation system on $[t_{N-2},t_{N-1})$ with terminal condition $\Th^{N-1}(t_{N-1},\cd\,,\cd)$ to get $V^{N-1}(\cd\,,\cd\,,\cd)$.

\ms

$\bullet$ Extend $\bar\Psi^\Pi(\cd\,,\cd\,,\cd)$ to $[t_{N-2},T]$, using $V^\Pi(\cd\,,\cd\,,\cd)$ (a concatenation of $V^{N-1}(\cd\,,\cd\,,\cd)$ and $V^N(\cd\,,\cd\,,\cd)$).

\ms

$\cds\cds$

\ms

\it Cycle $k$. \rm

\ms

$\bullet$ Solve representation PDE on $[t_k,T]$ to get $\Th^k(\cd\,,\cd\,,\cd)$, using $\bar\Psi^\Pi(\cd\,,\cd\,,\cd)$.

\ms

$\bullet$ Solve HJB equation system on $[t_{k-1},t_k)$ with terminal condition $\Th^k(t_k,\cd\,,\cd)$ to get $V^k(\cd\,,\cd\,,\cd)$.

\ms

$\bullet$ Extend $\bar\Psi^\Pi(\cd\,,\cd\,,\cd)$ to $[t_{k-1},T]$, using $V^\Pi(\cd\,,\cd\,,\cd)$ (a concatenation of $V^k(\cd\,,\cd\,,\cd),\cds,V^N(\cd\,,\cd\,,\cd)$).

\bs

The above procedure will end at the end of {\it Cycle 1}, which completes the construction of the sequences.

\ms

Having constructed $\bar\Psi^\Pi(\cd\,,\cd\,,\cd)$ on $[0,T]\times\dbR^n\times M$, for given $(x,i)\in\dbR^n\times M$, we solve the following closed-loop system on $[0,T]$:
\bel{closed-loop(Pi)}\left\{\2n\ba{ll}
\ns\ds d\bar X^\Pi(s)=b\big(s,\bar X^\Pi(s),\bar\a^\Pi(s),\bar\Psi^\Pi(s,\bar X^\Pi(s),\bar\a^\Pi(s))\big)ds\\ [2mm]
\ns\ds\qq\qq\qq+\si\big(s,\bar X^\Pi(s),\bar\a^\Pi(s),\bar\Psi^\Pi(s,\bar X^\Pi(s),\bar\a^\Pi(s))\big)dW(s),\qq s\in[0,T],\\
\ns\ds\bar X^\Pi(0)=x,\qq\bar\a^\Pi(0)=i.\ea\right.\ee
The local optimality of $\bar\Psi^\Pi(\cd\,,\cd\,,\cd)$ reads
$$\ba{ll}
\ns\ds\bar Y^k(t_{k-1})=J^k(t_{k-1},\bar X^\Pi(t_{k-1}),\bar\a^\Pi(t_{k-1});\bar\Psi^\Pi(\cd\,,\cd\,,\cd))
=V^\Pi\big(t_{k-1},\bar X^\Pi(t_{k-1}),\bar\a^\Pi(t_{k-1})\big)\\
\ns\ds\qq\qq\les J^k(t_{k-1},\bar X^\Pi(t_{k-1}),\bar\a^\Pi(t_{k-1});u(\cd)),\qq\forall u(\cd)\in\cU[t_{k-1},t_k],\ea$$
with $(\bar Y^k(\cd),\bar Z^k(\cd),\bar\G^k(\cd\,,\cd))$ being the adapted solution to the following BSDE:
\bel{BSDE(k)**}\left\{\2n\ba{ll}
\ns\ds dY^k(s)\1n=\1n-g\big(t_{k-1},s,\bar X^\Pi(s),\bar\a^\Pi(s),\bar Y^k(s),\bar Z^k(s),\1n\int_\dbR\2n\bar\G^k(s;\th)\pi(d\th),u(s)\big)ds\\
\ns\ds\qq\qq\qq\qq+\bar Z^k(s)dW(s)+\int_\dbR\bar\G^k(s-,\th)\wt N(ds,d\th),\qq s\in[t_{k-1},t_k),\\
\ns\ds\bar Y^k(t_k)=\Th^k\big(t_k,\bar X^\Pi(t_k),\bar\a^\Pi(t_k)\big).\ea\right.\ee
Note that if we solve the following (comparing with \rf{TH(k)}):
\bel{TH(k)*}\left\{\2n\ba{ll}
\ns\ds\Th^k_s(s,x,i)+\bar\cA^\Pi\Th^k(s,x,i)+g\big(t_{k-1},s,x,i,\Th^k(s,x,i),\\
\ns\ds\qq\qq\qq\Th^k_x(s,x,i)\si(s,x,i,\bar\Psi^\Pi(s,x,i)),[Q(x)\Th^k(s,x,\cd)]_i,\bar\Psi^\Pi(s,x,i)\big)=0;\\
\ns\ds\qq\qq\qq\qq\qq\qq\qq\qq\qq(s,x,i)\in[t_{k-1},T]\times\dbR^n\times M,\\
\ns\ds\Th^k(T,x,i)=h(t_{k-1},x,i),\qq(x,i)\in\dbR^n\times M,\ea\right.\ee
then the value function $V^k(\cd\,,\cd\,,\cd)$ is a restriction of $\Th^k(\cd\,,\cd\,,\cd)$:
\bel{V=TH}V^k(s,x,i)=\Th^k(s,x,i),\qq(s,x,i)\in[t_{k-1},t_k]\times\dbR^n\times M.\ee
Next, we define
\bel{Th(Pi)}\ba{ll}
\ns\ds\Th^\Pi(\t,s,x,i)=\sum_{k=1}^{N-1}\Th^k(s,x,i)I_{[t_{k-1},t_k)}(\t)
+\Th^N(s,x,i)I_{[t_{N-1},t_N]}(\t),\\
\ns\ds\qq\qq\qq\qq\qq\qq\qq\qq\qq0\les\t\les s\les T,~(x,i)\in\dbR^n\times M,\ea\ee
and
$$\ba{ll}
\ns\ds g^\Pi(\t,s,x,i,y,z,\g,u)=\sum_{k=1}^Ng(t^\Pi(\t),s,x,i,y,z,\g,u),\\
\ns\ds\qq\qq\qq\qq0\les\t\les s\les T,(x,i,y,z,\g,u)
\in\dbR^n\times M\times\dbR\times\dbR^d\times\dbR\times U.\ea$$
Then the following holds:
\bel{TH(Pi)}\left\{\2n\ba{ll}
\ns\ds\Th^\Pi_s(\t,s,x,i)+\bar\cA^\Pi\Th^\Pi(\t,s,x,i)+g^\Pi\big(\t,s,x,i,\Th^\Pi(\t,s,x,i),\\
\ns\ds\qq\qq\Th^\Pi_x(\t,s,x,i)\si(s,x,i,\bar\Psi^\Pi(s,x,i)),[Q(x)\Th^\Pi(\t,s,x,\cd)]_i,\bar\Psi^\Pi(s,x,i)\big)=0;\\
\ns\ds\qq\qq\qq\qq\qq\qq\qq\qq\qq(s,x,i)\in[0,T]\times\dbR^n\times M,\\
\ns\ds\Th^\Pi(\t,T,x,i)=h(t^\Pi(\t),x,i),\qq(x,i)\in\dbR^n\times M,\ea\right.\ee
\bel{V=TH}V^\Pi(s,x,i)=\Th^\Pi(t^\Pi(s),s,x,i),\qq(s,x,i)\in[0,T]\times\dbR^n\times M,\ee
and
\bel{barPsi}\ba{ll}
\ns\ds\bar\Psi^\Pi(s,x,i)=\psi\big(t^\Pi(s);s,x,i,\Th^\Pi(t^\Pi(s),x,\cd),\Th^\Pi_x(t^\Pi(s),x,i),
\Th^\Pi_{xx}(t^\Pi(s),x,i)\big),\\
\ns\ds\qq\qq\qq\qq\qq\qq\qq\qq\qq\qq(s,x,i)\in[0,T]\times\dbR^n\times M,\ea\ee
where $t^\Pi(\cd)$ is defined by \rf{t(Pi)}. We call $\bar\Psi^\Pi(\cd\,,\cd\,,\cd)$ an {\it approximate equilibrium strategy} associated with the partition $\Pi$.

\ms

\section{Equilibrium Strategy}

Now, we introduce the following further hypothesis.

\ms

{\bf (H4)} There exists a $\Th(\t,s,x,i)$ such that
$$\ba{ll}
\ns\ds\lim_{\|\Pi\|\to0}\(|\Th^\Pi(\t,s,x,i)\1n-\1n\Th(\t,s,x,i)|\1n+\1n
|\Th^\Pi_x(\t,s,x,i)\1n-\1n\Th_x(\t,s,x,i)|\1n+\1n|\Th_{xx}^\Pi(\t,s,x,i)\1n
-\1n\Th_{xx}(\t,s,x,i)|\)=0,\ea$$
uniformly for $(\t,s,x)$ in any compact sets.

\bs

By (H0)--(H4), taking $\|\Pi\|\to0$ in \rf{TH(Pi)} and \rf{barPsi}, we have the following:
\bel{Psi}\bar\Psi(s,x,i)=\psi(s;s,x,i,\Th(s,s,x,\cd),\Th_x(s,s,x,i),\Th_{xx}(s,s,x,i)),\ee
and $\Th(\t,t,x,i)$ satisfies the following non-linear PDE system, named {\it Equilibrium HJB equation}:
\bel{TH}\left\{\2n\ba{ll}
\ns\ds\Th_s(\t,s,x,i)+\bar\cA\Th(\t,s,x,i)+g\big(\t,s,x,i,\Th(\t,s,x,i),\\
\ns\ds\qq\qq\Th_x(\t,s,x,i)\si(s,x,i,\bar\Psi(s,x,i)),[Q(x)\Th(\t,s,x,\cd)]_i,\bar\Psi(s,x,i)\big)=0;\\
\ns\ds\qq\qq\qq\qq\qq\qq\qq\qq\qq(s,x,i)\in[0,T]\times\dbR^n\times M,\\
\ns\ds\Th(\t,T,x,i)=h(\t,x,i),\qq\qq(x,i)\in\dbR^n\times M,\ea\right.\ee
where
\bel{bar cA}\ba{ll}
\ns\ds\bar\cA\Th(\t,s,x,i)={1\over2}\tr\[\si\big(s,x,i,\bar\Psi(s,x,i)\big)^\top
\Th_{xx}(\t,s,x,i)\si\big(s,x,i,\bar\Psi(s,x,i)\big)\]\\
\ns\ds\qq\qq\qq\qq\qq+\Th_x(\t,s,x,i)b(s,x,i,\bar\Psi(s,x,i))+\big[Q(x)\Th(\t,s,x,\cd))\big]_i\,.\ea\ee
By the definition of $\dbH$ (see \rf{H}) and $\psi$ (see (H3)), we see that the above equilibrium HJB equation can also be written as
\bel{Th-HJB}\left\{\2n\ba{ll}
\ns\ds\Th_s(\t,s,x,i)+\dbH(\t,s,x,i,\Th(\t,s,x,\cd),\Th_x(\t,s,x,i),\Th_{xx}(\t,s,x,i),\bar\Psi(s,x,i))=0,\\
\ns\ds\qq\qq\qq\qq\qq\qq\qq\qq\qq\qq(s,x,i)\in[0,T]\times\dbR^n\times M,\\
\ns\ds\Th(\t,T,x,i)=h(\t,x,i),\qq\qq(x,i)\in\dbR^n\times M.\ea\right.\ee
From \rf{V=TH}, we see that $V^\Pi(\cd\,,\cd\,,\cd)$ converges to some $V(\cd\,,\cd\,,\cd)$ with
\bel{V=TH*}V(s,x,i)=\Th(s,s,x,i),\qq(s,x,i)\in[0,T]\times\dbR^n\times M,\ee
which is called the {\it equilibrium value function}. We call $\bar\Psi(\cd\,,\cd\,,\cd)$ defined by \rf{Psi}
an {\it equilibrium strategy} of our original Problem (N) on $[0,T]$.

\ms

Note that (H4) will be automatically satisfied if we are able to show that the equilibrium HJB equation is well-posed, which will be done in the next section, under proper conditions.

\ms

We now show that the equilibrium strategy $\bar\Psi(\cd\,,\cd\,,\cd)$ is {\it approximately locally optimal} in a suitable sense. Such a result can be viewed as a {\it verification theorem} for our equilibrium strategy. The main idea is taken from \cite{Wei-Yong-Yu2016}. But, due to the appearance of the regime switching process $\a(\cd)$, some delicate technical modifications will be introduced.

\ms

We call the following the {\it equilibrium system}:
\bel{equilibrium system}\left\{\2n\ba{ll}
\ns\ds d\bar X(s)=b\big(s,\bar X(s),\bar\Psi(s,\bar X(s),\bar\a(s))\big)ds
+\si\big(s,\bar X(s),\bar\Psi(s,\bar X(s),\bar\a(s))\big)dW(s),\qq s\in[0,T],\\
\ns\ds d\bar\a(s)=\int_\dbR\m(\bar X(s),\bar\a(s-),\th)N(ds,d\th),\qq s\in[0,T],\\
\ns\ds d\bar Y(t,s)=-g\big(t,s,\bar X(s),\bar Y(t,s),\bar Z(t,s),\int_\dbR\bar\G(t,s,\th)\pi(d\th),\bar\Psi(s,\bar
X(s),\bar\a(s))\big)ds\\
\ns\ds\qq\qq\qq+\bar Z(t,s)dW(s)+\int_\dbR\bar\G(t,s-,\th)\wt N(ds,d\th),\qq s\in[t,T],\\
\ns\ds\bar X(0)=x,\qq\bar\a(0)=i,\qq\bar Y(t,T)=h\big(t,\bar X(T),\bar\a(T)\big),\ea\right.\ee
It is known that if $(\bar X(\cd),\bar\a(\cd),\bar Y(\cd\,;\cd),\bar Z(\cd\,;\cd),\bar\G(\cd\,;\cd\,,\cd))$ is an adapted solution to the above system, then the following representation holds:
\bel{Y=Th*}\left\{\2n\ba{ll}
\ns\ds\bar Y(t;s)=\Th(t;s,\bar X(t;s),\bar\a(t;s)),\\
\ns\ds\bar Z(t;s)\1n=\1n\Th_x(t;s,\bar X(t;s),\bar\a(t;s))\si\big(s,\bar X(t;s),\bar\a(t;s),\bar\Psi(t;s,\bar X(t;s),\bar\a(t;s))\big),\\
\ns\ds\bar\G(t;s,\th)=\Th(t;s,\bar X(t;s),\bar\a(t;s)+\m(\bar X(t;s),\bar\a(t;s),\th))-\Th(t;s,\bar X(t;s),\bar\a(t;s)),\ea\right.\q s\in[\t,T].\ee
and the equilibrium value function is given by:
\bel{bar J}J(t,\bar X(t),\bar\a(t);\bar\Psi(\cd\,,\cd\,,\cd))=\bar Y(t,t)=\Th(t,t,\bar X(t))=V(t,\bar X(t)).\ee

We now assume that our equilibrium HJB equation system \rf{TH} has a  classical solution and all the involved functions in the equations are differentiable with bounded derivatives. We prefer not to get into the most general and technical situations so that the main clue will be clear in our presentation.

\ms

First, for any given $(x,i)\in\dbR^n\times M$, let $(\bar X(\cd),\bar\a(\cd))$ be the state process over $[0,T]$, under the equilibrium strategy $\bar\Psi(\cd\,,\cd\,,\cd)$, i.e.,
\bel{equilibrium FSDE}\left\{\2n\ba{ll}
\ds d\bar X(s)=b\big(s,\bar X(s),\bar\a(s),\bar\Psi(s,\bar X(s),\bar\a(s))\big)ds\\
\ns\ds\qq\qq\qq+\si\big(s,\bar X(s),\bar\a(s),\bar\Psi(s,\bar X(s),\bar\a(s))\big)dW(s),\qq s\in[0,T],\\
\ns\ds d\bar\a(s)=\int_\dbR\m(\bar X(s),\bar\a(s-),\th)N(ds,d\th),\qq s\in[0,T],\\
\ns\ds\bar X(0)=x,\qq\bar\a(0)=i.\ea\right.\ee
Next, let $t\in[0,T)$ be given and $\e>0$ be small so that $t+\e\les T$. Then the adapted solution $(\bar Y(t,\cd),\bar Z(t,\cd))$ of the BSDE in \rf{equilibrium system} over $[t,T]$ also satisfies the following on $[t,t+\e]$:
\bel{BSDE[t,t+e]}\left\{\2n\ba{ll}
\ds d\bar Y(t,s)=-g\big(t,s,\bar X(s),\bar Y(t,s),\bar Z(t,s),\int_\dbR\bar\G(t,s,\th)\pi(d\th),\bar\Psi(s,\bar
X(s),\bar\a(s))\big)ds\\
\ns\ds\qq\qq\qq+\bar Z(t,s)dW(s)+\int_\dbR\bar\G(t,s-,\th)\wt N(ds,d\th),\qq s\in[t,t+\e],\\
\ns\ds \bar Y(t,t+\e)=\Theta(t;t+\e,\bar X(t+\e),\bar\alpha(t+\e)).\ea\right.\ee
Take any $u(\cd)\in\cU[t,t+\e]$, we consider the state equation over $[t,T]$, with the initial triple $(t,\bar X(t),\bar\a(t))$, under the following control:
\bel{u+Psi}\big[u(\cd)\oplus\bar\Psi(\cd\,,\cd\,,\cd)](s)=\left\{\2n\ba{ll}
\ds u(s),\qq\qq s\in[t,t+\e),\\
\ns\ds\bar\Psi(s,\cdot,\cdot),\qq s\in[t+\e,T],\ea\right.\ee
Denote the corresponding state process by $(X^\e(\cd),\a^\e(\cd))$ which satisfies the following:
\bel{e-FSDE}\left\{\2n\ba{ll}
\ds dX^\e(s)=b\big(s,X^\e(s),\a^\e(s),u(s))\big)ds
+\si\big(s,X^\e(s),\a^\e(s),u(s))\big)dW(s),\qq s\in[t,t+\e),\\
\ns\ds dX^\e(s)=b\big(s,X^\e(s),\a^\e(s),\bar\Psi(s,X^\e(s),\a^\e(s))\big)ds\\
\ns\ds\qq\qq\qq+\si\big(s,X^\e(s),\a^\e(s),\bar\Psi(s,X^\e(s),\a^\e(s))\big)dW(s),
\q s\in[t+\e,T],\\
\ns\ds d\a^\e(s)=\int_\dbR\m(X^\e(s),\a^\e(s-),\th)N(ds,d\th),\qq s\in[t,T],\\
\ns\ds X^\e(t)=\bar X(t),\qq\a^\e(t)=\bar\a(t).\ea\right.\ee
The adapted solution of the corresponding BSDE over $[t,t+\e]$ is denoted by $(Y^\e(t,\cd),Z^\e(t,\cd))$ with the equation being as follows:
\bel{e-BSDE}\left\{\2n\ba{ll}
\ds dY^\e(t,s)=-g\big(t,s, X^\e(s),Y^\e(t,s), Z^\e(t,s),\int_\dbR\G^\e(t,s,\th)\pi(d\th),u(s)\big)ds\\
\ns\ds\qq\qq\qq+ Z^\e(t,s)dW(s)+\int_\dbR \G^\e(t,s-,\th)\wt N(ds,d\th),\qq s\in[t,t+\e],\\
\ns\ds Y^\e(t,t+\e)=\Theta(t;t+\e, X^\e(t+\e),\alpha^\e(t+\e)).\ea\right.\ee
Now, for the optimal control problem with the state equation \rf{e-FSDE} and the cost functional
$$J\big(t,\bar X(t),\bar\a(t);u(\cd)\oplus\bar\Psi(\cd\,,\cd\,,\cd)\big|_{[t+\e,T]}\big)=Y^\e(t,t),$$
let $\bar\Psi^\e(\cd\,,\cd\,,\cd)$ be the closed-loop representation of the optimal control. Thus,
\bel{J<J*}\ba{ll}
\ns\ds J\big(t,\bar X(t),\bar\a(t);\bar\Psi^\e(\cd\,,\cd\,,\cd)\oplus\bar\Psi(\cd\,,\cd\,,\cd)\big|_{[t+\e,T]}\big)\les J\big(t,\bar X(t),\bar\a(t);u(\cd)\oplus\bar\Psi(\cd\,,\cd\,,\cd)\big|_{[t+\e,T]}\big),\\
\ns\ds\qq\qq\qq\qq\qq\qq\qq\qq\forall u(\cd)\in\cU[t,t+\e],\ea\ee
where
\bel{ePsi}\big[\bar\Psi^\e(\cd\,,\cd\,,\cd)\oplus\bar\Psi(\cd\,,\cd\,,\cd)\big|_{[t+\e,T]}
\big](s)=\left\{\2n\ba{ll}
\ns\ds\bar\Psi^\e(s,\cd\,,\cd),\qq\qq s\in[t,t+\e),\\
\ns\ds\bar\Psi(s,\cd\,,\cd),\qq\qq s\in[t+\e,T].\ea\right.\ee
It is very important to know that
\bel{barPsi-e}\bar\Psi^\e(s,x,i)=\psi(t,s,x,i,\Th(t,s,x,\cd),\Th_x(t,s,x,i),\Th_{xx}(t,s,x,i)),
\q s\in[t,t+\e].\ee
Therefore, under our conditions (especially \rf{Lip-h-g}), one has
\bel{Psi-Psi}\ba{ll}
\ns\ds|\bar\Psi^\e(s,x,i)-\bar\Psi(s,x,i)|\\
\ns\ds=|\psi(t,s,x,\Th(t,s,x,\cd),\Th_x(t,s,x,i),\Th_{xx}(t,s,x,i))\\
\ns\ds\qq\qq-\psi(s,s,x,\Th(s,s,x,\cd),\Th_x(s,s,x,i),\Th_{xx}(s,s,x,i))|
\les K|t-s|\les K\e,\q s\in[t,t+\e].\ea\ee
This will play a crucial role below. Let the state process corresponding to \rf{ePsi} be denoted by $(\bar X^\e(\cd),\bar\a^\e(\cd))$. Then
\bel{e-FSDE*}\left\{\2n\ba{ll}
\ds d\bar X^\e(s)=b\big(s,\bar X^\e(s),\bar\a^\e(s),\bar\Psi^\e(s,\bar X^\e(s),\bar\a^\e(s))\big)ds\\
\ns\ds\qq\qq\qq+\si\big(s,\bar X^\e(s),\bar\a^\e(s),\bar\Psi^\e(s,\bar X^\e(s),\bar\a^\e(s))\big)dW(s),\qq s\in[t,t+\e),\\
\ns\ds d\bar X^\e(s)=b\big(s,\bar X^\e(s),\bar\a^\e(s),\bar\Psi(s,\bar X^\e(s),\bar\a^\e(s))\big)ds\\
\ns\ds\qq\qq\qq+\si\big(s,\bar X^\e(s),\bar\a^\e(s),\bar\Psi(s,\bar X^\e(s),\bar\a^\e(s))\big)dW(s),
\q s\in[t+\e,T],\\
\ns\ds d\bar\a^\e(s)=\int_\dbR\m(\bar X^\e(s),\bar\a^\e(s-),\th)N(ds,d\th),\qq s\in[t,T],\\
\ns\ds\bar X^\e(t)=\bar X(t),\qq\bar\a^\e(t)=\bar\a(t).\ea\right.\ee
The adapted solution of the corresponding BSDE over $[t,t+\e]$ is denoted by $(\bar Y^\e(t;\cd),\bar Z^\e(t;\cd))$ with the equation being as follows:
\bel{e-BSDE}\left\{\2n\ba{ll}
\ns\ds d\bar Y^\e(t;s)=-g\big(t,s,\bar X^\e(s),\bar Y^\e(t;s),\bar Z^\e(t;s),\int_\dbR\bar\G^\e(t,s,\th)\pi(d\th),\bar\Psi^\e(s,\bar X^\e(s),\bar\a^\e(s))\big)ds\\
\ns\ds\qq\qq\qq+\bar Z^\e(t;s)dW(s)+\int_\dbR\bar\G^\e(t;s-,\th)\wt N(ds,d\th),\qq s\in[t,t+\e],\\
\ns\ds\bar Y^\e(t;t+\e)=\Th(t;t+\e, \bar X^\e(t+\e),\bar\a^\e(t+\e)).\ea\right.\ee
With the above preparation, we are now ready to present the {\it approximately local optimality} of the equilibrium strategy.

\bt{localopti-2} \sl Let {\rm(H0)--(H4)} hold. Then the following holds:
\bel{lim J-J>0-2}\liminf_{\e\to0}{J\big( t,\bar X(t),\bar\a(t);u(\cd)\oplus\bar\Psi(\cd\,,\cd\,,\cd)\big|_{[t+\e,T]}\big)-J( t,\bar X(t),\bar\a(t);\bar\Psi(\cd\,,\cd)\big|_{[t,T]}\big)\over\e}\ges0.\ee

\et

To prove the above theorem, we need two lemmas. The following lemma gives the continuity of the switching diffusion with respect to its initial state.

\bl{} \sl Let {\rm(H0)--(H1)} hold. Let $(X_1(\cd),\a_1(\cd))$ and $(X_2(\cd),\a_2(\cd))$ be the solutions of the state equation on $[t,T]$, under the equilibrium strategy $\bar\Psi(\cd\,,\cd\,,\cd)$ with initial triples $(t,x_1,i_0)$ and $(t,x_2,i_0)$, respectively. Let
$$A_s:=\big\{\o\in\O\bigm|\a_1(r,\o)=\a_2(r,\o),~r\in[t,s]\big\}.$$
Then for some constant $K>0$, the following holds:
\bel{stability-0}\dbP(A_T^c)^2+\dbE\[\sup_{t\les r\les T}|X_1(r)- X_2(r)|^2I_{A_T}\]\les K|x_1-x_2|^2.\ee

\el

The above has exact the same form as \rf{P(Ac)+|X-X|} in Theorem \ref{well-posedness}. The only difference is that it was an open-loop case there and it is a closed-loop case here. The proof is essentially the same. Therefore, we omit it here.

\ms

The following  lemma gives an approximation to the difference between $(\bar X(\cd),\bar\a(\cd))$ and $(\bar X^\e(\cd),\bar\a^\e(\cd))$.

\bl{lemma 4.3.} \sl Let $(\bar X(\cd),\bar\a(\cd))$ and $(\bar X^\e(\cd),\bar\a^\e(\cd))$ be defined as above. Let
$$A^\e_s=\big\{\bar\a(\t)=\bar\a^\e(\t)\bigm|t\les\t\les s\big\}.$$
Then for any $r\ges1$, there exists a constant $K>0$ such that
\bel{Xdiff-ep*}\dbE\[\sup_{t\les s\les t+\e}\big|\bar X(s)-\bar X^\e(s)\big|^rI_{A_s^\e}\]\les K\e^{3r\over2},\qq\forall\e>0,\ee
and for any $p>2$,
\bel{ATe}\dbP\big[(A_{t+\e}^\e)^c\big]+\dbP\big[(A_T^\e)^c\big]\les K\e^{1+{p-2\over2(p+1)}},\qq\forall\e>0.\ee

\el
\it Proof. \rm Now we deal the closed-loop case here. Thus we have
$$\dbE_\t\[\sup_{\t\les s\les T}\(|\bar X(s)|^r+|\bar X(s)|^r\)\les K_r(1+|x|^r),\qq\forall r\ges1.$$
Similar to \rf{probaA-1} with $u(\cd)$ absent (since we have a closed-loop
case here), we have
\bel{}\dbP\big[(A^\e_{t+\e})^c\big]\les K\e\big(\e+\e^{p-2\over2}\d^{-p}+\d\big).\ee
We now for $p>2$, take $\d=\e^{p-2\over{2(p+1)}}$. Then the above leads to
\bel{}\dbP\big[(A^\e_{t+\e})^c\big]\les K\e^{3p\over2(p+1)}.\ee
Next, define $\t=\inf\{s\in[t,T]\bigm|\bar\a(s)\ne\bar\a^\e(s)\}$. Then it is standard that, making use the Lipschitz continuity in $u$ assumed in (H1), together with \rf{Psi-Psi},
$$\ba{ll}
\ns\ds\dbE_t\[\sup_{s\in[t,t\land\t]}|\bar X(s)-\bar X^\e(s)|^r\]\\
\ns\ds\les K\dbE_t\[\(\int_t^{(t+\e)\land\t}|b\big(s,\bar X(s),\bar\a(s),
\bar\Psi(s,\bar X(s),\bar\a(s))\big)-b\big(s,\bar X(s),\bar\a(s),
\bar\Psi^\e(s,\bar X(s),\bar\a(s))\big)|ds\)^r\\
\ns\ds\qq\qq+\(\int_t^{(t+\e)\land\t}|\si\big(s,\bar X(s),\bar\a(s),
\bar\Psi(s,\bar X(s),\bar\a(s))\big)-\si\big(s,\bar X(s),\bar\a(s),
\bar\Psi^\e(s,\bar X(s),\bar\a(s))\big)|^2ds\)^{r\over2}\]\\
\ns\ds\les K\dbE_t\[\int_t^{(t+\e)\land\t}|\bar\Psi(s,\bar X(s),\bar\a(s))-\bar\Psi^\e(s,\bar X(s),\bar\a(s))|ds\)^r\\
\ns\ds\qq\qq+\(\int_t^{(t+\e)\land\t}|\bar\Psi(s,\bar X(s),\bar\a(s))-\bar\Psi^\e(s,\bar X(s),\bar\a(s))|^2ds\)^{r\over2}\]\\
\ns\ds\les K\[\(\int_t^{t+\e}K\e ds\)^r+\(\int_t^{t+\e}K^2\e^2ds\)^{r\over2}\]\les K\(\e^{2r}+\e^{3r\over2}\)\les K\e^{3r\over2}.\ea$$
Hence,
\bel{Xdiff-ep}\dbE\[\sup_{t\les s\les t+\e}\big|\bar X(s)-\bar X^\e(s)\big|^rI_{A_s^\e}\]\les\dbE\[\sup_{t\les s\les t+\e}\big|\bar X(s\land\t)-\bar X^\e(s\land\t)\big|^r\]\les K\e^{3r\over2}.\ee
Compare $(\bar X(s),\bar\a(s))$ with $(\bar X^\e(s),
\bar\a^\e(s))$, for $s\in[t+\e,T]$. On $A_{t+\e}$, one has $\a(t+\e)=\bar\a^\e(t+\e)$. Thus by
\rf{Xdiff-ep},
$$\ba{ll}
\ns\ds\dbP\big[(A_T^\e)^c\big]=\dbP\big[(A_{t+\e}^\e)^c\big]
+\dbP\big[(A_T^\e)^c\cap A^\e_{t+\e}\big]\\
\ns\ds\qq\qq\les K\e^{3p\over2(p+1)}+\dbE\big\{\dbP_{t+\e}
\big[(A_T^\e)^c\cap A_{t+\e}\big]\big\}\\
\ns\ds\qq\qq\les K\e^{3p\over2(p+1)}+K\dbE\big[|\bar X(t+\e)-\bar X^\e(t+\e)|I_{A^\e_{t+\e}}\big]\les K\e^{3p\over2(p+1)}.\ea$$
This gives \rf{ATe} since ${3p\over2(p+1)}=1+{p-2\over2(p+1)}$. \endpf

\ms

\it Proof of Theorem \ref{localopti-2}. \rm Since $(\bar X(\cd),
\bar\a(\cd))$ and $(\bar X^\e(\cd),\bar\a^\e(\cd))$ follow the
same equation on $[t+\e,T]$, it is easy to see that for
$1<r<2$, similar to \rf{stability-3}--\rf{P(Ac)*},
$$\dbE_t\[\sup_{t+\e\les s\les T}|\bar X(s)
-\bar X^\e(s)|^rI_{A^\e_s}\]\les K
\dbE_t|\bar X(t+\e)-\bar X^\e(t+\e)|^r\les K\e^{3r\over2}.$$
On the other hand, since
$$\dbE_t\[\sup_{t+\e\les s\les T}\(|\bar X(s)|^p
+|\bar X^\e(s)|^p\)\]\les K_p,\qq\forall p\ges1,$$
we have, taking proper $\rho>1$,
$$\ba{ll}
\ns\ds\dbE_t\[\sup_{t+\e\les s\les T}|\bar X(s)
-\bar X^\e(s)|^rI_{(A^\e_T)^c}\]\\
\ns\ds\les K\Big\{1+\[\dbE_t\(\sup_{t+\e\les s\les T}
|\bar X(s)|^{r\rho}\)\]^{1\over\rho}+\[\dbE_t\(\sup_{t+\e\les s\les T}
|\bar X^\e(s)|^{r\rho}\)\]^{1\over\rho}\Big\}\(\dbP\big[(A^\e_T)\big]^c
\)^{\rho-1\over\rho}\\
\ns\ds\les K\e^{{3p\over2(p+1)}{\rho-1\over\rho}},\ea $$
with $p>4$ (see \rf{ATe}). By the $L^r$-estimate of BSDE on $[t+\e,T]$ for the case $r\in(1,2)$ (see Proposition 3 in \cite{Kruse-Popier2016,Kruse-Popier2017}) and the representation of $(Y(\cd),Z(\cd))$ in terms of $X(\cd)$,
$$\ba{ll}
\ns\ds\dbE_t\[\sup_{s\in[t+\e,T]}|\bar Y(t,s)-\bar Y^\e(t,s)|^r\]
\les K\dbE_t\[|h(t,\bar X(T))-h(t,\bar X^\e(T))|^r\\
\ns\ds\qq\qq+K\(\int_{t+\e}^T\big|g\big(t,s,\bar X(s),\bar\a(s),
\bar Y(t;s),\bar Z(t;s),\int_\dbR\bar\G(t;s,\th)\pi(d\th),
\bar\Psi(s,\bar X(s),\bar\a(s))\big) \\
\ns\ds\qq\qq\qq\qq-g\big(t,s,\bar X^\e(s),\bar\a^\e(s),
\bar Y(t;s),
\bar Z(t;s),\int_\dbR\bar\G(t;s,\th)\pi(d\th),\bar\Psi(s,\bar X^\e(s),\bar\a^\e(s))\big)\big|ds\)^r\]\\
\ns\ds\les K\[\dbE_t\(\sup_{t+\e\les s\les T}|\bar X(s)-\bar X^\e(s)|^rI_{A^\e_T}\)+\dbE_t\(\sup_{t+\e\les s\les T}|\bar X(r)
-\bar X^\e(s)|^rI_{(A^\e_T)^c}\)\]\\
\ns\ds\les K\(\e^{3r\over2}+\e^{3p(\rho-1)\over2(p+1)\rho}\).\ea$$
Use the same estimate on $[t,t+\e]$, note that $|\bar\Psi^\e(\cd,\cd,\cd)-\bar\Psi(\cd,\cd,\cd)|\les K\e$
(see \rf{Psi-Psi}), one has
$$\ba{ll}
\ns\ds\dbE_t\[\sup_{s\in[t,t+\e]}|\bar Y(t,s)
-\bar Y^\e(t,s)|^r\]
\les K\dbE_t\[|\bar Y(t,t+\e)-\bar Y^\e(t,t+\e)|^r\\
\ns\ds\qq\qq+\(\int_t^{t+\e}\big|g\big(t,s,\bar X(s),\bar\a(s),
\bar Y(t;s),\bar Z(t;s),\int_\dbR\bar\G(t;s,\th)\pi(d\th),
\bar\Psi(s,\bar X(s),\bar\a(s))\big) \\
\ns\ds\qq\qq\qq\qq-g\big(t,s,\bar X^\e(s),\bar\a^\e(s),\bar Y(t;s),
\bar Z(t;s),\int_\dbR\bar\G(t;s,\th)\pi(d\th),\bar\Psi^\e(s,\bar X^\e(s),\bar\a^\e(s))\big)\big|ds\)^r\]\\
\ns\ds\les K\(\e^{3r\over2}+\e^{{3p(\rho-1)\over2(p+1)\rho}}
+\e^{r-1}\int_t^{t+\e}|g(t,s,\bar X(s),\bar\a(s),
\bar\Psi(s,\bar X(s),\bar\a(s)))\\
\ns\ds\qq\qq\qq\qq\qq\qq\qq\qq\qq-g(t,s,\bar X^\e(s),\bar \a^\e(s),\bar\Psi^\e(s,\bar X^\e(s),\bar \a^\e(s)))|^rds\)\\
\ns\ds\les K\(\e^{3r\over2}+\e^{{3p(\rho-1)\over2(p+1)\rho}}
+\e^{r-1}\int_t^{t+\e}|\bar\Psi(s,\bar X^\e(s),\bar\a^\e(s)))-
\bar\Psi^\e(s,\bar X^\e(s),\bar\a^\e(s)))|^rds\\
\ns\ds\qq+\e^{r-1}\int_t^{t+\e}|g(t,s,\bar X(s),\bar\a(s),
\bar\Psi(s,\bar X(s),\bar\a(s)))-g(t,s,\bar X^\e(s),\bar \a^\e(s),\bar\Psi(s,\bar X^\e(s),\bar \a^\e(s)))|^rds\)\ea$$
$$\ba{ll}
\ns\ds\les K\[\e^{3r\over2}+\e^{{3p(\rho-1)\over2(p+1)\rho}}+\e^{2r}+\e^r\dbE_t\(\sup_{t\les s\les t+\e}|\bar X(s)
-\bar X^\e(s)|^r\big[I_{A_{t+\e}^\e}+I_{(A^\e_{t+\e})^c}\big]\)\]\\
\ns\ds\les K\Big\{\e^{3r\over2}+\e^{{3p(\rho-1)\over2(p+1)\rho}}+\e^{2r}+\e^{5r\over2}
+\e^r\[\dbE_t\(\sup_{t\les s\les t+\e}\big(|\bar X(s)|^{r\rho}
+|\bar X^\e(s)|^{r\rho}\)\]^{1\over\rho}
\(\dbP\big[(A^\e_{t+\e})^c\big]\)^{\rho-1\over\rho}\Big\}\\
\ns\ds\les K\Big\{\e^{3r\over2}+\e^{{3p(\rho-1)\over2(p+1)\rho}}+\e^{2r}+\e^{5r\over2}
+\e^r\e^{{3p(\rho-1)\over2(p+1)\rho}}\Big\}\les K\e^{{3p(\rho-1)\over2(p+1)\rho}}.\ea$$
In the above, we have suppressed $(\bar Y(t;s),\bar Z(t;s),\int_\dbR\bar\G(t;s,\th)\pi(d\th))$ in the second inequality. Consequently,
$$|\bar Y(t;t)-\bar Y^\e(t;t)|\les K\e^{{3p(\rho-1)\over2r(p+1)\rho}}.$$
Note that
$$\lim_{\rho\to\infty,r\to1}{3p(\rho-1)\over2r(p+1)\rho}
={3p\over2(p+1)}=1+{p-2\over2(p+1)}.$$
Hence, we may choose $1<r<2$ and $\rho>1$ so that
$$|\bar Y(t;t)-\bar Y^\e(t;t)|\les K\e^{1+{p-2\over4(p+1)}}.$$
Then, together with \eqref{J<J*}, our conclusion is proved. \endpf

\ms

Theorem \ref{localopti-2} is referred to as a {\it local optimality} of $\bar\Psi(\cd,\cd,\cd)$. Because of the above, it is proper for us to modify the definition of equilibrium strategy as follows.

\bde{equilibrium strategy} \rm A map $\Psi:[0,T]\times\dbR^n\times M\to U$ is called an equilibrium strategy of Problem (N) if the following hold:

\ms

(i) For any $(\t,\xi,\iota)\in\cD$, the following closed-loop system
\bel{closed-loop(tx)}\left\{\2n\ba{ll}
\ns\ds dX(s)=b(s,X(t),\a(s),\Psi(s,X(s),\a(s)))ds\\
\ns\ds\qq\qq\qq\qq+\si(s,X(s),\a(s),\Psi(s,X(s),\a(s)))dW(s),\qq s\in[\t,T],\\
\ns\ds d\a(s)=\int_{\dbR}\m(X(s),\a(s-),\th)N(ds,d\th),\qq s\in[\t,T],\\
\ns\ds X(\t)=\xi,\qq\a(\t)=\iota,\ea\right.\ee
admits a unique solution $(X(\cd),\a(\cd))\equiv\big(X(\cd\,;\t,\xi,\iota,\Psi(\cd)),\a(\cd\,;\t,\xi,\iota,\Psi(\cd))\big)$.

\ms

(ii) For any $(x,i)\in\dbR^n\times M$, let $(\bar X(\cd),\bar\a(\cd))$ be the solution to \rf{closed-loop(tx)} with $(\t,\xi,\iota)=(0,x,i)$. For any $0\les t<t+\e\les T$, \eqref{lim J-J>0-2} holds.

\ede

According to the above analysis, we obtain the following result.

\ms

\bt{Existence of equilibrium} \sl Let {\rm(H0)--(H4)} hold. Then $\bar\Psi(\cd\,,\cd\,,\cd)$ obtained in {\rm(H4)} is an equilibrium strategy of Problem {\rm(N)}.

\et

\section{Well-posedness of Equilibrium HJB Equation}

In this section, we will look at the well-posedness of equilibrium HJB equation \rf{Th-HJB}, which will complete our approach presented in previous sections. The major assumption that we will assume is the following:

\ms

{\bf(H5)} The diffusion of the state equation is independent of the control, i.e.,
\bel{si}\si(t,x,i,u)\equiv\si(t,x,i),\qq\forall(t,x,i,u)\in[0,T]\times\dbR^n\times M\times U.\ee
%
There exists a constant $\l\in(0,1]$ such that
\bel{a}\l I\les{1\over2}\si(t,x,i)\si(t,x,i)^\top\equiv a(t,x,i)\les{1\over\l}I,\qq\forall(t,x,i)\in[0,T]\times\dbR^n\times M.\ee

\ms

At the moment, more general case (for examples, $\si$ also depends on $u$, $\si$ is not non-degenerate, etc.) is left open.

\ms

Under (H5), we may rewrite \rf{TH} as follows (omitting the bar here):
\bel{5.3}\left\{\2n\ba{ll}
\ns\ds\Th_s(\t,s,x,i)+\tr[a(s,x,i)\Th_{xx}(\t,s,x,i)]\\
\ns\ds\qq+\Th_x(\t,s,x,i)b\big(s,x,i,\psi(s;s,x,i,\Th(s;s,x,\cd),\Th_x(s;s,x,i))
\big)+\big[Q(x)\Th(\t,s,x,\cd)\big]_i\\
\ns\ds\qq+g\big(\t,s,x,i,\Th(\t;s,x,\cd),\Th_x(\t;s,x,i)\si(s,x,i),
\big[Q(x)\Th(\t,s,x,\cd)\big]_i,\\
\ns\ds\qq\qq\qq\qq\qq\qq\qq\psi(s;s,x,i,\Th(s;s,x,\cd),\Th_x(s;s,x,i)\big)=0,\\
\ns\ds\qq\qq\qq\qq\qq\qq\qq\qq\qq0\les\t\les s\les T,~x\in\dbR^n,~i\in M,\\
\ns\ds\Th(\t;T,x,i)=h(\t;x,i),\qq(\t,x,i)\in[0,T]\times\dbR^n\times M.\ea\right.\ee
Inspired by the above, we redefine $b$ and $g$ as follows, for simplicity.
\bel{bg}\left\{\2n\ba{ll}
\ns\ds b(s,x,i,v,p):=b(s,x,i,\psi(s;s,x,i,v,p)),\\
\ns\ds g(\t,s,x,i,\n,\pi,v,p):=g\big(\t,s,x,i,\n,\pi\si(s,x,i),[Q(x)\n]_i,\psi(s;s,x,i,v,p)
\big).\ea\right.\ee
Then the above \rf{5.3} can be written as
\bel{5.4}\left\{\2n\ba{ll}
\ns\ds\Th_s(\t,s,x,i)+\tr[a(s,x,i)\Th_{xx}(\t,s,x,i)]+\Th_x(\t,s,x,i)
b\big(s,x,i,\Th(s;s,x,\cd),\Th_x(s;s,x,i)\big)\\
\ns\ds\q+\big[Q(x)\Th(\t,s,x,\cd)\big]_i+g\big(\t,s,x,i,\Th(\t;s,x,\cd),
\Th_x(\t;s,x,i),\Th(s;s,x,\cd),\Th_x(s;s,x,i)\big)=0,\\
\ns\ds\qq\qq\qq\qq\qq\qq\qq\qq\qq0\les\t\les s\les T,~x\in\dbR^n,~i\in M,\\
\ns\ds\Th(\t;T,x,i)=h(\t;x,i),\qq(\t,x,i)\in[0,T]\times\dbR^n\times M.\ea\right.\ee

To establish the well-posedeness of \eqref{5.4}, let us make some preparations. For $\b\in(0,1)$, let $C^\b(\dbR^n\times M)$ be the space of function $\f:\dbR^n\times M\to\dbR$ such that $x\mapsto\f(x,i)$ is
continuous, and
$$\|\f\|_\b:=\|\f\|_0+[\f]_\b<\infty,$$
where
$$\|\f\|_0=\sup_{(x,i)\in\dbR^n\times M}|\f(x,i)|,\q[\f]_\b=\sup_{x\ne y,
i\in M}{|\f(x,i)-\f(y,i)|\over|x-y|^\b}.$$
Further let $C^{1+\b}(\dbR^n\times M)$ and $C^{2+\b}(\dbR^n\times M)$ be the space of functions $\f:\dbR^n\times M\to\dbR$ such that
$$\|\f\|_{1+\b}=\|\f\|_0+\|\f_x\|_0+[\f_x]_\b<\infty,$$
and
$$\|\f\|_{2+\b}=\|\f\|_0+\|\f_x\|_0+\|\f_{xx}\|_0+[\f_{xx}]_\b<\infty,$$
respectively. Also let $L^\infty(0,T;C^\beta(\dbR^n\times M))$ be the set of all measurable functions $f:[0,T]\times\dbR^n\times M\to\dbR$ such that for fixed $t\in[0,T]$, $f(t,\cd\,,\cd)\in C^\b(\dbR^n\times M)$ with
$$\|f(\cd\,,\cd\,,\cd)\|_{L^\infty(0,T;C^\b(\dbR^n\times M ))}=\esssup_{t\in[0,T]}\|f(t,\cd\,,\cd)\|_\b<\infty.$$
Let $C([0,T];C^(\dbR^n\times M))$ be the continuous functions in $L^\infty(0,T;C^\b(\dbR^n\times M))$. Similarly, we can define
$C([0,T];C^{k+\b}(\dbR^n\times M))\subset L^\infty(0,T;C^{k+\b}(\dbR^n\times M))$.

\ms

We now introduce the following additional hypothesis.

\ms

{\bf (H6)} The maps $a:[0,T]\times\dbR^n\times M\to\dbS^n$ defined by \rf{a}, $b:[0,T]\times\dbR^n\times M\times\dbR^m\times\dbR^n\to\dbR^n$ and $g:[0,T]\times[0,T]\times\dbR^n\times M\times\dbR^m\times\dbR^n\to\dbR$, $h:[0,T]\times\dbR^n\times M\to\dbR$ are continuous and bounded. Moreover there exists a constant $L>0$ such that
$$\left\{\2n\ba{ll}
\ns\ds |h_x(\cdot,\cd,\cd)|\leq L,\\
\ns\ds|a(t,x_1,i)-a(t,x_2,i)|+|b(t,x_1,i,v_1,p_1)-b(t,x_2,i,v_2,p_2)|\les L(|x_1-x_2|+|v_1-v_2|+|p_1-p_2|),\\
\ns\ds|g(\t_1,t,x_1,i,\n_1,\pi_1,v_1,p_1)-g(\t_2,t,x_2,i,\n_2,\pi_2,v_2,p_2)|
+|h(\t_1,x_1,i)-h(\t_2,x_2,i)|\\
\ns\ds\qq\qq\qq\qq\les L\big(|\t_1-\t_2|+|x_1-x_2|+|\n_1-\n_2|+|\pi_1-\pi_2|+|v_1-v_2|+|p_1-p_2|\big),\\
\ns\ds|q_{ij}(t,x_1)-q_{ij}(t,x_2)|\les L|x_1-x_2|.\ea\right.$$

\bt{} \sl Let {\rm(H6)} holds. Then there exists a unique solution to \eqref{5.4}.

\et

\it Proof. \rm Let $\t\in[0,T)$ and $v(\cd\,,\cd\,,\cd)\in C([0,T];C^1(\dbR^n\times M))$. Consider following nonlinear PDE: for~$s\in[t,T]$ and each $i\in M $,
\bel{6-1}\left\{\2n\ba{ll}
\ns\ds\Th_s(\t;s,x,i)+\tr[a(s,x,i)\Th_{xx}(\t;t,x,i)]\\
\ns\ds~~+\Th_x(\t;s,x,i)b(s,x,i,v(s,x,\cd),v_x(s,x,i))+\sum_{j=1}^mq_{ij}(x)
\Th(\t;s,x,j)\\
\ns\ds~~+g(\t,s,x,i,\Th(\t;s,x,\cd),\Th_x(\t;s,x,i),v(s,x,\cd),v_x(s,x,i))=0,\\
\ns\ds\Th(\t;T,x,i)=h(\t;x,i);\ea\right.\ee
We will show that the mapping $\dbB:v(s,x,i)\mapsto\Th(s;s,x,i)$ in \eqref{6-1} is well-defined and is a contraction in some suitable space. Then the well-posedness of \rf{5.4} will follow.

\ms

\it Step 1. \rm The mapping $\dbB$ is well-defined.

\ms

For fixed $\Th^0(\t;\cd\,,\cd\,,\cd)\in C([0,T];C^1(\dbR^n\times M))$,
consider the following linear PDE:
\bel{6-12}\left\{\2n\ba{ll}
\ns\ds\Th_s(\t;s,x,i)+\tr[a(s,x,i)\Th_{xx}(\t;t,x,i)]\\
\ns\ds~~+\Th_x(\t;s,x,i)b(s,x,i,v(s,x,\cd),v_x(s,x,i))+\sum_{j=1}^mq_{ij}(x)
\Th(\t;s,x,j)\\
\ns\ds~~+g(\t,s,x,i,\Th^0(\t;s,x,\cd),\Th^0_x(\t;s,x,i),v(s,x,\cd),
v_x(s,x,i))=0,\\
\ns\ds\Th(\t;T,x,i)=h(\t;x,i);\ea\right.\ee
By \cite{Yin-Zhu2009} (p.53), we know that the above admits a unique solution $\Th(\cd\,;\cd\,,\cd,\cd)$. By (H6), the following defines a contraction mapping (see \cite{Friedman1964}):
$$\Th^0(\t;s,x,i)\mapsto\Th(\t;s,x,i):C([T-\d,T];C^1(\dbR^n\times M))\to C([T-\d,T];C^1(\dbR^n\times M)),$$
for some $\d>0$ which is an absolute constant. Thus this map admits a unique fixed point $\Th(\t;s,x,i)\in C([T-\d,T];C^1(\dbR^n\times M))$. Repeating such process on $[T-2\d,T-\d]$, $[T-3\d,T-2\d]$, and so on, we can prove that there is a unique solution of \eqref{6-1} in $C([0,T];C^1(\dbR^n\times M))$. Consequently, the mapping $\dbB$ is well-defined.

\ms

\it Step 2. \rm The mapping $\dbB:v(s,x,i)\mapsto\Th(s;s,x,i)$ defined by \eqref{6-1} is a contraction map from $C([T-\d,T];C^1(\dbR^n\times M))$ to itself for some small constant $\d$.

\ms

Under (H6), for each $\t\in[0,T)$ and $i\in M$, the fundamental solution $\F(t,x;s,y,i)$ for the following PDE
$$\Th_s(\t;s,x,i)+\tr[a(s,x,i)\Th_{xx}(\t;s,x,i)]=0$$
can be written as
$$\F(t,x;s,y,i)={1\over(4\pi(s-t))^{n\over2}\{\det[a(s,y,i)]\}^{1\over2}}
\exp\left\{{(x-y)^\top a(s,y,i)^{-1}(x-y)\over4(s-t)}\right\}.$$
Then $\Th(\t;t,x,i)$, the solution of \eqref{6-1}, satisfies the following:
\bel{789}\ba{ll}
\ns\ds\Th(\t;t,x,i)=\int_{\dbR^n}\F(t,x;T,y,i)h(\t;y,i)dy\\
\ns\ds~+\int_t^T\int_{\dbR^n}\F(t,x;s,y,i)\Th_x(\t;s,y,i)b(s,x,i,v(s,y,
\cd),v_x(s,y,i))dyds\\
\ns\ds~+\sum_{j=1}^m\int_t^T\int_{\dbR^n}\F(t,x;s,y,i)q_{ij}(y)\Th(\t;s,y,j)
dyds\\
\ns\ds~+\int_t^T\int_{\dbR^n}\F(t,x;s,y,i)g(\t,s,y,i,\Th(\t;s,x,\cd),
\Th_x(\t;s,x,i),v(s,y,\cd),v_x(s,y,i))dyds.\ea\ee
Direct calculation yields
$$\left\{\2n\ba{ll}
\ns\ds|\F(t,x;s,y,i)|\les K(s-t)^{-{n\over2}}\exp\left\{-\l{|x-y|^2\over4(s-t)}\right\},\\
\ns\ds|\F_x(t,x;s,y,i)|\les K(s-t)^{-{n+1\over2}}\exp\left\{-\l{|x-y|^2\over4(s-t)}\right\},\ea\right.$$
and
$$\F_y(t,x;s,y,i)=-\F_x(t,x;s,y,i)+\F(t,x;s,y,i)\rho(t,x;s,y,i),$$
where
$$\left\{\2n\ba{ll}
\ns\ds\rho(t,x;s,y,i)={(\det[a(s,y,i])_y\over2\det[a(s,y,i])}
+{\lan[a(s,y,i)^{-1}]_y(x-y),x-y\ran\over4(s-t)},\\
\ns\ds\lan[a(s,y,i)^{-1}]_y(x-y),x-y\ran=\begin{pmatrix}
      \lan[a(s,y,i)^{-1}]_{y_1}(x-y),x-y\ran\\
       \vdots \\
       \lan[a(s,y,i)^{-1}]_{y_n}(x-y),x-y\ran\end{pmatrix}.\ea\right.$$
It is easy to check under (P),
$$|\rho(t,x;s,y,i)|\les K(1+{|x-y|^2\over s-t}).$$

Now we want to prove that for any $(t,x,i)\in[0,T]\times\dbR^n\times M$,
\bel{b-Th}|\Th(\t;t,x,i)|+|\Th_x(\t;t,x,i)|\les K\big(1+\|h(\t,\cd\,,\cd)\|_1\big).\ee
By \eqref{789}, one has
\bel{7899}\ba{ll}
\ns\ds|\Th(\t;t,x,i)|\les\Big|\int_{\dbR^n}\F(t,x;T,y,i)h(\t;y,i)dy\Big|\\
\ns\ds\qq\qq+\Big|\int_t^T\int_{\dbR^n}\F(t,x;s,y,i)\Th_x (\t;s,y,i)b(s,y,i,v(s,y,\cd),v_x(s,y,i))dyds\Big|\\
\ns\ds\qq\qq+\Big|\sum_{j=1}^{m }\int_t^T\int_{\dbR^n}\F (t,x;s,y,i)q_{ij}(y)\Th(\t;s,y,j)dyds\Big|\\
\ns\ds\qq\qq+\Big|\int_t^T\3n\int_{\dbR^n}\F (t,x;s,y,i)g(\t,s,y,i,\Th(\t;s,x,\cd),\Th_x(\t;s,x,i),v(s,y,\cd),
v_x(s,y,i))dyds\Big|\\
\ns\ds\les K(1+\|h(\t,\cd\,,\cd)\|_0)+K\int_t^T\int_{\dbR^n}\F(t,x;s,y,i)\(1+|\Th_x (\t;s,y,i)|+\sum_{j=1}^m|\Th(\t;s,y,j)|\)dyds.\ea\ee
Integrating by parts we have
$$\ba{ll}
\ns\ds\Big|\int_{\dbR^n}\F_x(t,x;T,y,i)h(\t,y,i)dy\Big|\\
\ns\ds=\Big|-\int_{\dbR^n}\F_y(t,x;T,y,i)h(\t,y,i)dy-\int_{\dbR^n}\F (t,x;s,y,i)\rho(t,x;T,y,i)h(\t,y,i)dy\Big|\\
\ns\ds\les\Big|\int_{\dbR^n}\F(t,x;T,y,i)h_y(\t,y,i)dy\Big|
+\Big|\int_{\dbR^n}\F(t,x;T,y,i)\rho(t,x;T,y,i)h(\t,y,i)dy\Big|\\
\ns\ds\les K\big(1+\|h(\t,\cd\,,\cd)\|_1\big)\int_{\dbR^n}\({|y-x|^2\over T-t}+1\)\F(t,x;T,y,i)dy\les K\big(1+\|h(\t,\cd\,,\cd)\|_1\big).\ea$$
Therefore, we can conclude that
\bel{7890}\ba{ll}
\ns\ds|\Th_x(\t;t,x,i)|\les\Big|\int_{\dbR^n}\F_x(t,x;T,y,i)h(\t;y,i)dy\Big|\\
\ns\ds\q+\Big|\int_t^T\int_{\dbR^n}\F_x(t,x;s,y,i)\Th_x(\t;s,y,i)b(s,y,i,
v(s,y,\cd),v_x(s,y,i))dyds\Big|\\
\ns\ds\q+\Big|\sum_{j=1}^m\int_t^T\int_{\dbR^n}\F_x(t,x;s,y,i)q_{ij}(y)
\Th(\t;s,y,j)dyds\Big|\\
\ns\ds\q+\Big|\int_t^T\int_{\dbR^n}\F_x(t,x;s,y,i)g(\t,s,y,i,\Th(\t;s,x,\cd),
\Th_x(\t;s,x,i),v(s,y,\cd),v_x(s,y,i))dyds\Big|\\
\ns\ds\les K\big(1+\|h(\t,\cd\,,\cd)\|_1\big)\\
\ns\ds\q+K\int_t^T\int_{\dbR^n}\F_x(t,x;s,y,i)\(1+|\Th_x(\t;s,y,i)|
+\sum_{j=1}^m|\Th(\t;s,y,j)|\)dyds\\
\ns\ds\les K\big(1+\|h(\t,\cd\,,\cd)\|_1\big)\\
\ns\ds\q+K\int_t^T\int_{\dbR^n}\F_x(t,x;s,y,i)\(|\Th_x(\t;s,y,i)|+\sum_{j=1}^m |\Th(\t;s,y,j)|\)dyds.\ea\ee
Combining \eqref{7899} and \eqref{7890}, we have
$$\ba{ll}
\ns\ds\sum_{i=1}^m\Big[|\Th(\t;t,x,i)|+|\Th_x(\t;t,x,i)|\Big]\\
\ns\ds\les K\big(1+\|h(\t,\cd\,,\cd)\|_1\big)+K\int_t^T\int_{\dbR^n}
{\exp\{-\l{|x-y|^2\over4(s-t)}\}\over(s-t)^{n+1\over2}}\sum_{i=1}^m \[|\Th(\t;s,y,i)|+|\Th_x(\t;s,y,i)|\]dyds\\
\ns\ds\les K\big(1+\|h(\t,\cd\,,\cd)\|_1\big)+K\int_t^T(s-t)^{-{1\over2}}
\sup_{x\in\dbR^n}\sum_{i=1}^m\[|\Th(\t;s,x,i)|+|\Th_x(\t;s,x,i)|\]ds.\ea$$
By Grownwall's inequality, \eqref{b-Th} follows.
\ms

Now, let $v^k(\cd\,,\cd\,,i)\in C[0,T]\times\dbR^n$ for $k=0,1$ and each $i\in M$. Let $\Th^k(\t;\cd\,,\cd\,,\cd)$ be the corresponding solutions of \eqref{6-1}. Then
\bel{7-1}\ba{ll}
\ns\ds|\Th^0(\t;t,x,i)-\Th^1(\t;t,x,i)|\\
\ns\ds\q\les\int_t^T\int_{\dbR^n}\F(t,x;s,y,i)|\Th^0_x(\t;s,y,i)
b(s,y,i,v^0(s,y,\cd),v^0_x(s,y,i))\\
\ns\ds\qq\qq\qq\qq\qq-\Th^1_x(\t;s,y,i)b(s,y,i,v^1(s,y,\cd),v^1_x(s,y,i))|dyds\\
\ns\ds\q+\sum_{j=1}^m\int_t^T\int_{\dbR^n}\F(t,x;s,y,i)|q_{ij}(y)(\Th^0
(\t;s,y,j)-\Th^1(\t;s,y,j))|dyds\\
\ns\ds\q+\int_t^T\int_{\dbR^n}\F(t,x;s,y,i)|g(\t,s,y,i,\Th^0(\t;s,x,\cd),
\Th^0_x(\t;s,x,i),v^0(s,y,\cd),v^0_x(s,y,i))\\
\ns\ds~\qq\qq\qq\qq\qq-g(\t,s,y,i,\Th^0(\t;s,x,\cd),\Th^0_x(\t;s,x,i),v^0(s,y,\cd),
v^0_x(s,y,i))|dyds\\
\ns\ds\q\les\int_t^T\int_{\dbR^n}\F(t,x;s,y,i)|\Th^0_x(\t;s,y,i)\\
\ns\ds\qq\qq\qq\qq\qq \cd[b(s,y,i,v^0(s,y,\cd),v^0_x(s,y,i))-b(s,y,i,v^1(s,y,\cd),v^1_x(s,y,i))]|dyds\\
\ns\ds\q+\int_t^T\int_{\dbR^n}\F(t,x;s,y,i)|[\Th^0_x(\t;s,y,i)-\Th^1_x (\t;s,y,i)]b(s,y,i,v^1(s,y,\cd),v^1_x(s,y,i))|dyds\\
\ns\ds\q+\sum_{j=1}^m\int_t^T\int_{\dbR^n}\F(t,x;s,y,i)|q_{ij}(y)
[\Th^0(\t;s,y,j)-\Th^1(\t;s,y,j)]|dyds\\
\ns\ds\q+\int_t^T\int_{\dbR^n}\F(t,x;s,y,i)|g(\t,s,y,i,\Th^0(\t;s,x,\cd),
\Th^0_x(\t;s,x,i),v^0(s,y,\cd),v^0_x(s,y,i)) \\
\ns\ds\qq\qq\qq\qq\qq-g(\t,s,y,i,\Th^1(\t;s,x,\cd),\Th^1_x(\t;s,x,i),v^1(s,y,\cd),
v^1_x(s,y,i))|dyds\\
\ns\ds\q\les K\int_t^T\int_{\dbR^n}\F(t,x;s,y,i)\sum_{j=1}^m \[|\Th^0(\t;s,y,j)-\Th^1(\t;s,y,j)|+|\Th_x^0(\t;s,y,j)-\Th_x^1(\t;s,y,j)|\]dyds\\
\ns\ds\q+K\int_t^T\int_{\dbR^n}\sum_{j=1}^m(1+|\Th_x^0(\t;s,y,j)|
+|\Th^1(\t;s,y,j)|)|v^0(\t;s,y,j))-v^1(\t;s,y,j))|dyds\\
\ns\ds\q+K\int_t^T\int_{\dbR^n}\sum_{j=1}^m(1+|\Th_x^0(\t;s,y,j)|
+|\Th^1(\t;s,y,j)|)|v_x^0(\t;s,y,j))-v_x^1(\t;s,y,j))|dyds\\
\ns\ds\q\les K\int_t^T\int_{\dbR^n}\F(t,x;s,y,i)\Big\{\sum_{j=1}^m\[|\Th^0(\t;s,y,j)
-\Th^1(\t;s,y,j)|+|\Th_x^0(\t;s,y,j)-\Th_x^1(\t;s,y,j)|\]dyds\\
\ns\ds\q+K(1+\|h(\t;\cd)\|_1)\int_t^T\int_{\dbR^n}\sum_{j=1}^m\[|v^0(\t;s,y,j))
-v^1(\t;s,y,j))|+|v_x^0(\t;s,y,j))-v_x^1(\t;s,y,j))|\]dyds\\
\ns\ds\q\les K\int_t^T\int_{\dbR^n}\F(t,x;s,y,i)\sum_{j=1}^m\[|\Th^0(\t;s,y,j)
-\Th^1(\t;s,y,j)|+|\Th_x^0(\t;s,y,j)-\Th_x^1(\t;s,y,j)|\]dyds\\
\ns\ds\q+K(T-t)(1+\|h(\t;\cd\,,\cd)\|_1)\|v^0(\t;\cd\,,\cd\,,\cd))
-v^1(\t;\cd\,,\cd\,,\cd))\|_{L^\infty(0,T,C(\dbR^n\times M))}.\ea\ee
Similarly we can prove
\bel{7-2}\ba{ll}
\ns\ds|\Th_x^0(\t;t,x,i)-\Th_x^1(\t;t,x,i)|\\
\ns\ds\q\les K\int_t^T\int_{\dbR^n}\F_x(t,x;s,y,i)\sum_{j=1}^m\[|\Th^0(\t;s,y,j)
-\Th^1(\t;s,y,j)|+|\Th_x^0(\t;s,y,j)-\Th_x^1(\t;s,y,j)|\]dyds\\
\ns\ds\q+K(T-t)^{1\over2}\|v^0(\cd\,,\cd\,,\cd)-v^1(\cd\,,\cd\,,\cd)\|_{C([t,T]
\times\dbR^n\times M)}(1+\|h\|_{L^\infty([0,T]\times\dbR^n\times M)})\\
\ns\ds\q\les K\int_t^T(s-t)^{-{1\over2}}\sup_{x\in\dbR^n}\sum_{j=1}^m \[|\Th^0(\t;s,x,j)-\Th^1(\t;s,x,j)|+|\Th_x^0(\t;s,x,j)-\Th_x^1(\t;s,x,j)|\]ds\\
\ns\ds\q+K(T-t)^{1\over2}\|v^0(\cd\,,\cd\,,\cd)-v^1(\cd\,,\cd\,,\cd)\|_{C([t,T]
\times\dbR^n\times M)}(1+\|h\|_{L^\infty(0,T,C(\dbR^n\times M))}).\ea\ee
By Grownwall's inequality, \eqref{7-1} and \eqref{7-2} together yields that
\bel{7-3}\ba{ll}
\ns\ds \|\Th^0(\t;\cd\,,\cd\,,\cd)-\Th^1(\t;\cd\,,\cd\,,\cd)\|_{C([t,T]\times
\dbR^n\times M)}\\
\ns\ds\q\les K(T-t)^{1\over2}(1+\|h\|_{L^\infty(0,T,C(\dbR^n\times M))})\|v^0(\cd\,,\cd\,,\cd)-v^1(\cd\,,\cd\,,\cd)||_{C([t,T]\times\dbR^n\times M)}.\ea\ee
From the procedure above, we know that $K$ is independent of $\t$ and $t$. Hence in particular, for $V^k(t,x,i)=\Th^k(t,t,x,i)$,
\bel{7-3}\ba{ll}
\ns\ds\|V^0(\cd\,,\cd\,,\cd)-V^1(\cd\,,\cd\,,\cd)\|_{C([T-\d,T]\times
\dbR^n\times M)}\\
\ns\ds\q\les K\d^{1\over2}(1+\|h\|_{L^\infty([0,T]\times\dbR^n\times
M)})\|v^0(\cd\,,\cd\,,\cd)-v^1(\cd\,,\cd\,,\cd)\|_{C([T-\d,T]
\times\dbR^n\times M)}.\ea\ee
When $\d>0$ is small enough, we get a contraction mapping $v(\cd\,,\cd\,,\cd)\to V(\cd\,,\cd\,,\cd)$ on $C([T-\d,T];C^1(\dbR^n
\times M))$. Then this map admits a unique fixed point. Since we can
get similar estimates on $[T-2\d,T-\d]$, and so on, one can see that
there exists a unique fixed point on the whole space $C([0, T];C^1
(\dbR^n\times M))$. Then we get a solution $V(\cd\,,\cd\,,\cd)$ for \eqref{5.3} when $\t=t$. For different $\t$, \eqref{6-1} now is a non-degenerate parabolic equation with $v(t,x,i)=V(t,x,i)$ fixed.
The solution can be expressed as
$$\ba{ll}
\ns\ds\Th(\t;t,x,i)=\int_{\dbR^n}\F(t,x;T,y,i)h(\t,y,i)dy\\
\ns\ds\qq+\int_t^T\int_{\dbR^n}\F(t,x;s,y,i)\Th_x(\t;s,y,i)
b(s,y,V(s,y,\cd),V_x(s,y,i))dyds\\
\ns\ds\qq+\sum_{j=1}^m\int_t^T\int_{\dbR^n}\F(t,x;s,y,i)q_{ij}(y)
\Th(\t;s,y,j)dyds\\
\ns\ds\qq+\int_t^T\int_{\dbR^n}\F(t,x;s,y,i)g(\t,s,y,i,\Th(\t;s,y,
\cd),\Th_x(\t;s,y,i),V(s,y,\cd),V_x(s,y,i))dyds.\ea$$
Thus we have proved the well-posedeness of \eqref{5.3}.\endpf

\ms

Recalling that in previous section, in (H4), we assume that the
convergence of $\Th^\Pi(\cd\,;\cd\,,\cd\,,\cd)$ to $\Th(\cd\,;\cd\,,\cd\,,\cd)$. Now we will show that this convergence
is reasonable.

\ms

By the well-posedeness of \eqref{5.3}, we know that
$$|V^\Pi(s,x,i)-V(s,x,i)|+|V_x^\Pi(s,x,i)-V_x(s,x,i)|\les K\|\Pi\|.$$
Compare \eqref{TH(Pi)} and \eqref{5.3}, we have
\bel{7-4}\ba{ll}
\ns\ds\Th^\Pi(\t;t,x,i)-\Th(\t;t,x,i)=\int_{\dbR^n}\F(t,x;T,y,i)\big[h^\Pi(\t,y,i)
-h(\t,y,i)]dy\\
\ns\ds\q=\int_t^T\int_{\dbR^n}\F(t,x;s,y,i)\[\Th^\Pi_x(\t;s,y, i)b(s,y,V^\Pi(s,y,\cd),V^\Pi_x(s,y,i))\\
\ns\ds\qq\qq\qq\qq\qq-\Th_x(\t;s,y,i)b(s,y,V(s,y,\cd),V_x(s,y,i))
\]dyds\\
\ns\ds\q+\sum_{j=1}^m\int_t^T\int_{\dbR^n}\F (t,x;s,y,i)q_{ij}(y)(\Th^\Pi(\t;s,y,j)-\Th(\t;s,y,j))dyds\\
\ns\ds\q+\int_t^T\int_{\dbR^n}\F(t,x;s,y,i)\[g^\Pi(\t,s,y,i,\Th^\Pi
(\t;s,y,\cd),\Th^\Pi_x(\t;s,y,i),V^\Pi(s,y,\cd),V^\Pi_x(s,y,i))\\
\ns\ds\qq\qq\qq\qq\qq -g(\t,s,y,i,\Th(\t;s,y,\cd),\Th_x(\t;s,y,i),V(s,y,\cd),V_x(s,y,i))
\]dyds.\ea\ee
Thus we have
$$\ba{ll}
\ns\ds|\Th^\Pi(\t;t,x,i)-\Th(\t;t,x,i)|\\
\ns\ds\les K\|\Pi\|+K\int_t^T\int_{\dbR^n}\F(t,x;s,y,i)\[\sum_{j=1}^m |\Th^\Pi(\t;s,y,j)-\Th(\t;s,y,j)|+|\Th_x^\Pi(\t;s,y,i)
-\Th_x(\t;s,y,i)|\\
\ns\ds\qq\qq\qq\qq\qq\qq+\sum_{j=1}^m|V^\Pi(s,y,j)-V(s,y,j)|
+|V_x^\Pi(s,y,i)-V_x(s,y,i)|\]dyds\\
\ns\ds\les K\|\Pi\|+K\int_t^T\3n\int_{\dbR^n}\F(t,x;s,y,i)\[\sum_{j=1}^m |\Th^\Pi(\t;s,y,j)-\Th(\t;s,y,j)|+|\Th_x^\Pi(\t;s,y,i)-\Th_x(\t;s,y,i)|\]dyds.
\ea$$
Similarly we can prove
$$\ba{ll}
\ns\ds|\Th_x^\Pi(\t;t,x,i)-\Th_x(\t;t,x,i)|\\
\ns\ds\les K\|\Pi\|+K\int_t^T\3n\int_{\dbR^n}\F_x(t,x;s,y,i)\sum_{j=1}^m \[|\Th(\t;s,y,j)-\Th(\t;s,y,j)|+|\Th_x^\Pi(\t;s,y,j)-\Th_x(\t;s,y,j)|\]dyds.\ea$$
Thus we have
$$\ba{ll}
\ns\ds|\Th^\Pi(\t;t,x,i)-\Th(\t;t,x,i)|+|\Th_x^\Pi(\t;t,x,i)-\Th_x(\t;t,x,i)|\\
\ns\ds\les K\|\Pi\|+K\int_t^T\3n\int_{\dbR^n}{\exp\{-\l{|x-y|^2\over4(s-t)}\}\over
(s-t)^{n+1\over2}}\sum_{j=1}^m\[|\Th^\Pi(\t;s,y,j)-\Th(\t;s,y,j)|+|\Th_x^\Pi
(\t;s,y,j)-\Th_x(\t;s,y,j)|\]dyds\\
\ns\ds\les K\|\Pi\|+K\int_t^T(s-t)^{-{1\over2}}\sup_{x\in\dbR^n}\sum_{j=1}^m
\[|\Th^\Pi(\t;s,x,j)-\Th(\t;s,x,j)|+|\Th^\Pi_x(\t;s,x,j)-\Th_x(\t;s,x,j)|\]ds.
\ea$$
By Grownwall's inequality, the expected convergence follows.

\section{An Illustrative Example}

In this section, we present an example to illustrate the main results of the paper.

\ms

Let us first recall the classical Merton's problem. Consider the following SDE:
\bel{M1}\left\{\2n\ba{ll}
\ns\ds dX(s)=\big[bu(s)-c(s)\big]ds+\si u(s)dW(s),\q s\in[t,T],\\
\ns\ds X(t)=x,\ea\right.\ee
with payoff functional
\bel{P1}J(t,x;u(\cd),c(\cd))=\dbE_t\[\int_t^Tg(s)c(s)^\g ds+hX(T)^\g\],\ee
for some $\g\in(0,1)$.

In the above, we have assumed that the interest rate for the bank account is zero, the appreciation rate and the volatility of the risky asset are $b$ and $\si$, respectively; $X(\cd)$ is the wealth process,
$u(\cd)$ is the dollar amount invested in the risky asset, and
$c(\cd)$ is the consumption rate. The problem is to maximize the payoff functional $J(t,x;u(\cd),c(\cd))$, in which $g(\cd)$ and $h$ are certain (positive) weights. Now, we let $M=\{1,2\}$ with $1$
representing the ``bull market'' and $2$ representing the
``bear market'', respectively. We assume that the appreciation rate $b$ and the volatility $\si$ depend on the market index $i\in M$. Thus, the modified wealth process equation takes the following form:
\bel{M2}\left\{\2n\ba{ll}
\ns\ds dX(s)=\big[b(\a(s))u(s)-c(s)\big]ds+\si(\a(s))u(s)dW(s),\q s\in[t,T],\\
\ns\ds X(t)=x,\ea\right.\ee
with $\a(\cd)$ being an $M$-valued process.

\ms

According to \rf{transition2}, we have
$$\ba{ll}
\ns\ds\dbP\big(\a(s+\D s)=2\bigm|X(s)=x,\a(s)=1\big)
=q_{12}(x)\D s+o(\D s),\\
\ns\ds\dbP\big(\a(s+\D s)=1\bigm|X(s)=x,\a(s)=2\big)
=q_{21}(x)\D s+o(\D s).\ea$$
Hence, by assuming $q_{12}(\cd)$ to be increasing, we have that when
$X(s)$ is getting larger, the rate of $\a(s)$ jumping from 1
to 2 is getting larger, whereas, by assuming $q_{21}(\cd)$ to be
decreasing, we have that when $X(s)$ is getting larger, the
rate of $\a(s)$ jumping from 2 to 1 is getting smaller.
Therefore, by assuming such kind of monotonicity for $q_{12}(\cd)$ and $q_{21}(\cd)$, one accepts that state 1 is ``bearish'' and state 2 is
``bullish'', if $X(\cd)$ is regarded as an index for the market.
If $X(\cd)$ is the total wealth process, then this also means that
the more the capital, the more the ``bullish'', which matches rational people's intuition. From the above, one sees the interesting feature
hidden in the model of state-dependent regime switching process.

\ms

Next, we turn to time-inconsistent issue. Although we could keep the state-dependent regime-switching, for the simplicity of presentation,
let us only consider the case that the matrix $(q_{ij})$ is
independent of the state $x$. For the state equation \rf{M2} with
the payoff functional \rf{P1}, the optimal control problem is time-consistent. For comparison purpose, let us try to solve such a
problem first (note that the problem is a maximization problem). For this
problem, the Hamiltonian is as follows:
$$\ba{ll}
\ns\ds\dbH(s,x,i,v,\Bp,\BP,u,c)=\Bp_i[b(i)u-c]+{1\over2}\si(i)^2\BP_iu^2
+q_{i1}v_1+q_{i2}v_2+g(s)c^\g\\
\ns\ds\qq\qq\qq\qq\q={\si(i)^2\BP_i\over2}\(u^2+2{b(i)\Bp_i\over\si(i)^2\BP_i}u\)
+g(s)\(c^\g-{\Bp_i\over g(s)}c\)+q_{i1}v_1+q_{i2}v_2.\ea$$
Thus,
$$\sup_{u\in\dbR,c\ges0}\dbH(s,x,i,v,\Bp,\BP,u,c)=-{b(i)^2\Bp_i^2\over
2\si(i)^2\BP_i}+{g(s)^{1\over1-\g}\g^{\g\over1-\g}(1-\g)
\over\Bp_i^{\g\over1-\g}}+q_{i1}v_1+q_{i2}v_2,$$
with the maximum attained at the following:
\bel{bar u,c}\bar u=-{b(i)\Bp_i\over\si(i)^2\BP_i},\qq\bar c=\({\g g(s)\over\Bp_i}\)^{1\over1-\g}.\ee
Consequently, the HJB equation reads:
$$\left\{\1n\ba{ll}
\ds(V_i)_s-{b(i)^2[(V^i)_x]^2\over2\si(i)^2(V_i)_{xx}}
+{g(s)^{1\over1-\g}\g^{\g\over1-\g}(1-\g)
\over[(V_i)_x]^{\g\over1-\g}}+q_{i1}V_1+q_{i2}V_2=0,\qq(s,x)\in[0,T)
\times(0,\infty),\\
\ns\ds V_i(T,x)=hx^\g,\qq x\in(0,\infty),\\
\ns\ds V_i(s,0)=0,\qq s\in[0,T],\qq\qq i=1,2.\ea\right.$$
Let
$$V_i(s,x)=\f_i(s)x^\g.$$
Then
$$(V_i)_s=\f'_i(s)x^\g,\qq(V_i)_x=\g\f_i(s)x^{\g-1},\qq(V_i)_{xx}=\g(\g-1)\f_i(s)x^{\g-2}.$$
Hence,
$$\ba{ll}
\ns\ds0=(V_i)_s-{b(i)^2[(V^i)_x]^2\over2\si(i)^2(V_i)_{xx}}
+{g(s){1\over1-g}\g^{\g\over1-\g}(1-\g)
\over[(V_i)_x]^{\g\over1-\g}}+q_{i1}V_1+q_{i2}V_2\\ [4mm]
\ns\ds=\f_i'(s)x^\g-{b(i)^2\f_i(s)^2\g^2x^{2(\g-1)}\over2\si(i)^2\g(\g-1)
\f_i(s)x^{\g-2}}+{g(s)^{1\over1-\g}\g^{\g\over1-\g}(1-\g)\over\f_i(s)^{\g\over1-\g}
\g^{\g\over1-\g}x^{(\g-1){\g\over1-\g}}}+q_{i1}\f_1(s)x^\g+q_{i2}\f_2(s)
x^\g\\ [4mm]
\ns\ds=\[\f_i'(s)-{b(i)^2\f_i(s)\g\over2\si(i)^2(\g-1)}
+{(1-\g)g(s)^{1\over1-\g}\over\f_i(s)^{\g\over1-\g}}+q_{i1}\f_1(s)+q_{i2}\f_2(s)\]x^\g.
\ea$$
Therefore, $\f_i(\cd)$ should be the solution to the following:
\bel{f}\left\{\1n\ba{ll}
\ds\f_i'(s)+{\g b(i)^2\over2(1-\g)\si(i)^2}\f_i(s)
+(1-\g)g(s)^{1\over1-\g}\f_i(s)^{\g\over\g-1}+q_{i1}\f_1(s)+q_{i2}\f_2(s)=0,
\q s\in[0,T),\\
\ns\ds\f_i(T)=h,\qq\qq i=1,2.\ea\right.\ee
Note that
$$V_i(t,x)\ges J(t,x;0,0)=hx^\g,\qq\forall(t,x)\in[0,T]\times[0,\infty).$$
Thus, one can show that
$$\f_i(s)\ges h>0,\qq s\in[0,T].$$
Consequently,
$$\f_i'(s)=-\[{\g b(i)^2\over2(1-\g)\si(i)^2}\f_i(s)
+(1-\g)g(s)^{1\over1-\g}\f_i(s)^{\g\over\g-1}+q_{i1}\f_1(s)+q_{i2}\f_2(s)\]
\les0,$$
and (note ${\g\over\g-1}<0$)
$$\ba{ll}
\ns\ds\f_i'(s)=-\[{\g b(i)^2\over2(1-\g)\si(i)^2}\f_i(s)
+(1-\g)g(s)^{1\over1-\g}\f_i(s)^{\g\over\g-1}+q_{i1}\f_1(s)+q_{i2}\f_2(s)\]\\
\ns\ds\qq~\ges-K\(\f_1(s)+\f_2(s)+h^{\g\over\g-1}\),\qq i=1,2.\ea$$
Hence, one can obtain the well-posedness of system \rf{f}. If we want, it will be easily to get the solution by, say, Picard iteration. Then
the optimal strategy can be obtained as follows:
\bel{bar u,c*}\ba{ll}
\ns\ds\bar u=-{b(i)V_x(s,x,i)\over\si(i)^2V_{xx}(s,x,i)}=-{b(i)\f_i(s)\g x^{\g-1}\over\si(i)^2\f_i(s)\g(\g-1)x^{\g-2}}={b(i)\over(1-\g)\si(i)^2}x,\\
\ns\ds\bar c=\({\g g(s)\over V_x(s,x,i)}\)^{1\over1-\g}=\({\g g(s)\over \f_i(s)\g x^{\g-1}}\)^{1\over1-\g}=\({g(s)\over\f_i(s)}\)^{1\over1-\g}x.\ea\ee
This solves the problem.

\ms

Now, let us look at the following modified payoff functional:
\bel{modified J}J(t,x;u(\cd),c(\cd))=\dbE_t\[\int_t^Tg(t,s)c(s)^\g ds+h(t)X(T)^\g\].\ee
Then for fixed $\t\in[0,T)$, consider the state equation \rf{M2} on
$[t,T]$, $t\in[\t,T)$ with payoff functional
\bel{modified J*}J(\t;t,x;u(\cd),c(\cd))=\dbE_\t\[\int_\t^Tg(\t,s)c(s)^\g ds+h(\t)X(T)^\g\].\ee
Applying the dynamic programming method as above, we will obtain the {\it pre-committed optimal strategy} as follows:
\bel{bar u,c**}\ba{ll}
\ns\ds\bar u(\t;s,x,i)={b(i)\over(1-\g)\si(i)^2}x,\qq
\bar c(\t;s,x,i)=\({g(\t,s)\over\f_i(t;s)}\)^{1\over1-\g}x,\ea\ee
with
\bel{f*}\left\{\2n\ba{ll}
\ds{d\over ds}\f_i(\t;s)+{\g b(i)^2\over2(1-\g)\si(i)^2}\f_i(\t;s)
+(1-\g)g(\t,s)^{1\over1-\g}\f_i(\t;s)^{\g\over\g-1}\\ [1mm]
\ns\ds\qq\qq\qq\qq\qq\qq\qq+q_{i1}\f_1(\t;s)+q_{i2}\f_2(\t;s)=0,\q s\in[\t,T),\\ [1mm]
\ns\ds\f_i(\t;T)=h(\t),\qq\qq i=1,2.\ea\right.\ee
Clearly, the optimal strategy given by \rf{bar u,c**} is
time-inconsistent since there are two time variables $\t$ (initial
time) and $s$ (running time) are involved.

\ms

According to the theory that we have developed, the Hamiltonian should look like the following:
$$\ba{ll}
\ns\ds\dbH(\t,s,x,i,v,\Bp,\BP,u,c)=\Bp_i[b(i)u-c]+{1\over2}
\si(i)^2\BP_iu^2+q_{i1}v_1+q_{i2}v_2+g(\t,s)c^\g\\
\ns\ds\qq\qq\qq\qq\qq\;={\si(i)^2\BP_i\over2}\(u^2+2{b(i)\Bp_i\over\si(i)^2\BP_i}u\)
+g(\t,s)\(c^\g-{\Bp_i\over g(\t,s)}c\)+q_{i1}v_1+q_{i2}v_2,\ea$$
and
$$\sup_{u\in\dbR,c\ges0}\dbH(\t,s,x,i,v,\Bp,\BP,u,c)
=-{b(i)^2\Bp_i^2\over2\si(i)^2\BP_i}+{g(\t,s)^{1\over1-\g}\g^{\g\over1-\g}(1-\g)
\over\Bp_i^{\g\over1-\g}}+q_{i1}v_1+q_{i2}v_2,$$
with the map $\psi\equiv(\bar u,\bar c)^\top$ appeared in (H3)
taking the following form:
\bel{psi=bar u,c}\bar u=-{b(i)\Bp_i\over\si(i)^2\BP_i},\qq
\bar c=\({\g g(\t,s)\over\Bp_i}\)^{1\over1-\g},\ee
and the map $\bar\Psi(s,x,i)$ appeared in \rf{Psi} is as follows:
\bel{bar u(s,s)}\bar u(s,x,i)=-{b(i)\Th_x(s,s,x,i)\over\si(i)^2
\Th_{xx}(s,s,x,i)},\qq\bar c(s,x,i)=\({\g g(s,s)\over
\Th_x(s,s,x,i)}\)^{1\over1-\g}.\ee
Consequently, the equilibrium HJB equation should read as
$$\left\{\2n\ba{ll}
\ns\ds\Th_s(\t,s,x,i)+{b(i)^2\Th_x(s,s,x,i)^2
\over2\si(i)^2\Th_{xx}(s,s,x,i)^2}\Th_{xx}(\t,s,x,i)\\
\ns\ds\qq-\Th_x(\t,s,x,i)\[{b(i)^2
\Th_x(s,s,x,i)\over\si(i)^2\Th_{xx}(s,s,x,i)}+\({\g g(s,s)\over
\Th_x(s,s,x,i)}\)^{1\over1-\g}\]\\
\ns\ds\qq+g(\t,s)\({\g g(s,s)\over
\Th_x(s,s,x,i)}\)^{\g\over1-\g}+q_{i1}\Th(\t,s,x,1)+q_{i2}
\Th(\t,s,x,2)=0,\qq(s,x)\in[\t,T)\times(0,\infty),\\
\ns\ds\Th_i(\t,T,x)=h(\t)x^\g,\qq x\in(0,\infty),\\
\ns\ds\Th_i(\t,s,0)=0,\qq \t\les s\les T,\qq\qq i=1,2.\ea\right.$$
Now, we introduce the following ansatz:
$$\Th(\t,s,x,i)=\f(\t,s,i)x^\g.$$
Then
$$\ba{ll}
\ns\ds0=\f_s(\t,s,i)x^\g+{b(i)^2\f(s,s,i)^2\g^2x^{2(\g-1)}
\over2\si(i)^2\g^2(\g-1)^2\f(s,s,i)^2x^{2(\g-2)}}\f(\t,s,i)\g(\g-1)
x^{\g-2}\\ [1mm]
\ns\ds\qq-\f(\t,s,i)\g x^{\g-1}\[{b(i)^2
\f(s,s,i)\g x^{\g-1}\over\si(i)^2\f(s,s,i)\g(\g-1)x^{\g-2}}
+\({\g g(s,s)\over\f(s,s,i)\g x^{\g-1}}\)^{1\over1-\g}\]\\ [1mm]
\ns\ds\qq+g(\t,s)\({\g g(s,s)\over
\f(s,s,i)\g x^{\g-1}}\)^{\g\over1-\g}+q_{i1}\f(\t,s,1)x^\g+q_{i2}
\f(\t,s,2)x^\g\\ [2mm]
\ns\ds\q=\Big\{\f_s(\t,s,i)+{b(i)^2\g\over2\si(i)^2(\g-1)}\f(\t,s,i)
-\f(\t,s,i)\[{b(i)^2\g\over\si(i)^2(\g-1)}
+\g\({g(s,s)\over\f(s,s,i)}\)^{1\over1-\g}\]\\ [4mm]
\ns\ds\qq+g(\t,s)\({g(s,s)\over
\f(s,s,i)}\)^{\g\over1-\g}+q_{i1}\f(\t,s,1)+q_{i2}
\f(\t,s,2)\Big\}x^\g\\ [4mm]
\ns\ds\q=\Big\{\f_s(\t,s,i)+\f(\t,s,i)\[{\g b(i)^2\over2(1-\g)
\si(i)^2}
-\g\({g(s,s)\over\f(s,s,i)}\)^{1\over1-\g}\]\\ [3mm]
\ns\ds\qq+g(\t,s)\({g(s,s)\over
\f(s,s,i)}\)^{\g\over1-\g}+q_{i1}\f(\t,s,1)+q_{i2}
\f(\t,s,2)\Big\}x^\g.\ea$$
Hence, $\f(s,x,i)$ should satisfy the following:
\bel{6.13}\left\{\1n\ba{ll}
\ds\f_s(\t,s,i)+{\g b(i)^2\over2(1-\g)\si(i)^2}\f(\t,s,i)
-\g\f(\t,s,i)\({g(s,s)\over\f(s,s,i)}\)^{1\over1-\g}\\
\ns\ds\qq+g(\t,s)\({g(s,s)\over
\f(s,s,i)}\)^{\g\over1-\g}+q_{i1}\f(\t,s,1)+q_{i2}
\f(\t,s,2)=0,\qq0\les\t\les s\les T,\\
\ns\ds\f(\t,T,i)=h(\t),\qq0\les\t\les T,\qq\qq i=1,2.\ea\right.\ee
We point out that when $g(\t,s)\equiv g(s,s,)\equiv g(s)$, the above
is reduced to \rf{f}. Having solved \rf{6.13}, the equilibrium strategy is given by
\bel{bar u(s,s)*}\bar u(s,x,i)=-{b(i)x\over(1-\g)\si(i)^2},\qq\bar c(s,x,i)=\({g(s,s)\over\f(s,s,i)}\)^{1\over1-\g}x.\ee
One can further discuss the well-posedness of \rf{6.13}, which we
prefer not to get into more details.

\section{Concluding Remarks}

We have explored the time-inconsistent stochastic optimal control problem for the regime-switching SDEs, with state-dependent transition probability matrix for the regime-switching process. Due to such a state-dependence, the theory for the state equation has some unique interesting features, such as the continuous dependence of the solution on the initial state. Therefore, our results enrich the general theory of regime-switching diffusion processes. For our time-inconsistent optimal control problem, Equilibrium HJB equation has been successfully derived, and for an important special case (i.e., the diffusion of the state equation is independent of the control), the well-posedness of the equilibrium HJB equation has been established. Consequently, at least for this case, a time-consistent equilibrium strategy can be constructed in principal. An example is also presented to illustrate our main results.

\ms

We have seen that the theory is still in its infancy. Many questions are left open. Here are a partial list of the problems that we will continue to investigate. Of course, we welcome other interested researchers to join:

\ms

$\bullet$ The well-posedness of the equilibrium HJB equation for the case that the control appears in the diffusion of the state equation.

\ms

$\bullet$ The case that the diffusion in the state equation is possibly degenerate.

\ms

$\bullet$ Possible efficient numerical algorithm aspects.

\ms


\begin{thebibliography}{}


\bibitem{Bjork-Khapko-Murgoci2017} T.~Bj\"ork, M.~Khapko, and A.~Murgoci, \it On time-inconsistent stochastic control in continuous time, \sl Finance
    Stoch., \rm 21 (2017), 331--360.

\bibitem{Bjork-Murgoci2014} T.~Bj\"{o}rk and A.~Murgoci, \it A theory of Markovian time-inconsistent stochastic control in discrete time, \sl Finance Stoch, \rm 18 (2014), 545--592.

\bibitem{Bjork-Murgoci-Zhou2014} T.~Bj\"{o}rk, A.~Murgoci, and X.~Y.~Zhou, \it Mean-variance portfolio optimization with state-dependent risk aversion, \sl Math. Finance, \rm 24 (2014), 1--24.

\bibitem{Cohen-Elliott2008} S.~N.~Cohen and R.~J.~Elliott, \it Solutions of backward stochastic differential equations on Markov chains, \it Communications on Stochastic Analysis, \rm 2 (2008), 251--262.
\bibitem{Donnelly2011} C.~Donnelly, \it Sufficient stochastic maximum principle in a regime-switching diffusion model. \sl Appl. Math Optim., \rm 64 (2011), 155--169.

\bibitem{Donnelly-Heunis2012} C.~Donnelly and A.~J.~Heunis, \it Quadratic risk minimuzation in a regime-swtching model with portfolio constraints, \sl SIAM J. Control Optim., \rm 50 (2012), 2431--2461.

\bibitem{Duffie-Epstein1992a} D.~Duffie and L.~G.~Epstein, \it Stochastic differential utility, \sl Econometrica, \rm 60 (1992), 353--394.

\bibitem{Duffie-Epstein1992b} D.~Duffie and L.~G.~Epstein, \it Asset pricing with stochastic differential utility, \sl Review Financial Studies, \rm 5 (1992), 411--436.

\bibitem{Duffie-Lions1992} D.~Duffie and P.~L.~Lions, \it PDE solutions of stochastic differential utility, \sl J. Math. Econom., \rm 21 (1992), 577--606.

\bibitem{EI1969} S.~D.~Eidel'man SD. \sl Parabolic Systems, \rm North Holland Publishing Company, 1969.

\bibitem{Ekeland-Lazrak2010} I.~Ekeland and A.~Lazrak, \it The golden rule when preferences are time inconsistent, \sl Math. Financ. Econ., \rm 4 (2010), 29--55.

\bibitem{El Karoui-Peng-Quenez1997} N.~El Karoui, S.~Peng, and M.~C.~Quenez, \it Backward stochastic differential equations in finance, \sl Math. Finance, \rm  7 (1997), 1--71.

\bibitem{Framstad-Oksenda-Sulem2004} N.~C.~Framstad, B.~K.~$/\3n\3n\;{\rm O}$ksendal and A.~Sulem,
\it Sufficient stochastic maximum principle for the optimal control of jump diffusions and
applications to finance, \sl J. Optim. Theory Appl., \rm 121 (2004), 77--98.

\bibitem{Friedman1964} A.~Friedman, \sl Partial Differential Equations of Parabolic Type, \rm Prentice Hall, Inc., Englewood Cliffs, NJ, 1964.

\bibitem{Hamilton1989} J.~D.~Hamilton, \it A new approach to the economic analysis of nonstationry time series and the business cycle, \sl Econometrica, \rm 57 (1989), 357--384.

\bibitem{Heunis2015} A.~J.~Heunis, \it Utility maximization in a regime switching model wih convex portfolio constratints and margin requirements: Optimality relations and explicit solutions, \sl SIAM J. Control Optim., \rm 53 (2015), 2608--2656.

\bibitem{Hu-Jin-Zhou2012} Y.~Hu, H.~Jin, and X.~Y.~Zhou, \it Time-inconsistent stochastic linear-quadratic control, \sl SIAM J. Control Optim., \rm 50 (2012), 1548--1572.

\bibitem{Hume1739} D.~Hume, \sl A Treatise of Human Nature, \rm First edition, 1739; Reprint, Oxford Univ. Press, New York, 1978.

\bibitem{Karnam-Ma-Zhang2016} C.~Karnam, J.~Ma, and J.~Zhang, \it Dynamic approaches for some time inconsistent problems, \rm preprint.

\bibitem{Kraft-Seifried2014} H.~Kraft and F.~T.~Seifried, \it Stochastic differential utility as the continuous-time limit of recursive utility, \sl Journal of Economic Theory, \rm 151 (2014), 528--550.

\bibitem{Kruse-Popier2016} T.~Kruse and A.~Popier, \it BSDEs with monotone generator driven by Brownian and Poisson noises in a general filtration. \sl Stochastics, \rm 88(2016), 491--539.

\bibitem{Kruse-Popier2017} T.~Kruse and A.~Popier, \it $L^p$-solution for BSDEs with jumps in the case $p<2$. Corrections to the paper ``BSDEs with monotone generator driven by Brownian and Poisson noises in a general filtration'', \rm arXiv:1701.09071v1[math.PR] 31 Jan 2017.

\bibitem{Laibson1997} D.~Laibson, \it Golden eggs and hyperbolic discounting, \sl Qartely J. Econ., \rm 112 (1997), 443--477.

\bibitem{Lazrak2004} A.~Lazrak, \it Generalized stochastic differential utility and preference for information, \sl Ann. Appl. Probab., \rm  14 (2004), 2149--2175.

\bibitem{Lazrak-Quenez2003} A.~Lazrak and M.~C.~Quenez, \it A generalized stochastic differential utility, \sl Math. Oper. Res., \rm 28 (2003), 154--180.

\bibitem{Ma-Protter-Yong1994} J.~Ma, P.~Protter, and J.~Yong, \it Solving forward-backward stochastic differential equations explicitly - a four step scheme., \sl Probab. Theory Related Fields, \rm  98  (1994), 339--359.

\bibitem{Ma-Wu-Zhang-Zhang2015} J.~Ma, Z.~Wu, D.~Zhang, and J.~Zhang, \it On well-posedness of forward-backward SDEs --- a unified approach, \sl Ann. Appl. Probab., \rm 25 (2015), 2168--2214.

\bibitem{Ma-Yong1999} J.~Ma and J.~Yong, \sl Forward-Backward Stochastic Differential Equations and Their Applications, \rm Lecture Notes in Math., Vol. 1702, Springer-Verlag, 1999.

\bibitem{MN2010} J.~Marin-Solano and J.~Navas, \it Consumption and portfolio rules for time-inconsistent investors, \sl European J. Oper. Res., \rm 201 (2010), 860--872.

\bibitem{MS2011} J.~Marin-Solano and E.~V.~Shevkoplyas, \it Non-constant discounting and differential games with random time horizon, \sl Automatica, \rm 47 (2011), 2626--2638.

\bibitem{Oksendal-Sulem2005} B.~K.~$/\3n\2n\;{\rm O}$ksendal and A.~Sulem, \sl Applied Stochastic Control of Jump Diffusions, \rm Springer-Verlag, Berlin, 2005.

\bibitem{Palacios-Huerta2003} I.~Palacios-Huerta, \it Time-inconsistent preferences in Adam Smith and David Hume, \sl History of Political Ecomony, \rm 35 (2003), 241--268.




\bibitem{Pardoux-Peng1990} E.~Pardoux and S.~Peng, \it Adapted solution of backward stochastic equation, \sl Syst. \& Control Lett., \rm 14 (1990),
    55--61.

\bibitem{Pollak1968} R.~A.~Pollak, \it Consistent planning, \sl Review of Economic Studies, \rm 35 (1968), 185--199.

\bibitem{Situ2006} R.~Situ, \sl Theory of Stochastic Differential Equations with Jumps and Applications: Mathematical and Analytical Techniques with Applications to Engineering, \rm Springer Science \& Business Media, 2006.

\bibitem{Smith1759} A.~Smith, \sl The Theory of Moral Senriments, \rm First edition, 1759; Reprint, Oxford Univ. Press, 1976.

\bibitem{Sotomayor-Cadenillas2009} L.~R.~Sotomayor and A.~Cadenilla, \it Explicit solutions of consumption-investment problems in financecial markets with regime switching, \sl Math Finance, \rm 19 (2009), 251--279.

\bibitem{Strotz1955} R.~H.~Strotz, \it Myopia and inconsistency in dynamic utility maximization, \sl Review of Economic Studies, \rm 23 (1955), 165--180.

\bibitem{Wei2017} J.~Wei, \it Time-inconsistent optimal control problems with regime-switching, \sl Math. Control Related Fields, \rm 7 (2017), 585--622.?

\bibitem{Wei-Yong-Yu2016} Q.~Wei, J.~Yong, and Z.~Yu, \it Time-inconsistent recrusive stochastic optimal control problems, \sl SIAM J. Control Optim., \rm 55 (2017), 4156--4201.


\bibitem{Yin-Zhu2009} G.~Yin and C.~Zhu, \sl Hybrid Switching Diffusions: Properties and Applications, \rm Vol. 63, Springer Science \&  Business Media, 2009.



\bibitem{Yong2012b} J. Yong, \it Time-inconsistent optimal control problems and the equilibrium HJB equation, \sl Math. Control Relat. Fields, \rm 2 (2012), no. 3, 271-329.


\bibitem{Yong2014} J. Yong, \it Time-inconsistent optimal control problems, \sl Proceedings of 2014 ICM,
    Section 16. Control Theory and Optimization, \rm 947--969.

\bibitem{Yong2015} J.~Yong, \it Linear-quadratic optimal control problems for mean-field stochastic differential equations --- time-consistent solutions, \sl Trans. AMS, \rm 369 (2017), 5467--5523.

\bibitem{Yong-Zhou1999} J.~Yong and X.~Y.~Zhou, \sl Stochastic Control: Hamiltonian Systems and HJB Equations, \rm Springer-Verlag, New York, 1999.

\bibitem{Zhang2017} J.~Zhang, \sl Backward Stochastic Differential Equations --- From Linear to Fully Nonlinear Theory, \rm Springer, 2017.

\bibitem{Zhang-Elliott-Siu2012} X.~Zhang, R.~J.~Elliott, and T.~K.~Siu, \it A stochastic maximum principle for a Markov regime-swtching jump-diffusion model and its application to finance, \sl SIAM J. Control Optim., \rm 50 (2012), 964--990.

\bibitem{Zhou-Yin2003} X.~Y.~Zhou and G.~Yin, \it Markowitz's mean-varaince portfolio selection with regime swtching: A contionus-time model, \sl SIAM J. Control Optim., \rm 42 (2003), 1466--1482.


\end{thebibliography}
\end{document}